\documentclass[]{scrartcl}

\usepackage[utf8]{luainputenc}
\usepackage[USenglish]{babel}
\usepackage{csquotes}

\usepackage[a4paper,top=27mm,bottom=20mm,inner=25mm,outer=20mm]{geometry}

\usepackage[%
  backend=bibtex,bibencoding=ascii,
  style=numeric-comp,
  giveninits=true, uniquename=init, 
  natbib=true,
  url=true,
  doi=true,
  isbn=false,
  backref=false,
  maxnames=99,
  ]{biblatex}
\usepackage{amsmath}
\allowdisplaybreaks
\numberwithin{equation}{section}
\usepackage{amssymb}
\usepackage{commath}
\usepackage{mathtools}
\usepackage{bbm}
\usepackage{nicefrac}
\usepackage{subdepth}

\usepackage{siunitx}
\usepackage{algorithm}
\usepackage{algpseudocode}

\usepackage{amsthm}
\usepackage{thmtools}
\usepackage{etoolbox}
\makeatletter
\patchcmd{\thmt@setheadstyle}
 {\bgroup\thmt@space}
 {\thmt@space}
 {}{}
\patchcmd{\thmt@setheadstyle}
 {\egroup\fi}
 {\fi}
 {}{}
\makeatother
\declaretheoremstyle[
  bodyfont=\normalfont\itshape,
  headformat=\NAME\ \NUMBER\NOTE,
]{myplain}
\declaretheoremstyle[
  headformat=\NAME\ \NUMBER\NOTE,
]{mydefinition}
\newcommand{\envqed}{{\lower-0.3ex\hbox{$\triangleleft$}}}

\usepackage[plainpages=false,pdfpagelabels,hidelinks,unicode]{hyperref}

\usepackage{color}
\usepackage{graphicx}
\usepackage[small]{caption}
\usepackage{subcaption}

\begingroup\expandafter\expandafter\expandafter\endgroup
\expandafter\ifx\csname pdfsuppresswarningpagegroup\endcsname\relax
\else
\fi

\usepackage{booktabs}
\usepackage{rotating}
\usepackage{multirow}

\usepackage{enumitem}

\usepackage{ifluatex}
\ifluatex
  \usepackage[no-math]{fontspec}
\else
  \usepackage[T1]{fontenc}
\fi
\usepackage{newpxtext,newpxmath}

\usepackage{xparse}

\let\epsilon\varepsilon
\let\phi\varphi
\let\rho\varrho

\usepackage{xspace}

\newcommand{\eg}[0]{{e.g.\@}\xspace}
\newcommand{\ie}[0]{{i.e.\@}\xspace}
\newcommand{\ca}[0]{{ca.\@}\xspace}

\providecommand\R{}
\renewcommand{\R}{\mathbb{R}}

\newcommand{\bhat}{\widehat{b}}
\newcommand{\uhat}{\widehat{u}}
\newcommand{\qhat}{\widehat{q}}
\newcommand{\tol}{\ensuremath{\texttt{tol}}}
\newcommand{\atol}{\ensuremath{\tau_a}}
\newcommand{\rtol}{\ensuremath{\tau_r}}

\newcommand{\cfl}{\nu}
\newcommand{\maxspeed}{\lambda_\mathrm{max}}

\NewDocumentCommand{\RK}{o m O{\the\numexpr#2-1\relax} m O{} O{} o}{%
  \IfValueTF{#1}{#1}{RK}%
  #2(#3)#4%
  \ifblank{#6}{}{\textsubscript{F}}%
  \ifblank{#5}{}{[#5]}%
  \IfValueT{#7}{#7}%
}

\newcommand{\ESstarp}{3S*\textsubscript{+}}
\newcommand{\ssp}[2]{\RK[SSP]{#2}{#1}[\ESstarp]}

\newcommand{\dt}{\Delta t}
\newcommand{\dx}{\Delta x}

\newcommand{\ndims}{\ensuremath{d}}

\newcommand{\polydeg}{\ensuremath{p}}

\renewcommand{\vec}[1]{\pmb{#1}}

\NewDocumentCommand{\opD}{m+g}{%
  \IfNoValueTF{#2}
    {D_{#1}}
    {D_{#1,#2}}%
}
\NewDocumentCommand{\opDsplit}{m+g}{%
  \IfNoValueTF{#2}
    {\widetilde{D}_{#1}}
    {\widetilde{D}_{#1,#2}}%
}
\NewDocumentCommand{\opM}{g}{%
  \IfNoValueTF{#1}
    {M}
    {M_{#1}}%
}
\NewDocumentCommand{\opQ}{g}{%
  \IfNoValueTF{#1}
    {Q}
    {Q_{#1}}%
}
\NewDocumentCommand{\opI}{g}{%
  \IfNoValueTF{#1}
    {I}
    {I_{#1}}%
}
\NewDocumentCommand{\opV}{g}{%
  \IfNoValueTF{#1}
    {V}
    {V_{#1}}%
}
\NewDocumentCommand{\opB}{g}{%
  \IfNoValueTF{#1}
    {B}
    {B_{#1}}%
}
\NewDocumentCommand{\opR}{g}{%
  \IfNoValueTF{#1}
    {R}
    {R_{#1}}%
}
\NewDocumentCommand{\opN}{m+g}{%
  \IfNoValueTF{#2}
    {N_{#1}}
    {N_{#1,#2}}%
}

\NewDocumentCommand{\fnum}{g}{%
  \IfNoValueTF{#1}
    {f^{\mathrm{num}}}
    {f^{\mathrm{num,#1}}}%
}
\NewDocumentCommand{\vecfnum}{g}{%
  \IfNoValueTF{#1}
    {\vec{f}^{\mathrm{num}}}
    {\vec{f}^{\mathrm{num,#1}}}%
}
\NewDocumentCommand{\vecfcorr}{g}{%
  \IfNoValueTF{#1}
    {\vec{f}^{\mathrm{corr}}}
    {\vec{f}^{\mathrm{corr,#1}}}%
}
\NewDocumentCommand{\fvol}{g}{%
  \IfNoValueTF{#1}
    {f^{\smash{\mathrm{vol}}}}
    {f^{\smash{\mathrm{vol,#1}}}}%
}

\newcommand{\orcid}[1]{ORCID:~\href{https://orcid.org/#1}{#1}}
\usepackage{authblk}

\newenvironment{keywords}{\par\textbf{Key words.}}{\par}
\newenvironment{AMS}{\par\textbf{AMS subject classification.}}{\par}

\title{On error-based step size control for discontinuous Galerkin methods for compressible fluid dynamics}

\author[1]{Hendrik~Ranocha\thanks{\orcid{0000-0002-3456-2277}}}
\affil[1]{Applied Mathematics, University of Hamburg, Germany}

\author[2]{Andrew R.~Winters\thanks{\orcid{0000-0002-5902-1522}}}
\affil[2]{Department of Mathematics; Applied Mathematics, Linköping University, Sweden}

\author[3]{Hugo Guillermo Castro\thanks{\orcid{0000-0003-1715-1238}}}
\author[3]{Lisandro Dalcin\thanks{\orcid{0000-0001-8086-0155}}}
\affil[3]{King Abdullah University of Science and Technology (KAUST), Extreme Computing Research Center (ECRC), Computer Electrical and Mathematical Science and Engineering Division (CEMSE), Thuwal, 23955-6900, Saudi Arabia}

\author[4,5]{Michael Schlottke-Lakemper\thanks{\orcid{0000-0002-3195-2536}}}
\affil[4]{Applied and Computational Mathematics, RWTH Aachen University, Germany}
\affil[5]{High-Performance Computing Center Stuttgart, University of Stuttgart, Germany}

\author[6,7]{Gregor J. Gassner\thanks{\orcid{0000-0002-1752-1158}}}
\affil[6]{Department of Mathematics and Computer Science, University of Cologne, Germany}
\affil[7]{Center for Data and Simulation Science, University of Cologne, Germany}

\author[3]{Matteo Parsani\thanks{\orcid{0000-0001-7300-1280}}}

\date{February 5, 2023} 

\makeatletter
\hypersetup{pdfauthor={Hendrik Ranocha et al.}} 
\hypersetup{pdftitle={On error-based step size control for discontinuous Galerkin methods for compressible fluid dynamics}} 
\makeatother

\begin{document}

\maketitle

\begin{abstract}
\noindent
  We study temporal step size control of explicit Runge-Kutta methods for
compressible computational fluid dynamics (CFD), including the Navier-Stokes
equations and hyperbolic systems of conservation laws such as the Euler equations.
We demonstrate that error-based approaches
are convenient in a wide range of applications and compare them to more classical
step size control based on a Courant-Friedrichs-Lewy (CFL) number. Our numerical
examples show that error-based step size control is easy to use, robust, and efficient,
e.g., for (initial) transient periods, complex geometries, nonlinear shock
capturing approaches, and schemes that use nonlinear entropy projections.
We demonstrate these properties for problems ranging from well-understood
academic test cases to industrially relevant large-scale computations with two
disjoint code bases, the open source Julia packages Trixi.jl with OrdinaryDiffEq.jl
and the C/Fortran code SSDC based on PETSc.

\end{abstract}



\section{Introduction}
\label{sec:introduction}

Hyperbolic conservation laws are at the heart of many physical models in science
and engineering, e.g., fluid dynamics describing the airflow around an airplane,
acoustics, and space plasma modeling. Usually, the systems need to be solved
numerically. A common approach is to apply the method of lines, i.e., to apply
a spatial semidiscretization first and solve the resulting ordinary differential
equation (ODE) numerically. For this task, explicit Runge-Kutta (RK) methods
are among the most commonly used schemes because of their efficiency and parallel
scalability \cite{ka.sh:05,kolev2021efficient,kopriva2013assessment}.

Since demanding simulations require a significant amount of compute resources,
there is a strong interest in developing efficient numerical methods. Here,
efficiency can be measured by the required compute time, energy, and human time.
The time- and energy-to-solution depend on the algorithms and their
implementation as well as the hardware \cite{vincent2016towards,becker2017exaflop,heinisch2020towards}
and are often a focus of high performance computing (HPC) studies. The
required human time depends on the robustness of the overall numerical method
and the sensitivity with respect to parameters that may need to be tuned to
ensure robustness and/or efficiency. In this article, we study error-based step
size control for compressible flow problems to demonstrate that it is efficient
in all three aspects.

We use entropy-dissipative semidiscretizations to ensure robustness
\cite{fisher2013high,gassner2016split,sjogreen2018high,rojas2021robustness}
but do not investigate specific implementation techniques discussed elsewhere
\cite{maier2021efficient,guermond2021implementation,ranocha2021efficient}.
Instead, we focus on step size control of time integration methods applied to
ODEs resulting from spatial semidiscretizations of compressible flow problems.
The goal is to adapt the time step size $\dt$ such that it is as big as possible
while still satisfying accuracy and stability requirements. For explicit RK methods
applied to (high-order) semidiscretizations, stability requirements are usually
more restrictive. Thus, a commonly used approach is to adapt the time step size
according to a Courant-Friedrichs-Lewy (CFL) number chosen by the user. This involves
estimating local wave speeds and mesh spacing, which can be demanding for complex
systems and curved high-order meshes in multiple space dimensions required in
practice. The optimal CFL number depends on the space and time discretizations
and possibly on the specific problem
\cite{saad2022stable,kubatko2008time,chalmers2020robust}.
An alternative, widely used in the context of general-purpose ODE solvers, is
an error-based step size control. There, a local error estimate is obtained from
an embedded method and fed into a controller adapting the time step size according
to user-specified tolerances. This approach has also been used successfully for
partial differential equations
\cite{berzins1995temporal,ware1995adaptive,ketcheson2020more,ranocha2021optimized}
but appears to be less-widely applied in the CFD community.

In this article, we systematically study error-based step size control
for compressible CFD problems and compare it to CFL-based approaches.
Specifically, we use the RK methods and controllers developed in
\cite{ranocha2021optimized} in two different software environments.
Most of the examples use the open source packages OrdinaryDiffEq.jl
\cite{rackauckas2017differentialequations} for ODE solvers and Trixi.jl
\cite{ranocha2022adaptive,schlottkelakemper2021purely} for spatial
semidiscretizations implemented in Julia \cite{bezanson2017julia}.
Some large-scale, industrially-relevant CFD simulations are implemented in
SSDC \cite{parsani2021ssdc}, which is built on top of the
Portable and Extensible Toolkit for Scientific computing (PETSc) \cite{petsc313},
its mesh topology abstraction (DMPLEX) \cite{knepley2009mesh}, and its scalable
ODE/DAE solver library \cite{abhyankar2018petscts}.
All source code and instructions required to reproduce the numerical experiments
using open source packages is available online in our reproducibility repository
\cite{ranocha2022errorRepro}.

In the following, we briefly review RK methods and step size control techniques in
Section~\ref{sec:rk}. Afterwards, we study the robustness and sensitivity to
user-supplied parameters under a change of mesh structure of both step size
control approaches in Section~\ref{sec:meshes}. Section~\ref{sec:dg_fv} focuses
on the effect of nonlinear shock capturing schemes. In Section~\ref{sec:cold_start},
we discuss the impact of an initial transient period, e.g., in the cold start
of a simulation initialized with free stream values. Thereafter, we study
step size control in the presence of a change of variables in
Section~\ref{sec:GaussSBP}. Next, we discuss the convenience of error-based
step size control for exploratory research of new systems in
Section~\ref{sec:new_systems}. Finally, we summarize and discuss our findings
in Section~\ref{sec:summary}.

\section{Runge-Kutta methods and step size control}
\label{sec:rk}

We use the method of lines approach and first discretize a conservation
law in space. In particular, we focus on spatial semidiscretizations
using collocated discontinuous Galerkin methods (and related shock capturing schemes), see, e.g.,
\cite{fisher2013high,gassner2016split,sjogreen2018high,rojas2021robustness}.
The resulting ordinary differential equation
\begin{equation}
\label{eq:ode}
\begin{aligned}
  \od{}{t} u(t) &= f(t, u(t)),
  && t \in (0,T),
  \\
  u(0) &= u^{0},
\end{aligned}
\end{equation}
is then solved using numerical time integration methods. Here, we use explicit
Runge-Kutta methods with an embedded error estimator of the form
\cite{hairer2008solving,butcher2016numerical}
\begin{equation}
\label{eq:RK-step}
\begin{aligned}
  y_i
  &=
  u^n + \dt_n \sum_{j=1}^{i-1} a_{ij} \, f(t_n + c_j \dt_n, y_j),
  \qquad i \in \{1, \dots, s\},
  \\
  u^{n+1}
  &=
  u^n + \dt_n \sum_{i=1}^{s} b_{i} \, f(t_n + c_i \dt_n, y_i),
  \\
  \uhat^{n+1}
  &=
  u^n + \dt_n \sum_{i=1}^{s} \bhat_{i} \, f(t_n + c_i \dt_n, y_i) + \bhat_{s+1} f(t_{n+1}, u^{n+1}).
\end{aligned}
\end{equation}
Here, $u^{n+1}$ is the numerical solution of the main method, $\uhat^{n+1}$ is
the embedded methods solution used
to obtain the local error estimate $u^{n+1} - \uhat^{n+1}$,
and $y_i$ are the stage values.
RK methods are parameterized by their Butcher tableau
\begin{equation}
\label{eq:butcher}
\begin{array}{c | c}
  c & A
  \\ \hline
    & b^T
  \\
    & \bhat^T
\end{array}
\end{equation}
where $A = (a_{ij})_{i,j=1}^s \in \R^{s \times s}$ is strictly lower-triangular
for explicit methods, $b, c \in \R^s$, and $\bhat \in \R^{s+1}$.
If $\bhat_{s+1} \ne 0$, the RK pair \eqref{eq:RK-step} uses the so-called
first-same-as-last (FSAL) idea to include the predicted value at the new time
in the embedded error estimator \cite{dormand1980family}. This can increase
the performance of the embedded method and comes at no additional cost if the
step is accepted since $f(t_{n+1}, u^{n+1})$ must be computed as the first stage
of the next step anyway.

All methods considered in this article are applied in local extrapolation mode,
i.e., a main method of order $q$ is coupled with an embedded error estimator
of order $\qhat = q - 1$. Whenever possible, we make use of low-storage formats
such as \ESstarp\ \cite{ketcheson2010runge,ranocha2021optimized}. We use the
same naming convention of RK methods as in \cite{ranocha2021optimized}, i.e.,
\RK[NAME]{$q$}[$\qhat$]{$s$} indicates an $s$-stage method of order $q$ with
embedded error estimator of order $\qhat$. This base name is followed by
additional information on low-storage properties and a subscript ``F'' for
FSAL methods.

\subsection{Time integration loop and callbacks}

We use implementations of the time integration methods available in
OrdinaryDiffEq.jl \cite{rackauckas2017differentialequations} in Julia
or PETSc/TS \cite{abhyankar2018petscts}.
The basic steps for the time integration loop are given in Algorithm~\ref{alg:time_loop_crude}.
At first, a preliminary time step is performed with the given RK method.
If error-based step size control is activated, the embedded error estimator
is also computed; it is used to update the time step size and determine whether
the time step is rejected. After accepting a time step, callbacks are activated.
We use these callback mechanisms for adaptive mesh refinement (AMR) and, if activated, CFL-based step
size control.

\begin{algorithm}[htb]
  \caption{High-level overview of a time integration step in the Julia package
           OrdinaryDiffEq.jl.}
  \label{alg:time_loop_crude}
  \begin{enumerate}
    \item Compute preliminary time step update $u^{n+1}$ and the
          embedded solution $\uhat^{n+1}$ if error-based step size
          control is activated using \eqref{eq:RK-step}
    \label{itm:step}
    \item If error-based step size control is activated:
          \begin{itemize}
            \item Let the step size controller update the time step size
            \item If the time step is rejected, repeat step \ref{itm:step} with
                  the updated time step size
          \end{itemize}
    \item Accept the new time step and update internal caches
    \item Activate callbacks
    \item Proceed to the next time step (i.e., increment the time step number
          $n$ and go to step \ref{itm:step})
  \end{enumerate}
\end{algorithm}

A more detailed overview of the time integration loop is given below in
Algorithm~\ref{alg:time_loop_more_detailed}, discussed after introducing the
step size control mechanisms.

\subsection{CFL-based step size control}

Explicit time integration methods for first-order hyperbolic conservation laws
are subject to a CFL time step restriction of the
form $\dt \leq C \dx / \maxspeed$ \cite{courant1928partiellen}. However, it is non-trivial
to provide sharp estimates of the terms in this restriction, e.g., an appropriate value for $\Delta x$ on curved meshes for
high-order discontinuous Galerkin methods in multiple space dimensions.
In general, CFL-based step size control requires user estimates of the maximal
local wave speeds $\maxspeed(u^n_i)$ at the degree of freedom $u^n_i$ and the
corresponding local mesh spacing $\dx_i$. Then, the time step is determined as
\begin{equation}
\label{eq:cfl-dt}
  \dt_n = \cfl \, \min_i \frac{\dx_i}{\maxspeed(u^n_i)},
\end{equation}
where $\cfl$ is the desired CFL number. For discontinuous Galerkin methods
using polynomials of degree $\polydeg$ on (possibly) curved grids for the
linear advection equation $\partial_t u + (a \cdot \partial_x) u = 0$,
we follow
\cite{jahdali2021optimized,ranocha2021optimized}
and use the estimate
\begin{equation}
\label{eq:cfl-factor}
  \frac{\dx_i}{\maxspeed(u^n_i)}
  =
  \frac{2}{\polydeg + 1}
  \frac{J_i}{\sum_{j=1}^\ndims \vert (J \partial_{x} \xi^j)_i \cdot a \vert},
\end{equation}
where $\ndims$ is the number of spatial dimensions,
$J_i$ the determinant of the grid Jacobian $\partial_{x} \xi$ at node $i$,
$(J \partial_{x} \xi^j)_i$ the contravariant basis vector in direction $j$
at node $i$ \cite[Chapter~6]{kopriva2009implementing}, and $2$ is the size of the reference element. For nonlinear conservation
laws, $a$ is replaced by appropriate estimates of the local wave speeds.

\subsection{Error-based step size control}

We perform error-based step size control based on controller ideas from
digital signal processing
\cite{gustafsson1988pi,gustafsson1991control,soderlind2002automatic,soderlind2003digital,soderlind2006time}.
In particular, we use PID controllers that select a new time step based on
\begin{equation}
\label{eq:PID}
  \dt_{n+1} = \kappa\Bigl(%
                  \epsilon_{n+1}^{\beta_1 / k}
                  \epsilon_{n }^{\beta_2 / k}
                  \epsilon_{n-1}^{\beta_3 / k} \Bigr) \dt_{n},
\end{equation}
where $k = \min(q, \qhat) + 1$, $q$ is the order of the main method, and $\qhat$ is the order of the
embedded method. Because $\qhat = q  - 1$ in our case, we have that $k = \qhat + 1 = q$.
Moreover, $\beta_i$ are the controller
parameters, $\kappa$ is a step size limiter, and
\begin{equation}
  \epsilon_{n+1} = \frac{1}{w_{n+1}},
  \qquad
  w_{n+1} = \left( \frac{1}{m} \sum_{i=1}^{m} \left( \frac{u_i^{n+1} - \uhat_i^{n+1}}{\atol + \rtol \max\{ \vert u_i^{n+1}\vert, \vert\widetilde{u}_i\vert \}} \right)^2 \right)^{1/2},
\end{equation}
where $m$ is the total number of degrees of freedom in $u$, and $\atol, \rtol$
are the absolute and relative error tolerances. In OrdinaryDiffEq.jl,
$\widetilde{u} = u^{n}$ while PETSc uses $\widetilde{u} = \uhat^{n+1}$.
Unless stated otherwise, we use equal tolerances $\atol = \rtol = \tol$ and the
default step size limiter $\kappa(a) = 1 + \arctan(a - 1)$ \cite{soderlind2006adaptive}.

The decision whether the step shall be accepted or rejected is determined by the
size of the factor multiplying the time step size $\dt_{n}$ in \eqref{eq:PID}.
The default option used for all numerical experiments is to accept a step
whenever the limited factor is at least $0.9^2$. Otherwise, the step is
rejected and the time step size $\dt_{n}$ is set as the new predicted step size
\eqref{eq:PID}.
Another common strategy is to decide whether a step should be rejected
based only on the current local error estimate. Söderlind et al.
\cite{soderlind2003digital,soderlind2006adaptive} argue why the approach we
described here can be beneficial to reduce the amount of (unnecessary but
expensive) step rejections.

To make the step size control fully automatic, we use the estimate of the
initial step size implemented in OrdinaryDiffEq.jl from the algorithm
described by Hairer et al.~\cite[p.~169]{hairer2008solving}.
$\epsilon_{n}$ and $\epsilon_{n-1}$ are initially set to one (for $n = 0$).

A more detailed overview of the time integration loop including additional
aspects about step size control and step rejections is shown in
Algorithm~\ref{alg:time_loop_more_detailed}.
Parameters such as the threshold $0.9^2$ used to determine whether a step
should be accepted or the limiting threshold $w_\mathrm{min}$ used in
Algorithm~\ref{alg:time_loop_more_detailed} are based on recommendations in
reference works such as \cite[Section~IV.2]{hairer2008solving} and
\cite[Section~371]{butcher2016numerical}; we always use their default values
and do not consider them as user-facing parameters that need to be chosen
manually.

\begin{algorithm}[htbp]
  \caption{High-level overview of a time integration loop using explicit
           Runge-Kutta methods with embedded error estimator and PID controllers
           in OrdinaryDiffEq.jl and deviations in PETSc/TS.
           }
  \label{alg:time_loop_more_detailed}
  \begin{algorithmic}[1]
    \State Initialize time $t \gets 0$, time step number $n \gets 0$,
           and initial state $u^0$
    \State Initialize PID controller with $\epsilon_0 \gets 1$, $\epsilon_{-1} \gets 1$
    \State Initialize $\dt_{0} = \widetilde{\dt}$ with a given value
           or the algorithm of \cite[p.~169]{hairer2008solving}
    \State Initialize \texttt{accept\_step} $\gets$ \texttt{false}
    \While{$t < T$}
        \If{\texttt{accept\_step}}
          \Comment{callbacks can change \texttt{accept\_step}}
            \State \texttt{accept\_step} $\gets$ \texttt{false} \Comment{prepare for the next step}
            \State $t \gets \widetilde{t}$
            \State $\dt_{n+1} \gets \widetilde{\dt}$

          \State $n \gets n + 1$
        \Else
            \State $\dt_{n} \gets \widetilde{\dt}$ 
        \EndIf


          \If{$t + \dt_{n} > T$}
            \State $\dt_{n} \gets T - t$
          \EndIf

        \State Compute $u^{n+1}$ and $\uhat^{n+1}$ with time step $\dt_n$
        \State $w_{n+1} \gets \sqrt{ {\displaystyle \frac{1}{m} \sum_{i=1}^m} \left( \frac{u^{n+1}_i - \uhat^{n+1}_i}{\atol + \rtol \max\bigl\{ \bigl\vert u^{n+1}_i\bigr\vert, \bigl\vert \widetilde{u}_i\bigr\vert \bigr\}} \right)^2 }$ \State\Comment{$\widetilde{u} = u^{n}$ in OrdinaryDiffEq.jl, $\widetilde{u} = \uhat^{n+1}$ in PETSc}

          \State $w_{n+1} \gets \max\{w_{n+1}, w_\mathrm{min}\}$
              \State\Comment{$w_\mathrm{min} \approx \num{2.2e-16}$ in OrdinaryDiffEq.jl, $= \num{1.0e-10}$ in PETSc}
          \State $\epsilon_{n+1} \gets 1 / w_{n+1}$
          \State \texttt{dt\_factor} $\gets \kappa\Bigl( \epsilon_{n+1}^{\beta_1 / k}
                                                         \epsilon_{n  }^{\beta_2 / k}
                                                         \epsilon_{n-1}^{\beta_3 / k} \Bigr)$
            \Comment{default $\kappa(a) = 1 + \arctan(a - 1)$}
          \State $\widetilde{\dt} \gets$ \texttt{dt\_factor} $\cdot \dt_{n}$

        \If{\texttt{dt\_factor} $\ge$ \texttt{accept\_safety}}
            \Comment{default \texttt{accept\_safety} $= 0.81$}
          \State \texttt{accept\_step} $\gets$ \texttt{true}
        \Else
          \State \texttt{accept\_step} $\gets$ \texttt{false}
        \EndIf

        \If{\texttt{accept\_step}}

          \State $\widetilde{t} \gets t + \dt_{n}$
          \If{$\widetilde{t} \approx T$} \Comment{within 100 units in last place in OrdinaryDiffEq.jl}
            \State $\widetilde{t} \gets T$
          \EndIf


            \State Apply callbacks \Comment{AMR, CFL-based control in OrdinaryDiffEq.jl}
        \EndIf
    \EndWhile

  \end{algorithmic}
\end{algorithm}

\subsection{Importance of controller parameters}

The RK pair \eqref{eq:RK-step} and the PID controller \eqref{eq:PID} must be
developed together to obtain good results. In particular, semidiscretizations of
conservation laws limited by stability should be integrated using a
combination of RK pair and controller leading to stable step size control,
based on the test equation \cite{hall1988analysis}, see also
\cite[Section~IV.2]{hairer2010solving} and \cite{kennedy2000low}. This has been
demonstrated for spectral element discretizations for compressible flows in
\cite{ranocha2021optimized}. In this article, we use the optimized controllers
of \cite{ranocha2021optimized}.

\subsection{Representative Runge-Kutta methods used in this article}
\label{sec:RK-representative}

There is a vast amount of literature on RK methods. Many schemes have been
designed as general-purpose methods for low- or high-tolerance applications
or even specifically for hyperbolic conservation laws. While we have performed
numerical experiments with various schemes, we restrict this article to the
following representative set of methods.
First, we consider the general-purpose method
\begin{itemize}
  \item \RK[BS]{3}{3}[][FSAL],
        the third-order, three-stage method of \cite{bogacki1989a32}
        (\texttt{BS3()} in \mbox{OrdinaryDiffEq.jl}, \texttt{3bs} in PETSc)
\end{itemize}
This method has been shown to be among the most efficient --- if not the most
efficient --- general purpose method for the CFD problems in which we are interested
\cite{ranocha2021optimized}.
Next, we consider methods optimized for spectral element semidiscretizations
of compressible CFD, namely
\begin{itemize}
  \item \RK[RDPK]{3}{5}[\ESstarp][FSAL],
        the third-order, five-stage method of \cite{ranocha2021optimized}
        (\texttt{RDPK3SpFSAL35()} in OrdinaryDiffEq.jl)

  \item \RK[RDPK]{4}{9}[\ESstarp][FSAL],
        the fourth-order, nine-stage method of \cite{ranocha2021optimized}
        (\texttt{RDPK3SpFSAL49()} in OrdinaryDiffEq.jl)
\end{itemize}
Finally, we consider the strong stability preserving (SSP) method
\begin{itemize}
  \item \ssp43,
        the third-order, four-stage SSP method of
        \cite{kraaijevanger1991contractivity}
        with embedded method of \cite{conde2022embedded} and
        efficient implementation and step size controller of \cite{ranocha2021optimized}
        (\texttt{SSPRK43()} in OrdinaryDiffEq.jl)
\end{itemize}

\section{Robustness under change of mesh structure}
\label{sec:meshes}

Both CFL- and error-based step size control come with parameters that must
be chosen by the user, either a CFL factor $\cfl$ or a tolerance $\tol$.
However, the sensitivity with respect to these parameters differs significantly.
The CFL factor $\nu$ influences the time step sizes --- and thus the efficiency ---
linearly, while the tolerance $\tol$ has roughly a logarithmic influence.
Thus, it is arguably easier to choose an appropriate tolerance than
an optimal CFL factor.

Moreover, in the stability limited regime, it is often convenient to use error-based
step size control since there is usually a range of tolerances resulting in the
same number of right hand side (RHS) evaluations a manually optimized CFL-control could achieve at best.
Furthermore, an optimal CFL factor $\nu$ depends on the mesh. Thus, it can vary
when introducing curved coordinates compared to a uniform grid.
This was already demonstrated for a 2D linear advection equation in
\cite[Section~3]{ranocha2021optimized}. Here, we extend this demonstration to
the nonlinear compressible Euler equations of an ideal gas in $\ndims = 3$ space
dimensions given by
\begin{equation}
  \partial_t \begin{pmatrix}
    \rho \\
    \rho v_i \\
    \rho e
  \end{pmatrix}
  + \sum_{j=1}^\ndims \partial_{x_j} \begin{pmatrix}
    \rho v_j \\
    \rho v_i v_j + p \delta_{ij} \\
    (\rho e + p) v_j
  \end{pmatrix}
  = 0,
\end{equation}
where $\rho$ is the density, $v$ the velocity, $\rho e$ the total energy, and
$p$ the pressure given by
\begin{equation}
  p = (\gamma - 1) \left( \rho e - \frac{1}{2} \rho v^2 \right),
\end{equation}
with the ratio of specific heats $\gamma = 1.4$.
Specifically, we consider the classical inviscid Taylor-Green vortex with initial data\footnote{We use a superscript $0$ to denote initial data.}
\begin{equation}
\begin{gathered}
%
%
  \rho^0 = 1, \;
  v_1^0 =  \sin(x_1) \cos(x_2) \cos(x_3), \;
  v_2^0 = -\cos(x_1) \sin(x_2) \cos(x_3), \;
  v_3^0  = 0, \\
  p^0 = \frac{\rho^0}{\mathrm{Ma}^2 \gamma} + \rho^0 \frac{\cos(2 x_1) \cos(2 x_3) + 2 \cos(2 x_2) + 2 \cos(2 x_1) + \cos(2 x_2) \cos(2 x_3)}{16},
\end{gathered}
\end{equation}
where $\mathrm{Ma} = 0.1$ is the Mach number.
We consider the domain $[-\pi, \pi]^3$ with periodic boundary
conditions and meshes with $8$ elements per coordinate direction. We apply
entropy-dissipative discontinous Galerkin (DG) methods with polynomials of degree $\polydeg = 3$ using
a local Lax-Friedrichs flux at interfaces and the flux of
\cite{ranocha2020entropy,ranocha2018thesis,ranocha2021preventing}
in the volume terms. We integrate the semidiscretizations in the time interval
$[0, 10]$.

\begin{figure}[htb]
\centering
  \includegraphics[width=0.7\textwidth]{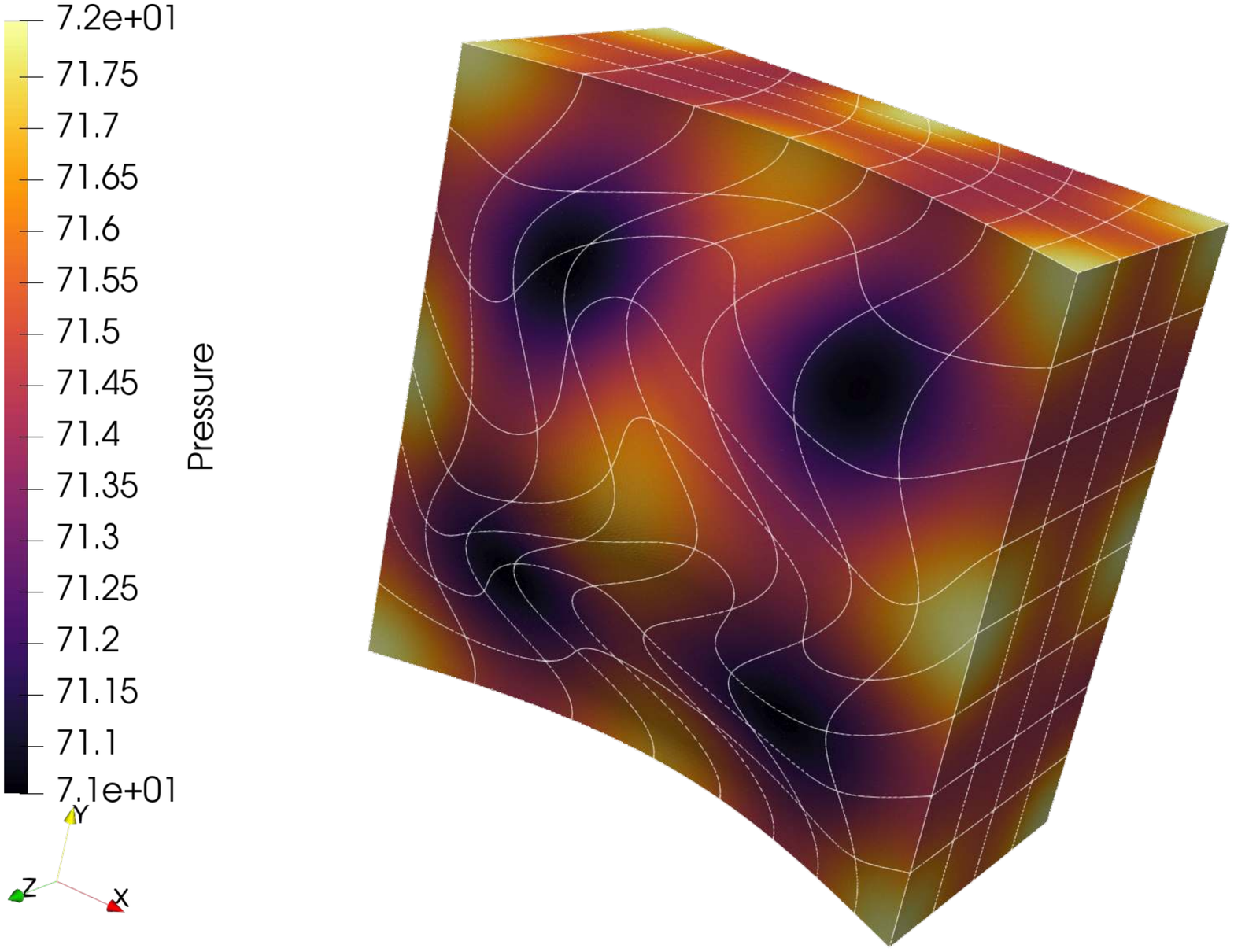}
  \caption{Initial pressure of the compressible Euler equations for the inviscid
           Taylor-Green vortex on a slice of the curved mesh.}
  \label{fig:taylor_green_vortex_StructuredMesh_pressure}
\end{figure}

We consider two meshes, a uniform Cartesian mesh and a curved mesh
that heavily warps the Cartesian reference coordinates $(\xi, \eta, \zeta) \in [-1, 1]^3$
to the desired mesh in physical coordinates $(x, y, z)$ with the mapping
\begin{equation}
%
%
%
%
%
\begin{gathered}
\begin{aligned}
  y &= \eta  + \frac{L_y}{8} \bigl( \cos(3 \pi (\xi   - c_x) / L_x)
                                    \cos(  \pi (\eta  - c_y) / L_y)
                                    \cos(  \pi (\zeta - c_z) / L_z \bigr),
  \\
  x &= \xi   + \frac{L_x}{8} \bigl( \cos(  \pi (\xi   - c_x) / L_x)
                                    \cos(4 \pi (y     - c_y) / L_y)
                                    \cos(  \pi (\zeta - c_z) / L_z \bigr),
  \\
  z &= \zeta + \frac{L_y}{8} \bigl( \cos(  \pi (x     - c_x) / L_x)
                                    \cos(2 \pi (y     - c_y) / L_y)
                                    \cos(  \pi (\zeta - c_z) / L_z \bigr),
\end{aligned}
  \\
  L_x = L_y = L_z = 2 \pi, \quad
  c_x = c_y = c_z = 0,
\end{gathered}
\end{equation}
that has been adapted from \cite{chan2019efficient}.
The curved mesh and the initial pressure are visualized
in Figure~\ref{fig:taylor_green_vortex_StructuredMesh_pressure}.

\begin{figure}[htb]
\centering
  \begin{subfigure}{0.48\linewidth}
    \includegraphics[width=\textwidth]{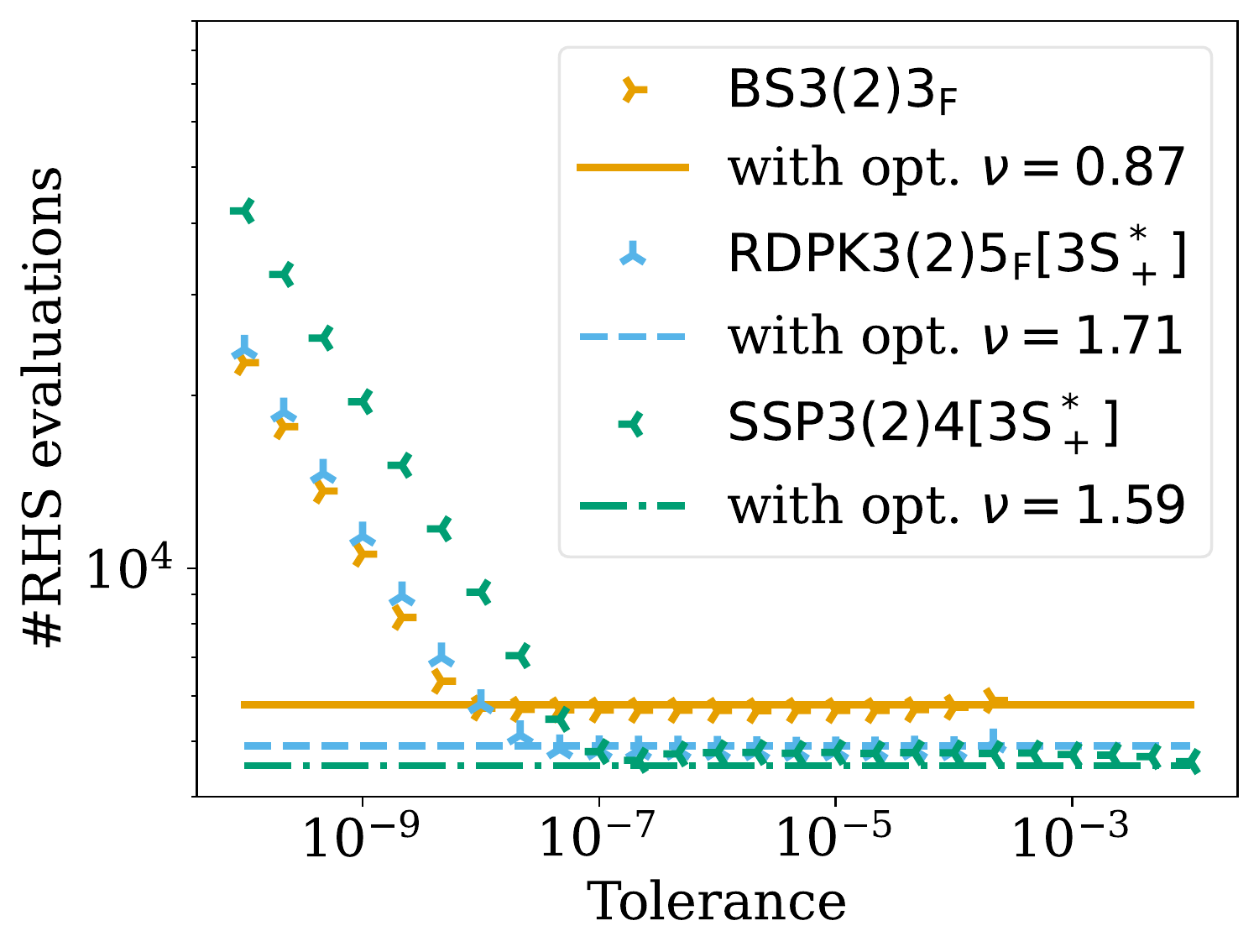}
    \caption{Cartesian mesh.}
  \end{subfigure}%
  \hspace*{\fill}
  \begin{subfigure}{0.48\linewidth}
    \includegraphics[width=\textwidth]{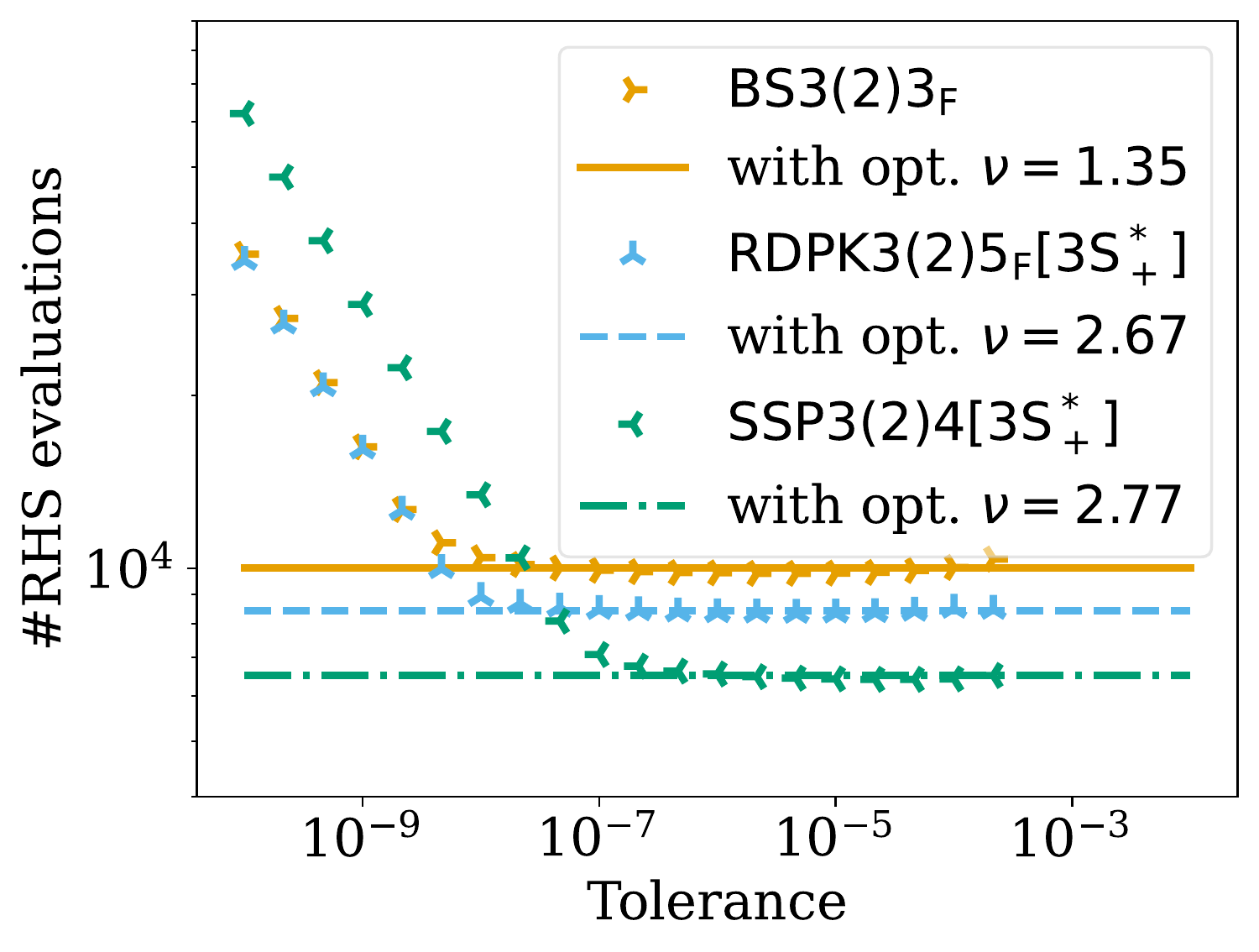}
    \caption{Curved mesh.}
  \end{subfigure}%
  \caption{Number of ODE RHS evaluations for representative RK methods using
           error-based (markers) and CFL-based (lines) step size control.
           The CFL factors were tuned manually to be as large as possible without
           crashing the simulation. Note that the optimized CFL factors are between
           \SI{55}{\percent} and \SI{74}{\percent} larger on the curved grid.}
  \label{fig:taylor_green_vortex}
\end{figure}

The required number of RHS evaluations for the selected, representative RK methods are shown in
Figure~\ref{fig:taylor_green_vortex}, for CFL-based and error-based step
size control. All RK methods show the same qualitative behavior on both grids.
The step size for this problem is largely limited by stability, resulting
in a range of tolerances yielding the optimal number of RHS evaluations one can
also achieve by tuning the CFL factor manually (maximizing $\cfl$ to three
significant digits so that the simulations do not crash). This range of tolerances
spans up to five orders of magnitude. For looser tolerances, the simulations may
crash due to positivity issues, typically in the first few time step. Stricter
tolerances let the controller detect an accuracy-limited regime and increase the
number of RHS evaluations accordingly. This happens at looser tolerances for
\ssp43 than for the other methods since \ssp43 is only optimized for the SSP
coefficient while the other methods are also constructed to minimize the principal
truncation error. Since all methods are of the same order of accuracy, the slope
of the increasing number of RHS evaluations is the same.

However, there is a significant difference between CFL-based and error-based step
size control when switching the mesh types. Indeed, the CFL factor can be
increased between \SI{55}{\percent} (for \RK[BS]{3}{3}[][FSAL]) and
\SI{74}{\percent} (for \ssp43) on the curved mesh compared to the equivalent
Cartesian mesh. In particular, choosing an optimized CFL factor from the uniform
mesh results in the same amount of efficiency loss when the CFL factor is not
tuned again manually from scratch. In contrast, basically the same range of
tolerances can be used to get the optimal number of RHS evaluations on both
grids for the error-based step size control.

\section{Change of linear stability restrictions in nonlinear schemes}
\label{sec:dg_fv}

There are many shock capturing approaches for DG methods, e.g.,
artificial viscosity \cite{persson2006sub,guermond2011entropy},
replacing DG elements by finite volume subcells \cite{dumbser2014posteriori,sonntag2014shock},
modal filtering \cite{meister2013extended,ranocha2018stability},
or specially constructed invariant domain preserving methods
\cite{pazner2021sparse,guermond2018second,kuzmin2021entropy}.
Here, we use the a priori convex blending of high-order DG elements with finite
volume subcells described in \cite{hennemann2021provably}.

\subsection{Linear CFL restrictions}
\label{sec:spectra}

Shock capturing approaches in nonlinear schemes will typically change the linear
CFL restriction on the time step for explicit time integration schemes. Here,
we demonstrate this for some third-order schemes of the representative RK
methods introduced in Section~\ref{sec:RK-representative}.

We consider the linear advection equation
\begin{equation}
  \partial_t u(t,x) + \sum_{j=1}^2 a^j \partial_j u(t,x) = 0
\end{equation}
in the domain $[-1, 1]^2$ with periodic boundary conditions and velocity
$a = (1, 1)^T / \sqrt{2}$. We discretize the domain using $8^2$ uniform elements
with polynomials of degree $\polydeg = 3$ on Legendre-Gauss-Lobatto nodes.
We choose a fixed blending parameter $\alpha = 0.5$ in the shock capturing scheme
of \cite{hennemann2021provably} with the standard discontinuous Galerkin collocation spectral
element method (DGSEM), e.g., \cite{kopriva2009implementing}, and local
Lax-Friedrichs flux for the finite volume subcells. We compute the spectra of
the standard DGSEM scheme and the shock capturing scheme. Due to the fixed choice
of the blending parameter, both semidiscretizations are linear and we compute
their spectra numerically.

\begin{figure}[!htb]
\centering
  \begin{subfigure}{0.33\linewidth}
    \includegraphics[width=\textwidth]{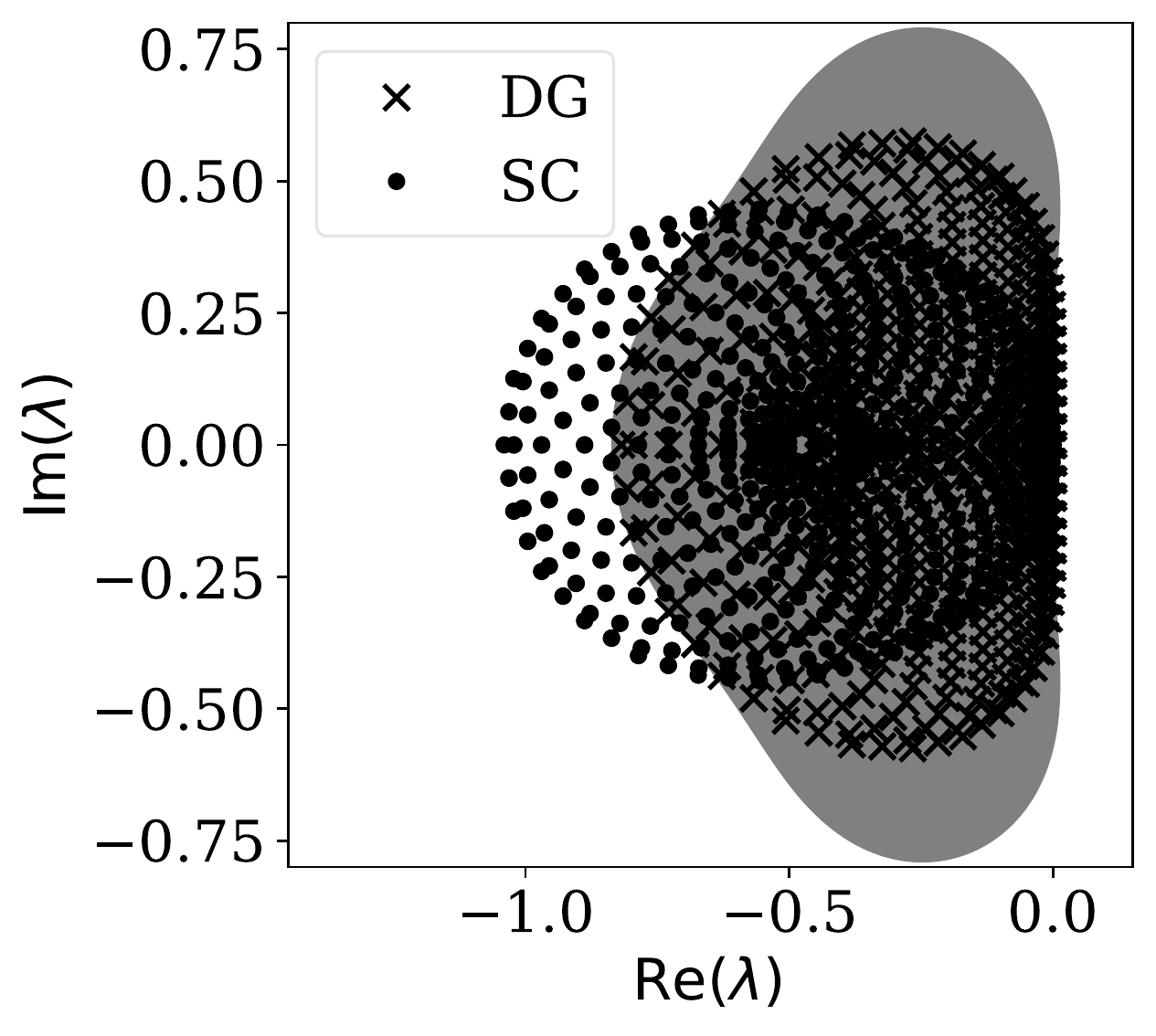}
    \caption{\RK[BS]{3}{3}[][FSAL].}
  \end{subfigure}%
  \hspace*{\fill}
  \begin{subfigure}{0.33\linewidth}
    \includegraphics[width=\textwidth]{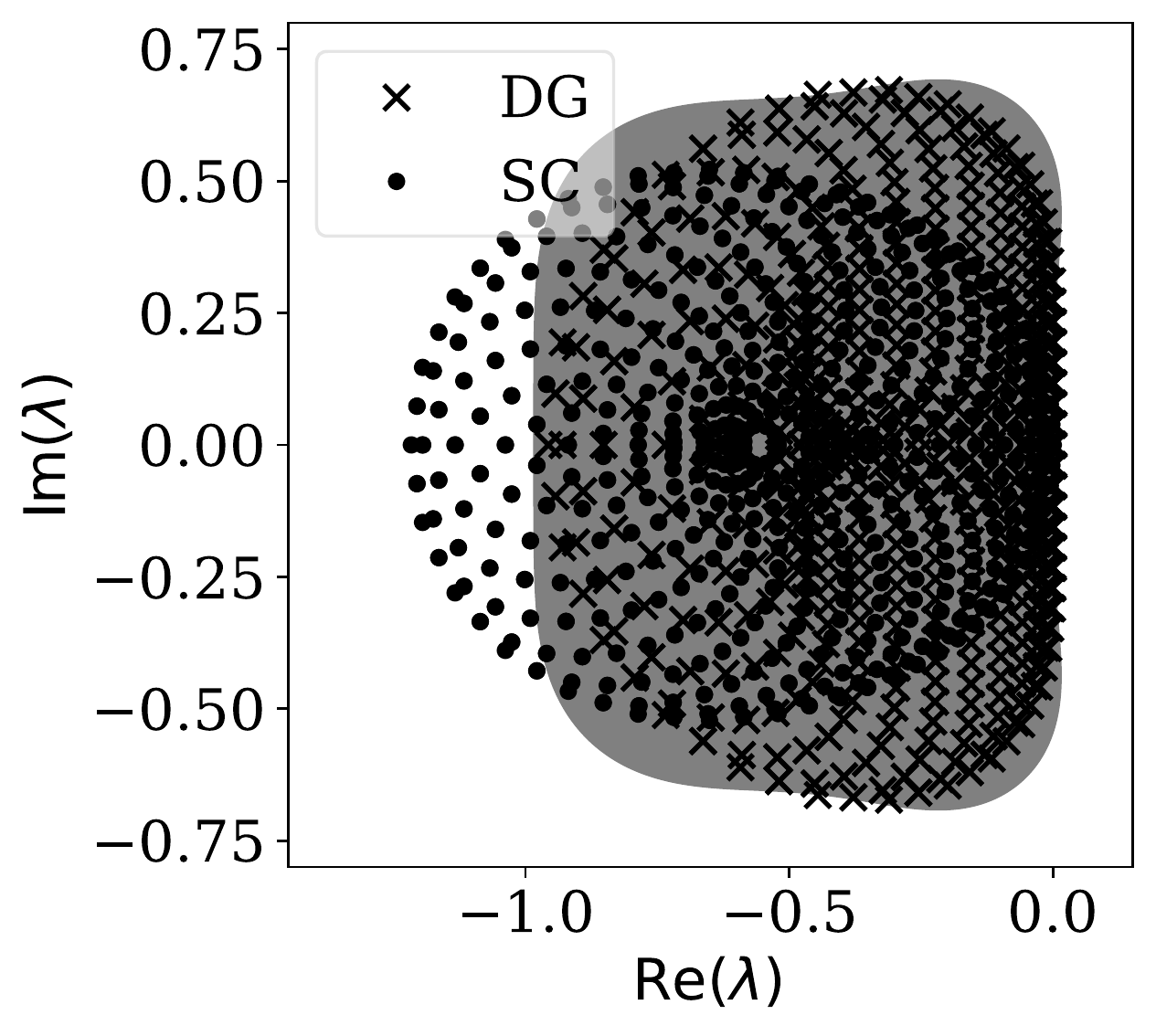}
    \caption{\RK[RDPK]{3}{5}[\ESstarp][FSAL].}
  \end{subfigure}%
  \hspace*{\fill}
  \begin{subfigure}{0.33\linewidth}
    \includegraphics[width=\textwidth]{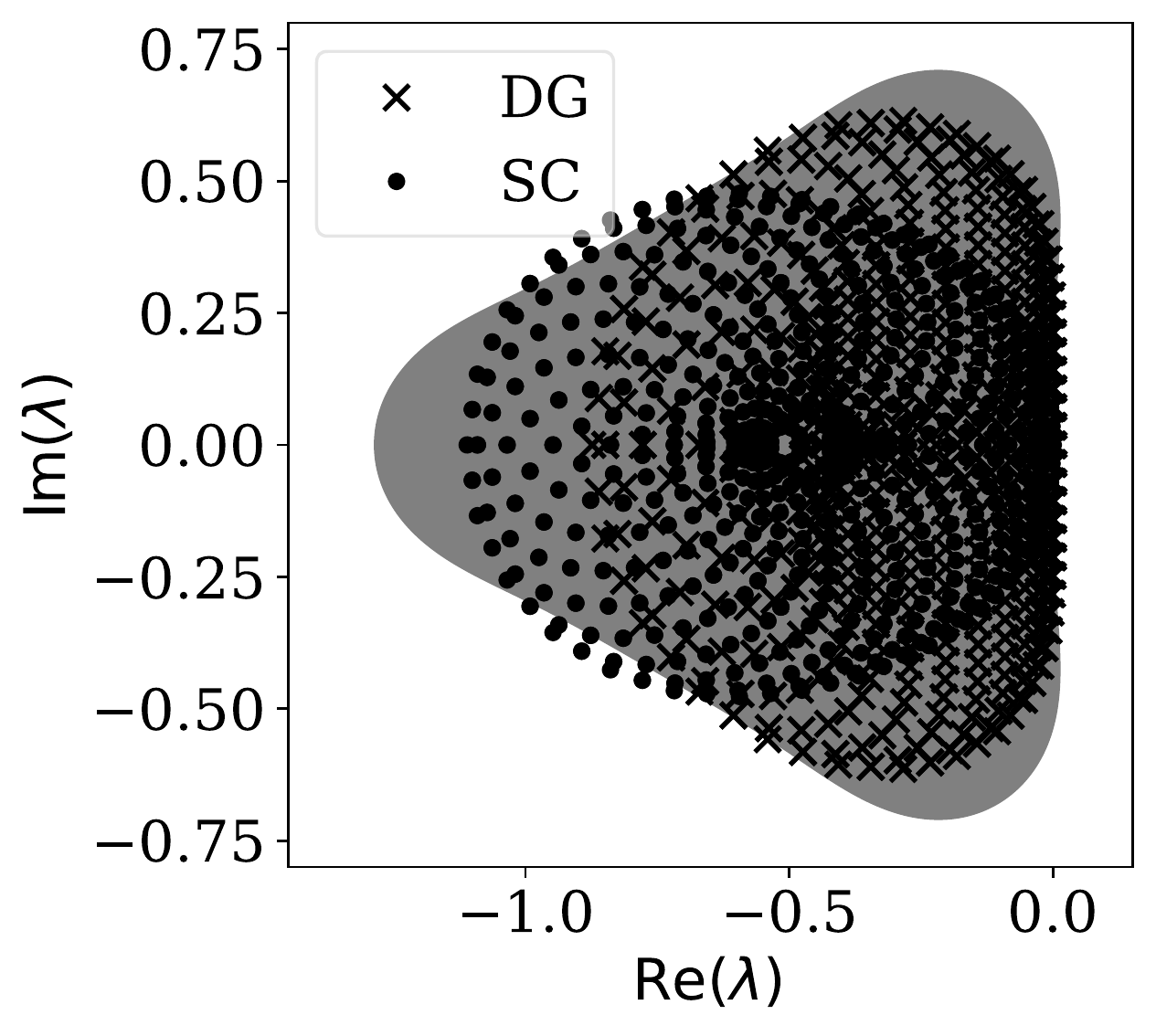}
    \caption{\ssp43.}
  \end{subfigure}%
  \caption{Spectra of DGSEM and shock capturing (SC)
           semidiscretizations of the linear advection equation are embedded
           into the stability regions of some representative Runge-Kutta methods.
           The stability regions are scaled by the effective number of stages of
           the RK methods (taking the FSAL property into account). The spectra
           of the semidiscretizations are scaled by the maximal factor so that
           the standard DGSEM spectrum is in the region of absolute stability.}
  \label{fig:spectra-DGSEM-SC}
\end{figure}

The spectra embedded into the stability regions of the representative RK methods
are visualized in Figure~\ref{fig:spectra-DGSEM-SC}. To make it easier to compare
the RK methods, the stability regions are scaled by the effective number of
RK stages, taking the FSAL property into account.
The spectra of the semidiscretizations are scaled by the maximal factor so that
the standard DGSEM spectrum is in the region of absolute stability of the RK method.
Clearly, the (scaled) spectra of the shock capturing semidiscretizations are
partially outside of the stability regions. Thus, one can expect a more restrictive
CFL condition than for the standard DGSEM schemes.

These linear CFL restrictions predict an approximately one quarter bigger CFL number
for standard DGSEM discretizations compared to the shock capturing variants for
\RK[BS]{3}{3}[][FSAL] and \RK[RDPK]{3}{5}[\ESstarp][FSAL]. For these time integration methods,
the SC spectra extend significantly to the left-half of the complex plane outside
of the stability regions, in particular near the negative real axis.
The effect is qualitatively similar but quantitatively less pronounced for \ssp43
(\ca \SI{5}{\percent} instead of more than \SI{20}{\percent}).

\subsection{Application to nonlinear magnetohydrodynamics: Orszag-Tang vortex}
\label{sec:orszag_tang}

Next, we apply the representative time integration methods to
entropy-dissipative semidiscretizations of the ideal
generalized Lagrange multiplier (GLM) magnetohydrodynamics (MHD) equations
\cite{derigs2018ideal,bohm2018entropy}
for the Orszag-Tang vortex setup \cite{orszag1979small}.
Specifically, we use the initial condition
\begin{equation}
\begin{gathered}
  \rho^0 = 1, \quad
  v_1^0 = -\sin(2\pi x_2), \quad
  v_2^0 = \sin(2\pi x_1), \quad
  v_3^0 = 0, \quad
  p^0 = 1 / \gamma, \\
  B_1^0 = -\sin(2\pi x_2) / \gamma, \quad
  B_2^0 = \sin(4\pi x_1) / \gamma, \quad
  B_3^0 = 0, \quad
  \psi^0 = 0,
\end{gathered}
\end{equation}
in the domain $[0, 1]^2$ with periodic boundary conditions, where $\rho$ is the
density, $v$ the velocity, $p$ the pressure, $B$ the magnetic field, and $\psi$
the Lagrange multiplier to control divergence errors.
The ideal GLM-MHD equations with ratio of specific heats
$\gamma = 5 / 3$ and divergence cleaning speed $c_h = 1$ are discretized on a
uniform mesh with $2^6 = 64$ elements per
coordinate direction and polynomials of degree $\polydeg = 3$. We use the
shock capturing method of \cite{hennemann2021provably} with the entropy-conservative
numerical flux of \cite{hindenlang2019entropy} and a local Lax-Friedrichs flux
at interfaces and for finite volume subcells, both with the nonconservative
Powell source term. We use the product of density and pressure as the shock indicator
variable with blending parameters $\alpha_{\mathrm{max}} = 0.5$ and $\alpha_{\mathrm{min}}
= 0.001$.

\begin{figure}[!htb]
\centering
  \begin{subfigure}{0.31\linewidth}
    \includegraphics[width=\textwidth]{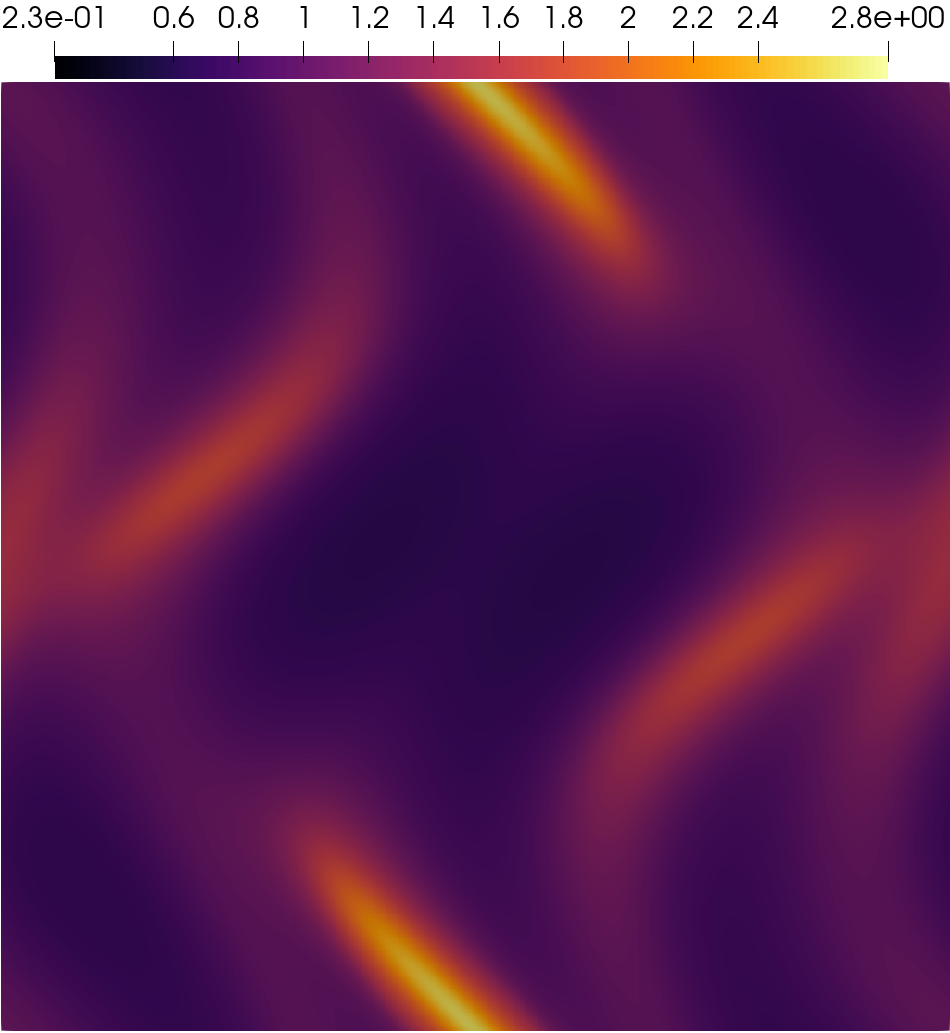}
    \caption{Density at $t \approx 0.11$.}
  \end{subfigure}%
  \hspace*{\fill}
  \begin{subfigure}{0.31\linewidth}
    \includegraphics[width=\textwidth]{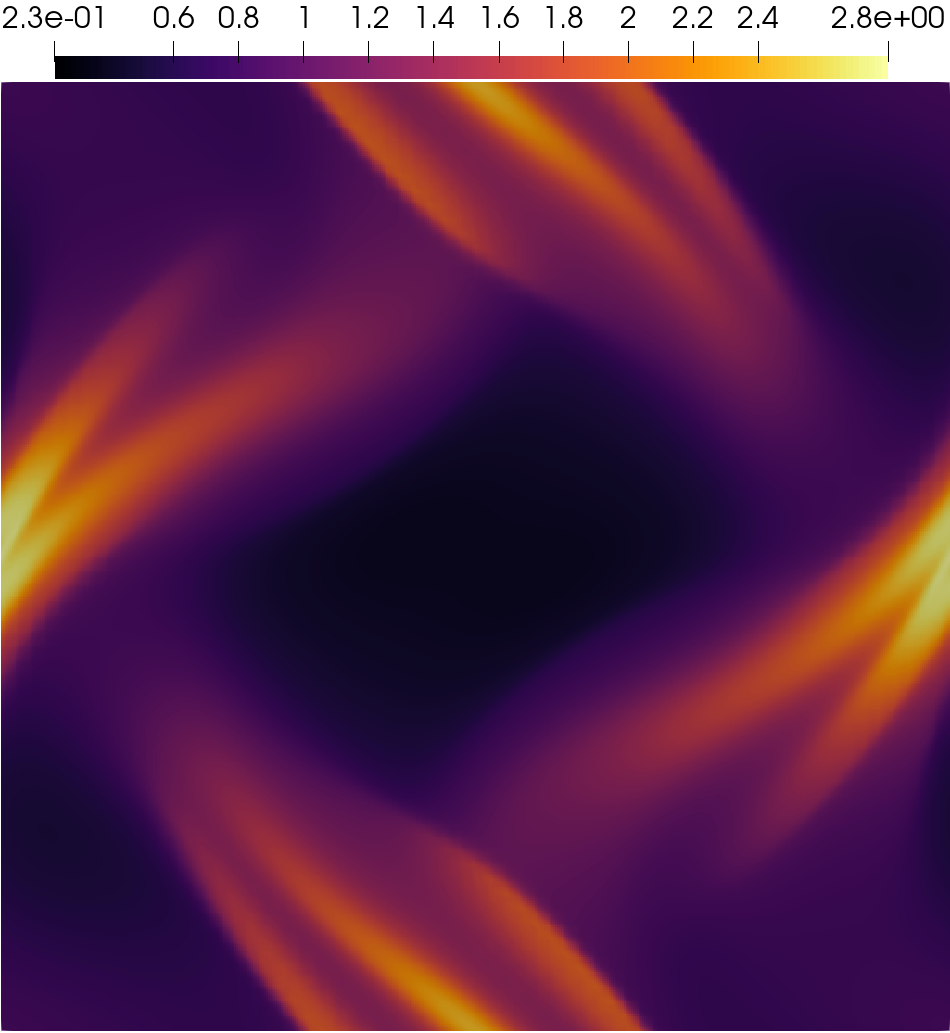}
    \caption{Density at $t \approx 0.19$.}
  \end{subfigure}%
  \hspace*{\fill}
  \begin{subfigure}{0.31\linewidth}
    \includegraphics[width=\textwidth]{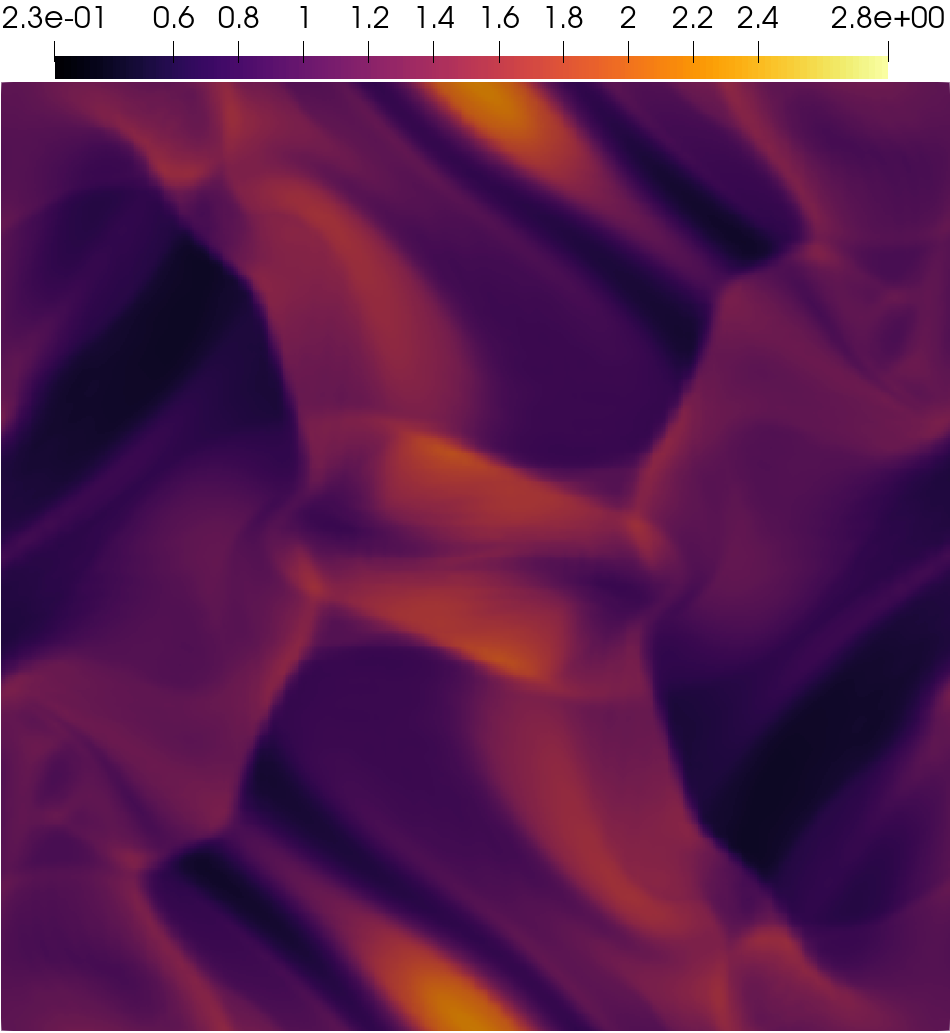}
    \caption{Density at $t = 0.50$.}
  \end{subfigure}%
  \\
  \begin{subfigure}{0.31\linewidth}
    \includegraphics[width=\textwidth]{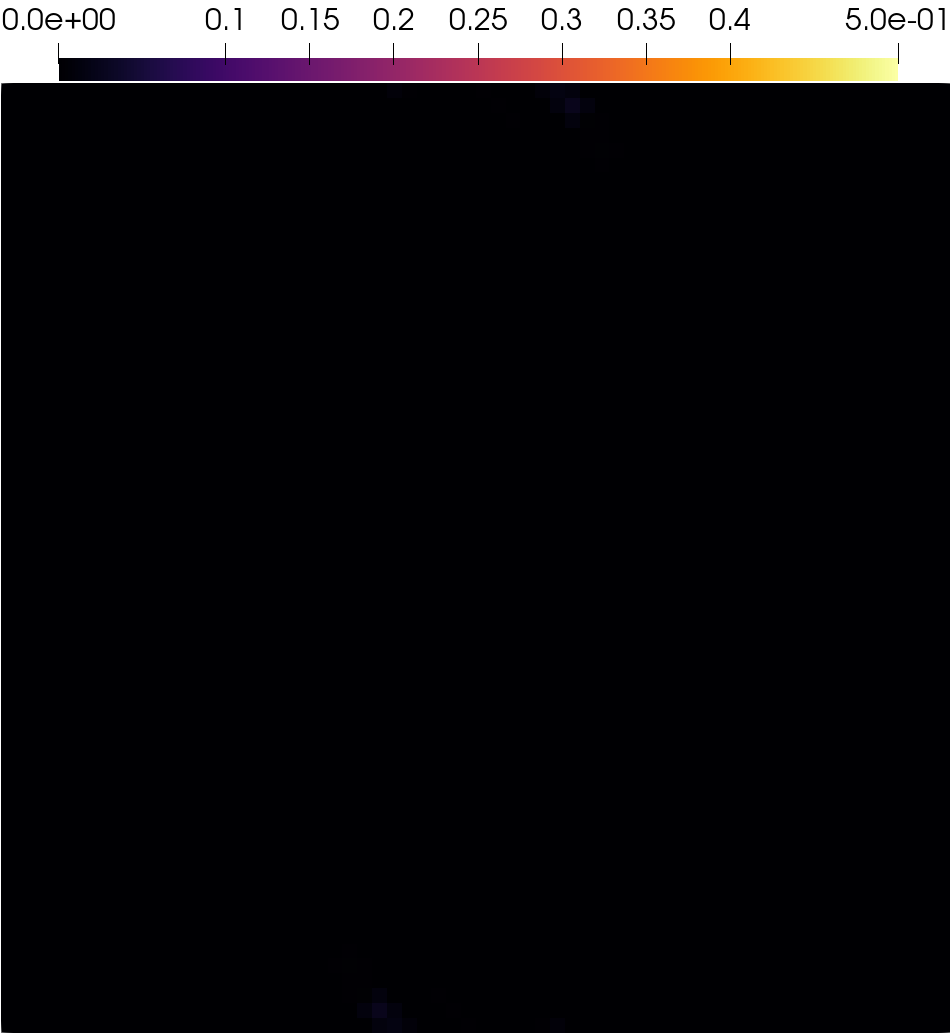}
    \caption{SCI at $t \approx 0.11$.}
  \end{subfigure}%
  \hspace*{\fill}
  \begin{subfigure}{0.31\linewidth}
    \includegraphics[width=\textwidth]{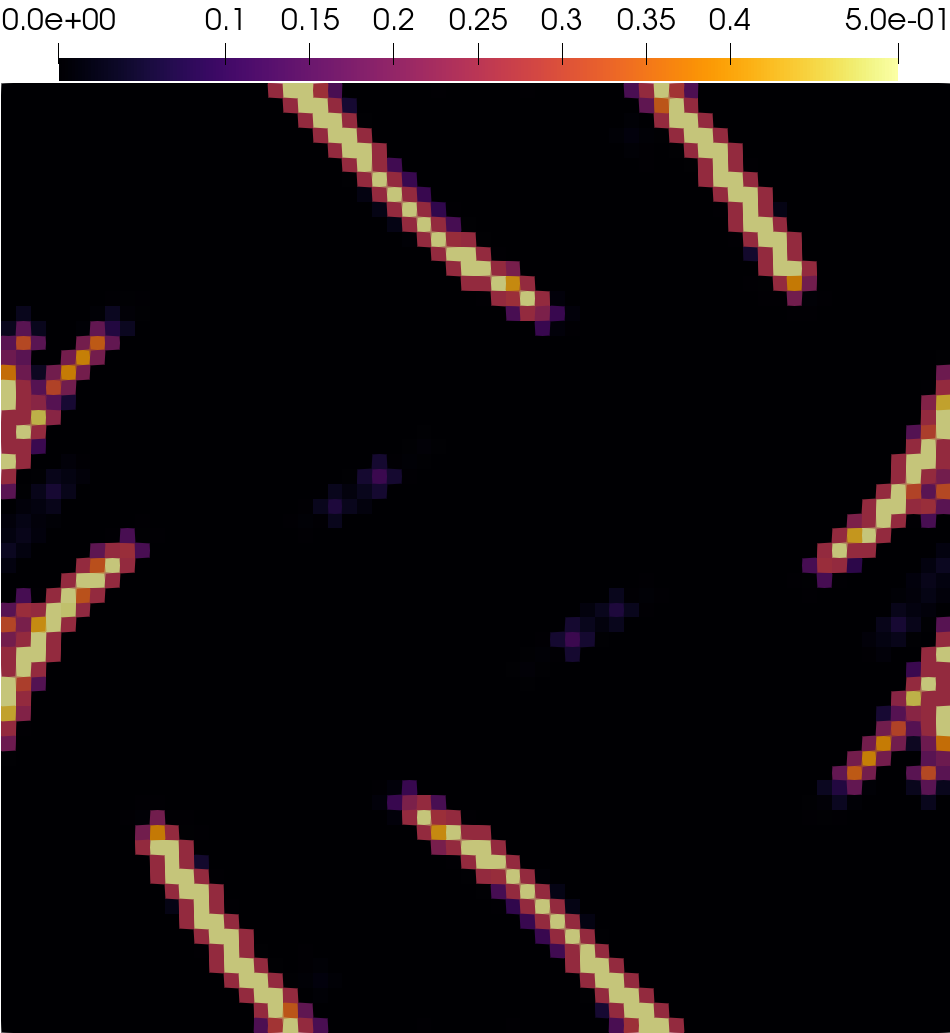}
    \caption{SCI at $t \approx 0.19$.}
  \end{subfigure}%
  \hspace*{\fill}
  \begin{subfigure}{0.31\linewidth}
    \includegraphics[width=\textwidth]{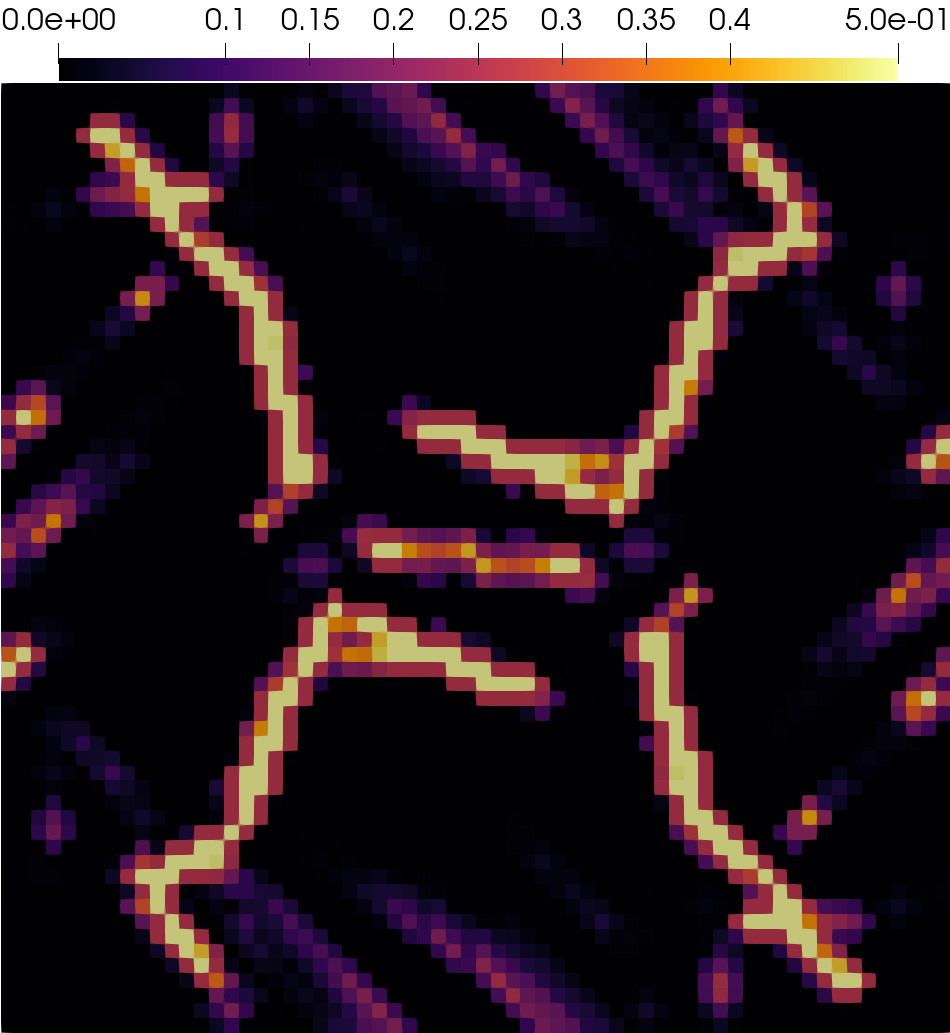}
    \caption{SCI at $t = 0.50$.}
  \end{subfigure}%
  \caption{Evolution of the density and the shock capturing indicator (SCI)
           computed using \ssp43 and entropy-dissipative shock capturing
           semidiscretizations of the ideal GLM-MHD equations for the
           Orszag-Tang vortex.}
  \label{fig:orszag_tang_snapshots}
\end{figure}

As shown in Figure~\ref{fig:orszag_tang_snapshots}, the flow and its numerical
approximation remain smooth in its initial phase. Between $t = 0.1$
and $t = 0.2$, the shock capturing indicator detects troubled cells and activates
the finite volume shock capturing mechanism. At the final time, several shocklets
are visible in the approximation.

\begin{figure}[!htb]
\centering
  \begin{subfigure}{0.33\linewidth}
    \includegraphics[width=\textwidth]{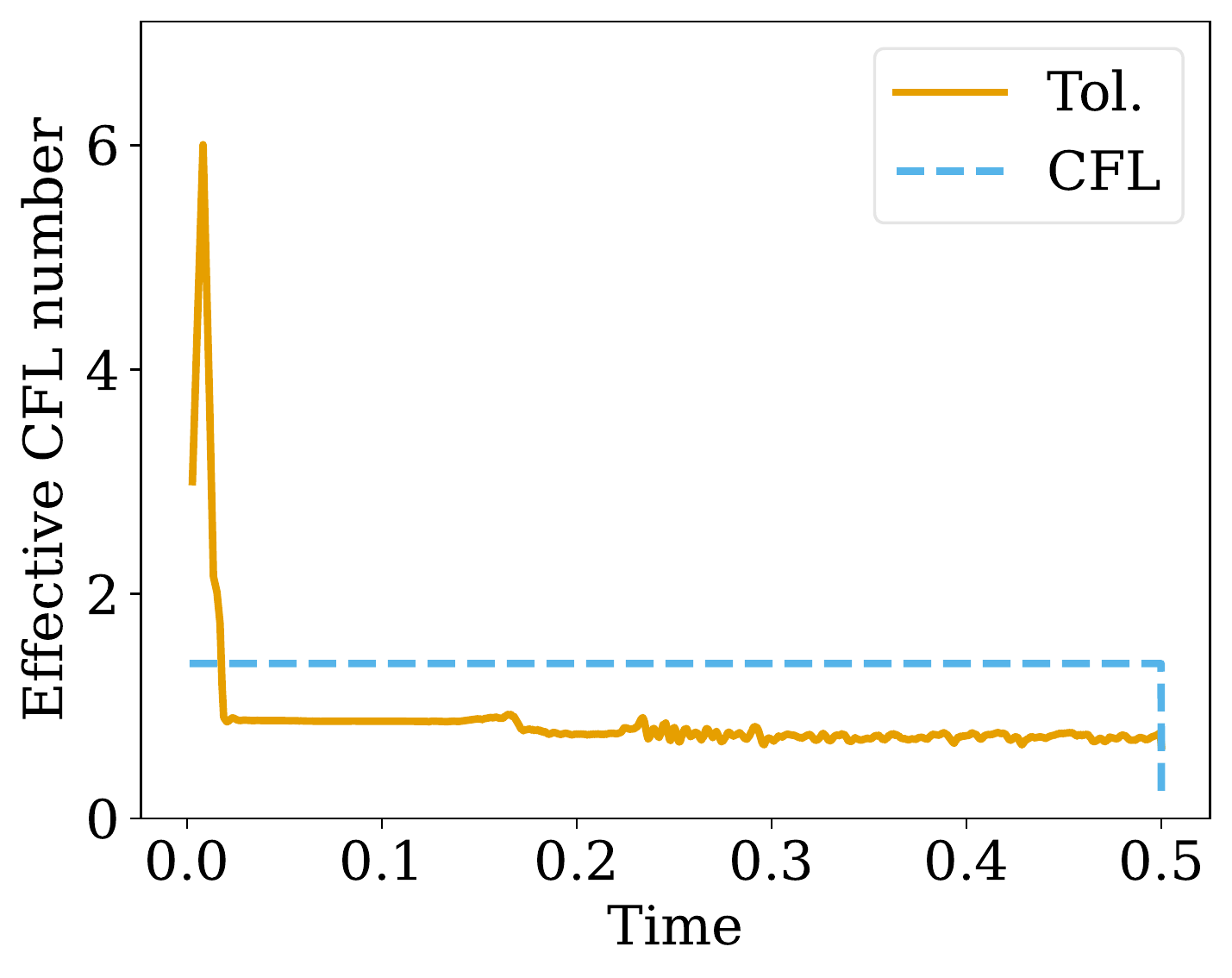}
    \caption{\RK[BS]{3}{3}[][FSAL].}
  \end{subfigure}%
  \hspace*{\fill}
  \begin{subfigure}{0.33\linewidth}
    \includegraphics[width=\textwidth]{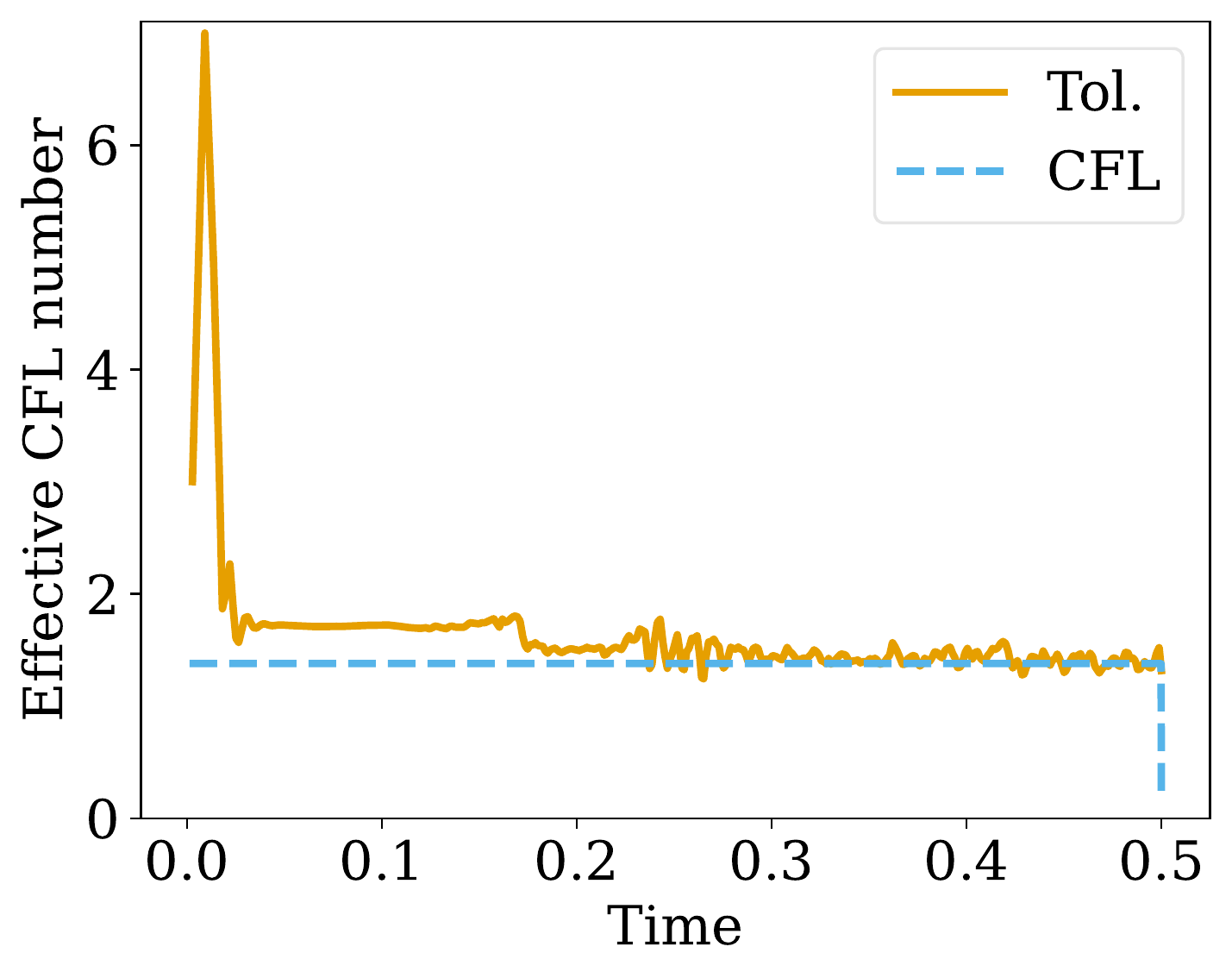}
    \caption{\RK[RDPK]{3}{5}[\ESstarp][FSAL].}
  \end{subfigure}%
  \hspace*{\fill}
  \begin{subfigure}{0.33\linewidth}
    \includegraphics[width=\textwidth]{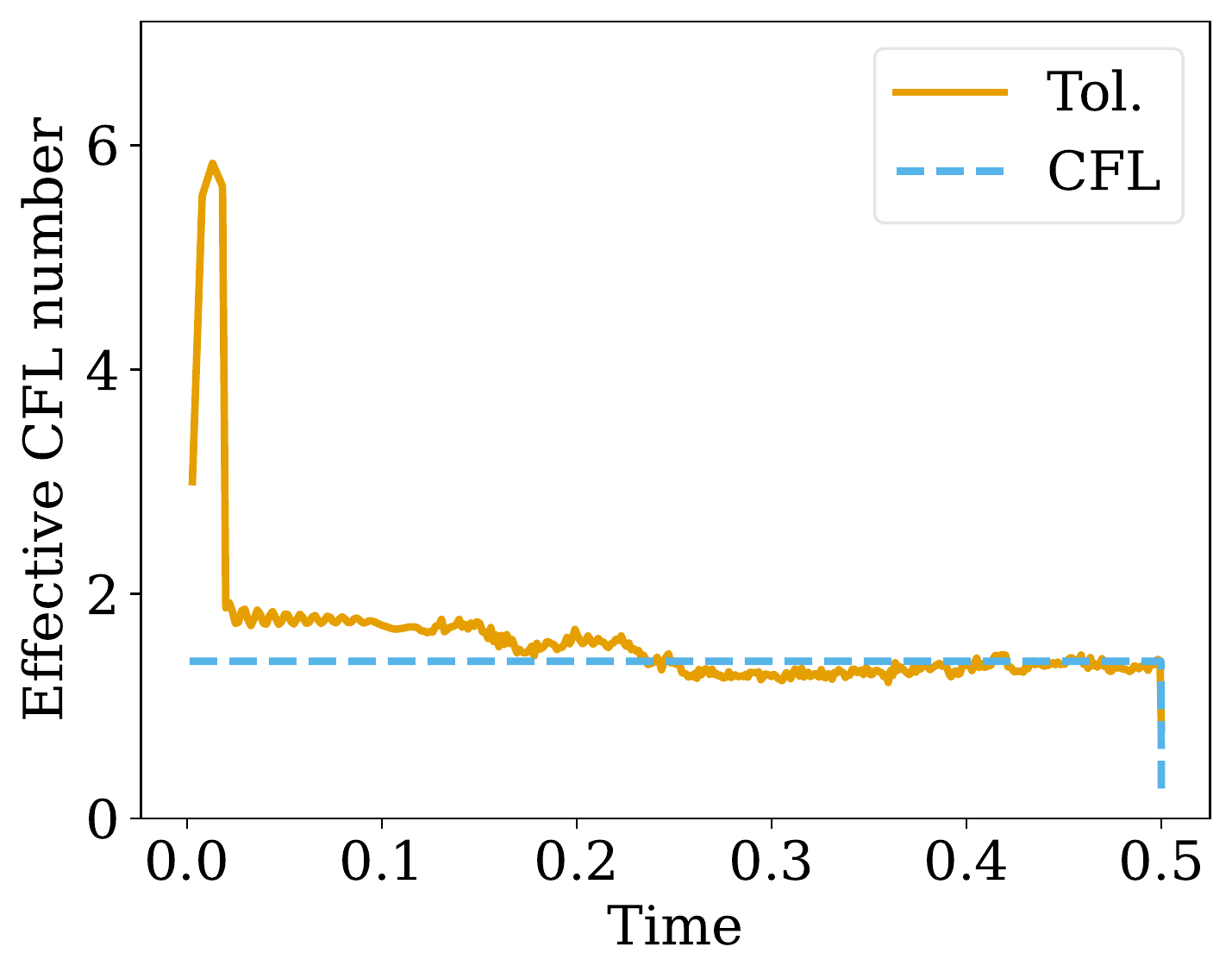}
    \caption{\ssp43.}
  \end{subfigure}%
  \\
  \begin{subfigure}{0.33\linewidth}
    \includegraphics[width=\textwidth]{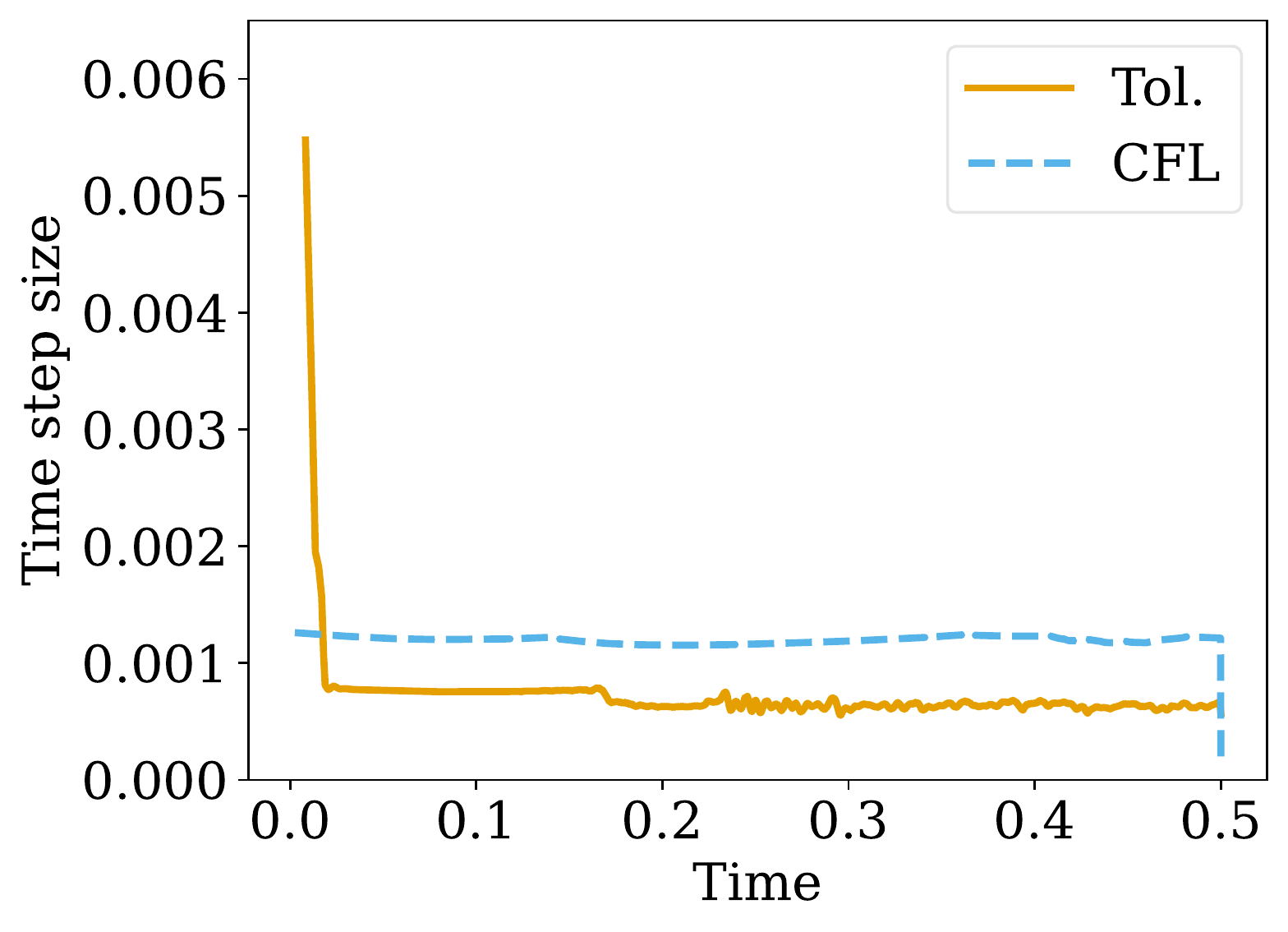}
    \caption{\RK[BS]{3}{3}[][FSAL].}
  \end{subfigure}%
  \hspace*{\fill}
  \begin{subfigure}{0.33\linewidth}
    \includegraphics[width=\textwidth]{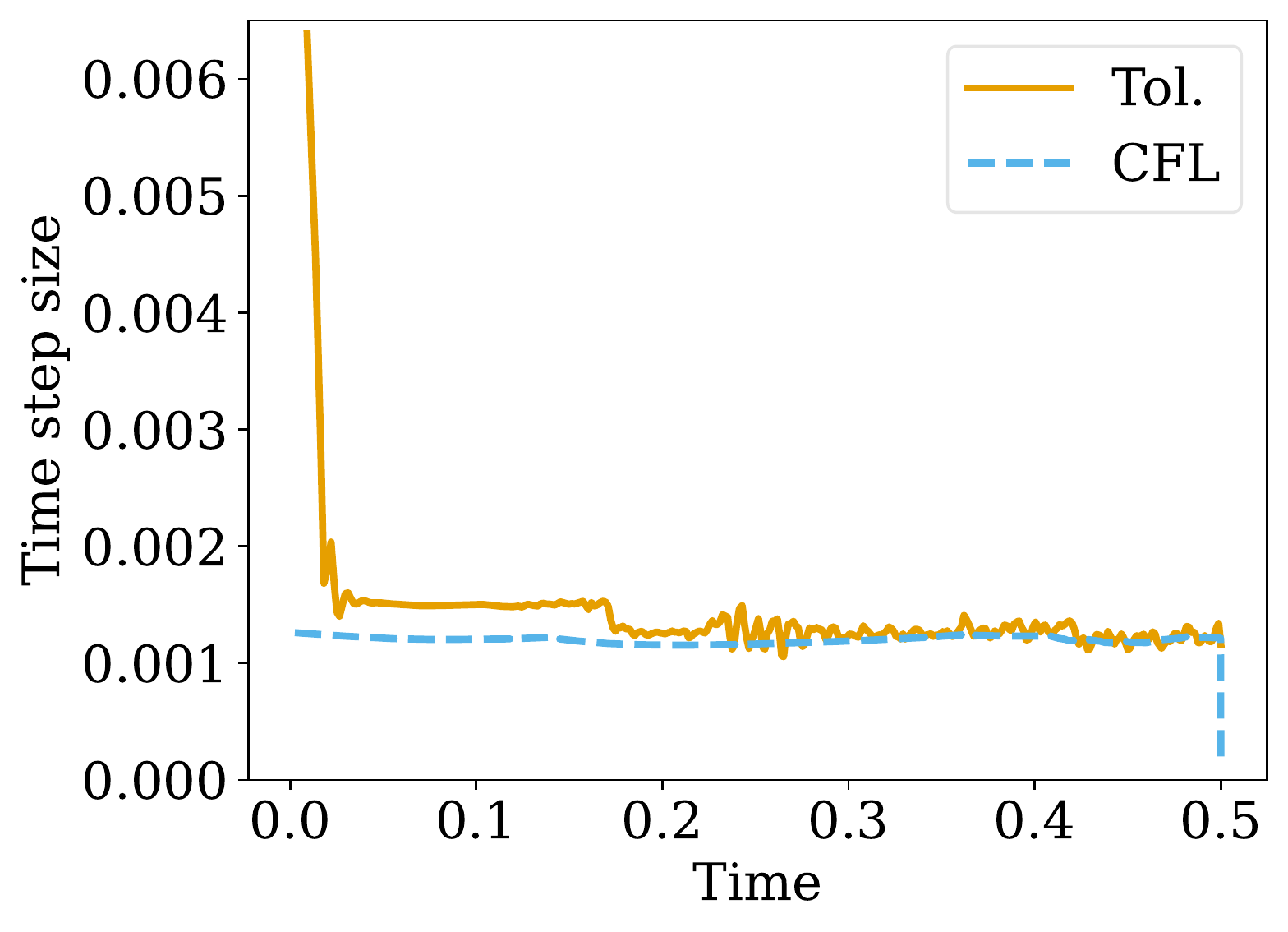}
    \caption{\RK[RDPK]{3}{5}[\ESstarp][FSAL].}
  \end{subfigure}%
  \hspace*{\fill}
  \begin{subfigure}{0.33\linewidth}
    \includegraphics[width=\textwidth]{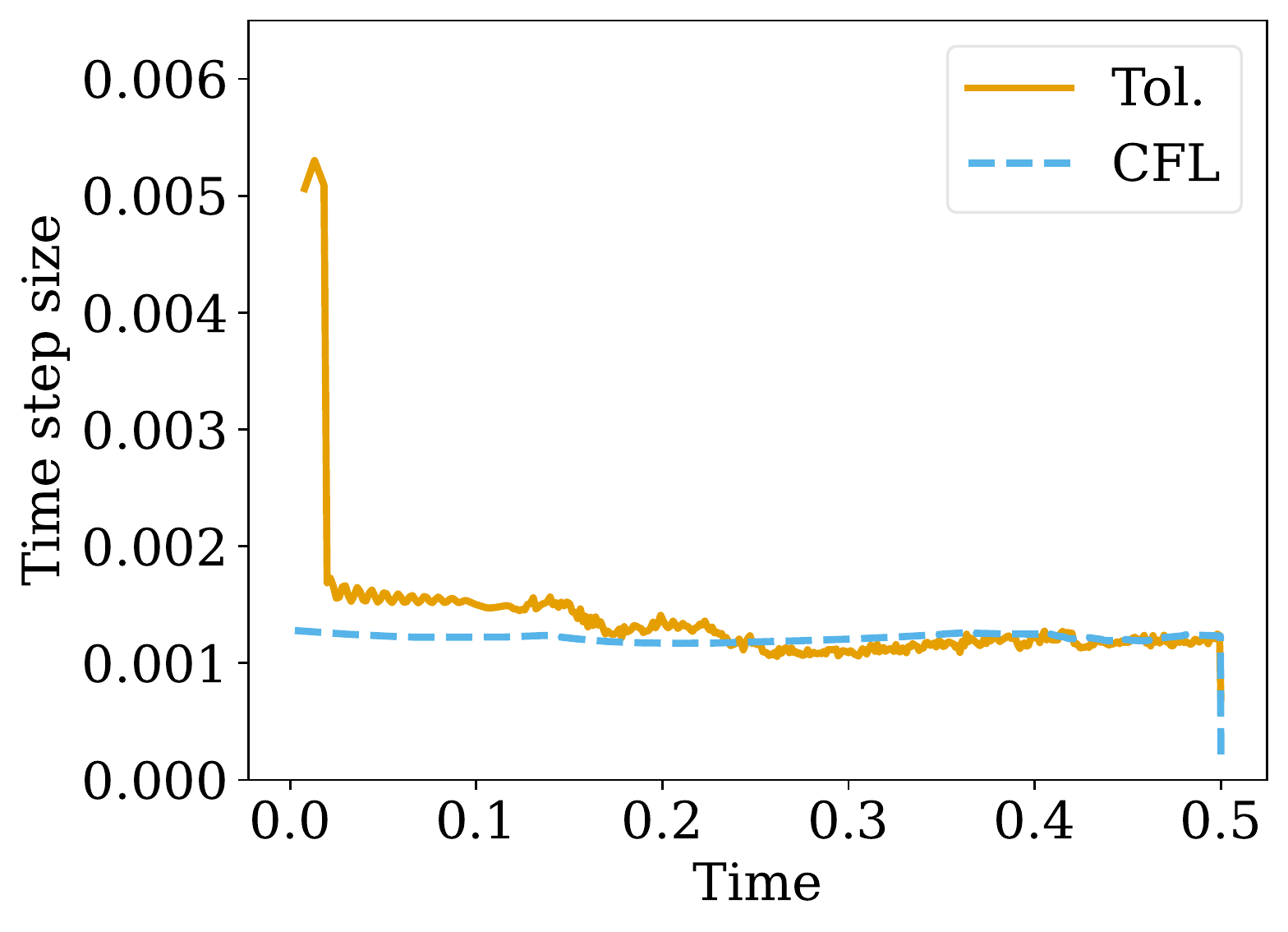}
    \caption{\ssp43.}
  \end{subfigure}%
  \caption{Effective CFL numbers (first row) and time step sizes (second row)
           of some representative RK methods applied to entropy-dissipative
           shock capturing semidiscretizations of the ideal GLM-MHD equations
           for the Orszag-Tang vortex.}
  \label{fig:orszag_tang}
\end{figure}

We recorded the time step sizes and effective CFL numbers, i.e., the CFL factor computed a
posteriori based on the time step size chosen by the error-based step size controller,
after every accepted
time step for some representative RK methods. The tolerance is set to $\tol = 10^{-4}$
in all cases with error-based step size control. For CFL-based step size control,
we bisected the maximal CFL number without blow-up to three significant digits.
The results are shown in Figure~\ref{fig:orszag_tang}.
Clearly, the time step size of the CFL-based approaches does not vary significantly
in time. Thus, the estimated CFL restriction are nearly constant. However, the
nonlinear schemes activate shock capturing mechanisms between $t = 0.1$ and $t = 0.2$.
Thus, the estimated linear CFL restriction should become more severe as discussed
in Section~\ref{sec:spectra}. This can be observed also for error-based step size
control; all RK methods use larger time steps initially until the shock capturing
part is activated. Then, error-based step size control yields approximately the
maximal step size also used by CFL-based step size control with manually tuned
CFL number. In particular, the initial time step size is approximately
one quarter bigger with error-based step size control (after the controller has
adapted the time step size from the automated initial guess).
During the full simulation, error-based step size control results in the
following performance improvements based on the number of function evaluations
(see also Table~\ref{tab:orszag_tang}):
\SI{13}{\percent} for \RK[BS]{3}{3}[][FSAL],
\SI{11}{\percent} for \RK[RDPK]{3}{5}[\ESstarp][FSAL], and
\SI{6}{\percent} for \ssp43.

While the quantitative numbers for
\RK[BS]{3}{3}[][FSAL] and \RK[RDPK]{3}{5}[\ESstarp][FSAL]
fit very well to the linear analysis, one could expect that the initial advantage
of \ssp43 with error-based step size control should be smaller, since the stability
region is less optimal for high-order DGSEM and relatively better suited for low-order
shock capturing methods. First, we do not necessarily expect that such a linear
analysis is quantitatively correct for general nonlinear problems, although it
was shown to behave very well in applications \cite{parsani2013optimized}.
Second, the shock capturing mechanism also changes the behavior of the ODE RHS,
usually reducing the smoothness of the spatially discrete system, which also
affects step size control.

\begin{table}[!htb]
\sisetup{output-exponent-marker=}
\sisetup{scientific-notation=fixed, fixed-exponent=0}
\centering
\caption{Performance of representative RK methods with default error-based and
         manually tuned CFL-based step size controllers:
         Number of function evaluations (\#FE), accepted steps (\#A),
         and rejected steps (\#R) for the Orszag-Tang vortex with entropy-dissipative
         shock capturing semidiscretizations.}
\label{tab:orszag_tang}
\setlength{\tabcolsep}{0.75ex}
\begin{tabular*}{\linewidth}{@{\extracolsep{\fill}}c *2c r@{\hskip 0.5ex}rr@{\hskip 0.5ex}r@{\hskip 1ex}cr@{\hskip 0.5ex}r@{\hskip 0.5ex}r} 
  \toprule
  Scheme & $\beta$ & $\tol$/$\cfl$ & \multicolumn{1}{c}{\#FE} & \multicolumn{1}{c}{\#A} & \multicolumn{1}{c}{\#R} \\ 
  \midrule

  \RK[BS]{3}{3}[][FSAL]            & $(0.60, -0.20, 0.00)$
     &  $\tol = 10^{-4}$ & $ 2187$ & $  724$ & $    4$ \\ 
   & &  $\cfl =  0.69  $ & $ 2506$ & $  835$ &         \smallskip\\ 
  \RK[RDPK]{3}{5}[\ESstarp][FSAL]  & $(0.70, -0.23, 0.00)$
     &  $\tol = 10^{-4}$ & $ 1863$ & $  368$ & $    4$ \\ 
   & &  $\cfl =  1.38  $ & $ 2091$ & $  418$ &         \smallskip\\ 
  \ssp43                           & $(0.55, -0.27, 0.05)$
     &  $\tol = 10^{-4}$ & $ 1550$ & $  384$ & $    3$ \\ 
   & &  $\cfl =  1.40  $ & $ 1648$ & $  412$ &         \\ 
  \bottomrule
\end{tabular*}
\end{table}

As expected, the method \RK[RDPK]{3}{5}[\ESstarp][FSAL] optimized for spectral
element discretizations performs better than the general-purpose method
\RK[BS]{3}{3}[][FSAL] and uses \SI{15}{\percent} less function evaluations
(with error-based step size control).
For this problem with a significant amount of shock capturing, \ssp43 with
error-based step size performs even better and uses additionally \SI{16}{\percent}
less function evaluations than \RK[RDPK]{3}{5}[\ESstarp][FSAL].

In summary, the results of this section demonstrate that error-based
step size control can react robustly and efficiently to varying linear CFL
restrictions in nonlinear schemes. In particular, they usually require less
parameter tuning while still achieving optimal performance. For this problem
with a significant amount of shock capturing mechanisms, the strong stability
preserving method \ssp43 performs well.

\section{Initial transient period: cold start of a simulation}
\label{sec:cold_start}

When performing a cold start of a demanding CFD simulation, \eg, by initializing
the flow around objects with free stream values, the initial transient period
often requires smaller time steps than the fully developed simulation.
For demanding simulations with positivity issues for high-order spatial
schemes used in this section, we apply SSP time integration methods.

\subsection{Double Mach reflection of a strong shock}
\label{sec:double_mach_reflection}

\begin{figure}[!htb]
\centering
  \begin{subfigure}{\linewidth}
    \includegraphics[width=\textwidth]{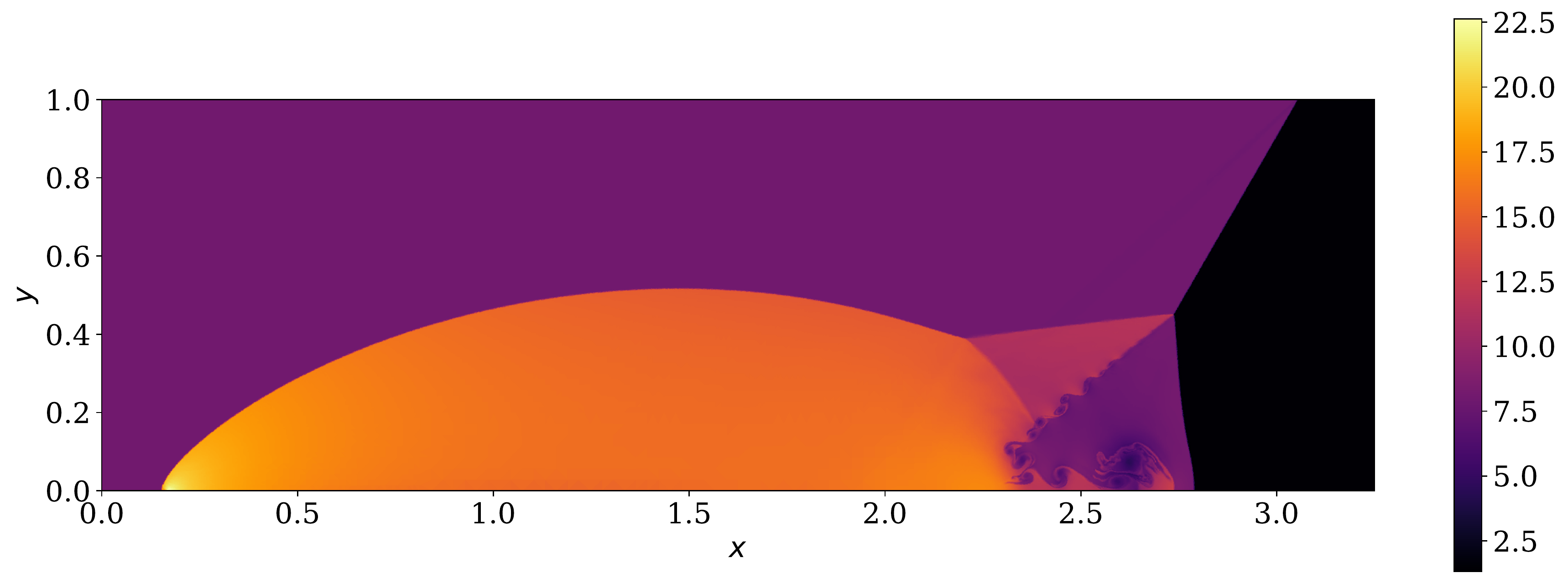}
    \caption{Density.}
  \end{subfigure}%
  \\
  \begin{subfigure}{\linewidth}
    \includegraphics[width=\textwidth]{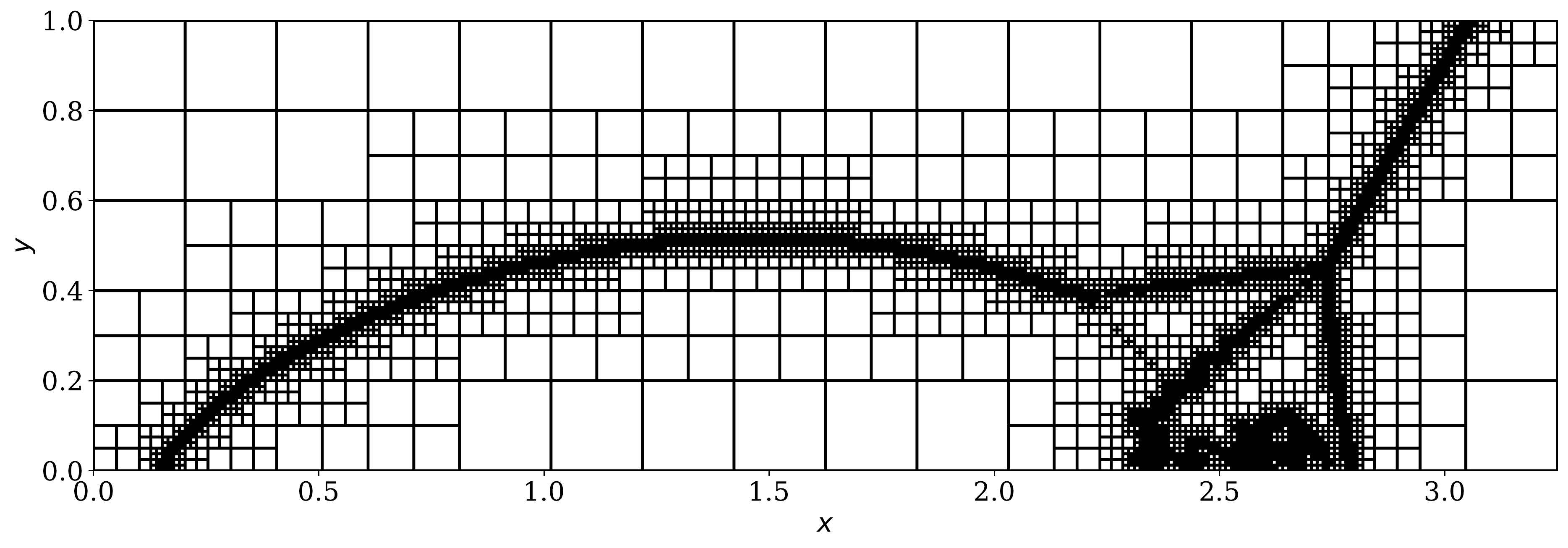}
    \caption{Mesh.}
  \end{subfigure}%
  \caption{Density and mesh at the final time of the simulation of the
           Woodward-Colella double Mach reflection using \ssp43 applied to
           entropy-dissipative shock capturing semidiscretizations with
           positivity-preserving limiters and adaptive mesh refinement.}
  \label{fig:double_mach_reflection_snapshot}
\end{figure}

First, we consider the double Mach reflection of Woodward and Colella \cite{woodward1984numerical}.
We set up an initial grid with 80 uniform elements that is adaptively refined
during the simulation using \texttt{p4est} \cite{burstedde2011p4est}. We use
the shock capturing approach of \cite{hennemann2021provably} with DG elements
using polynomials of degree $\polydeg = 4$ and the entropy-conservative numerical
flux of \cite{ranocha2020entropy,ranocha2018thesis,ranocha2021preventing} in the
volume terms and a local Lax-Friedrichs flux at interfaces and finite volume
subcells. AMR is triggered every two time steps and the positivity-preserving
limiter of \cite{zhang2011maximum} for density and pressure is applied after
every Runge-Kutta stage. We apply \ssp43 to integrate the system in time
$t \in [0, 0.2]$. The complete setup can be found in the reproducibility
repository \cite{ranocha2022errorRepro}.
Figure~\ref{fig:double_mach_reflection_snapshot} shows the density and the mesh
at the final time.

\begin{figure}[!htb]
\centering
  \begin{subfigure}{0.49\linewidth}
    \includegraphics[width=\textwidth]{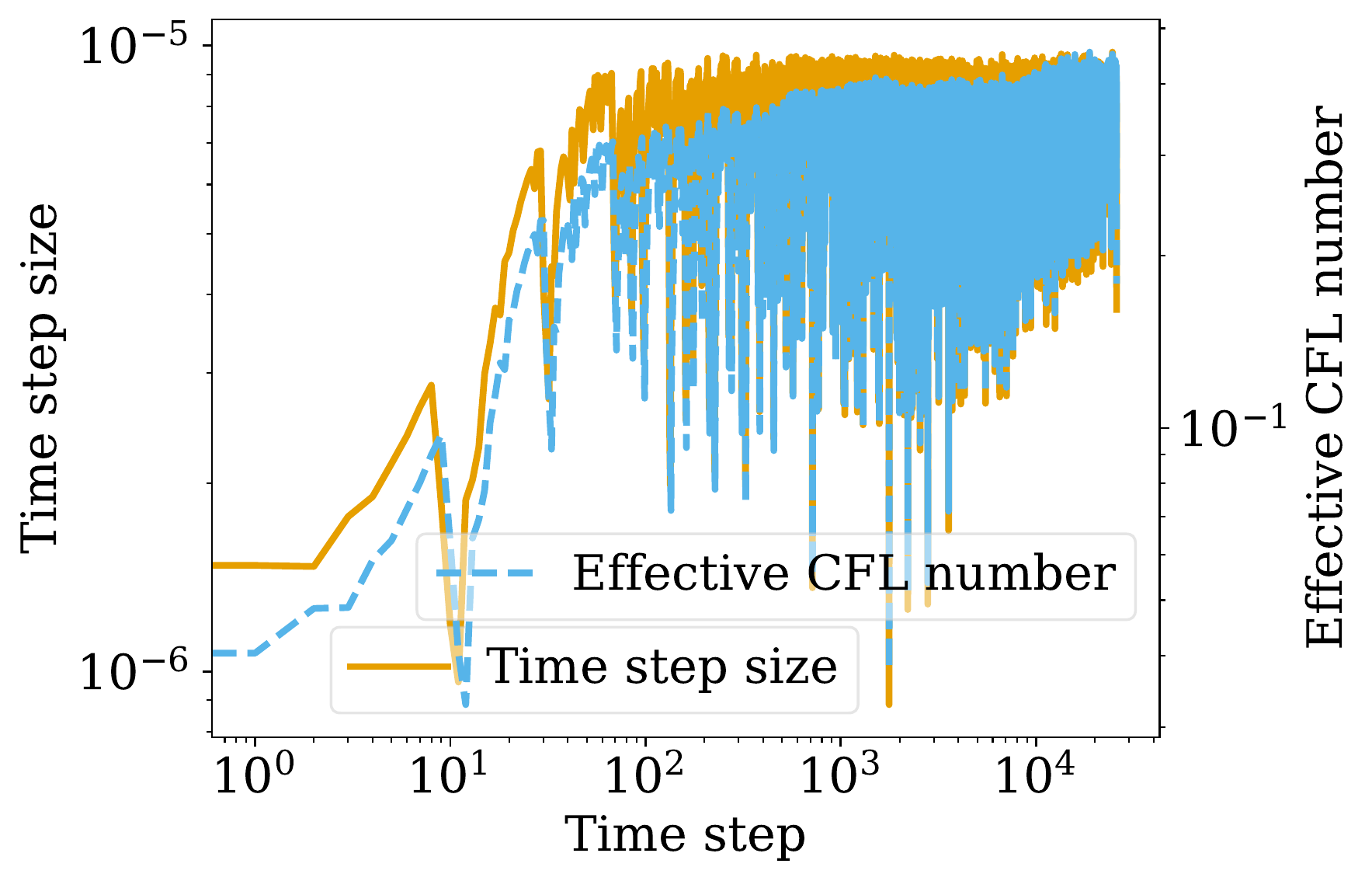}
    \caption{Dependence on time step number.}
  \end{subfigure}%
  \hspace*{\fill}
  \begin{subfigure}{0.49\linewidth}
    \includegraphics[width=\textwidth]{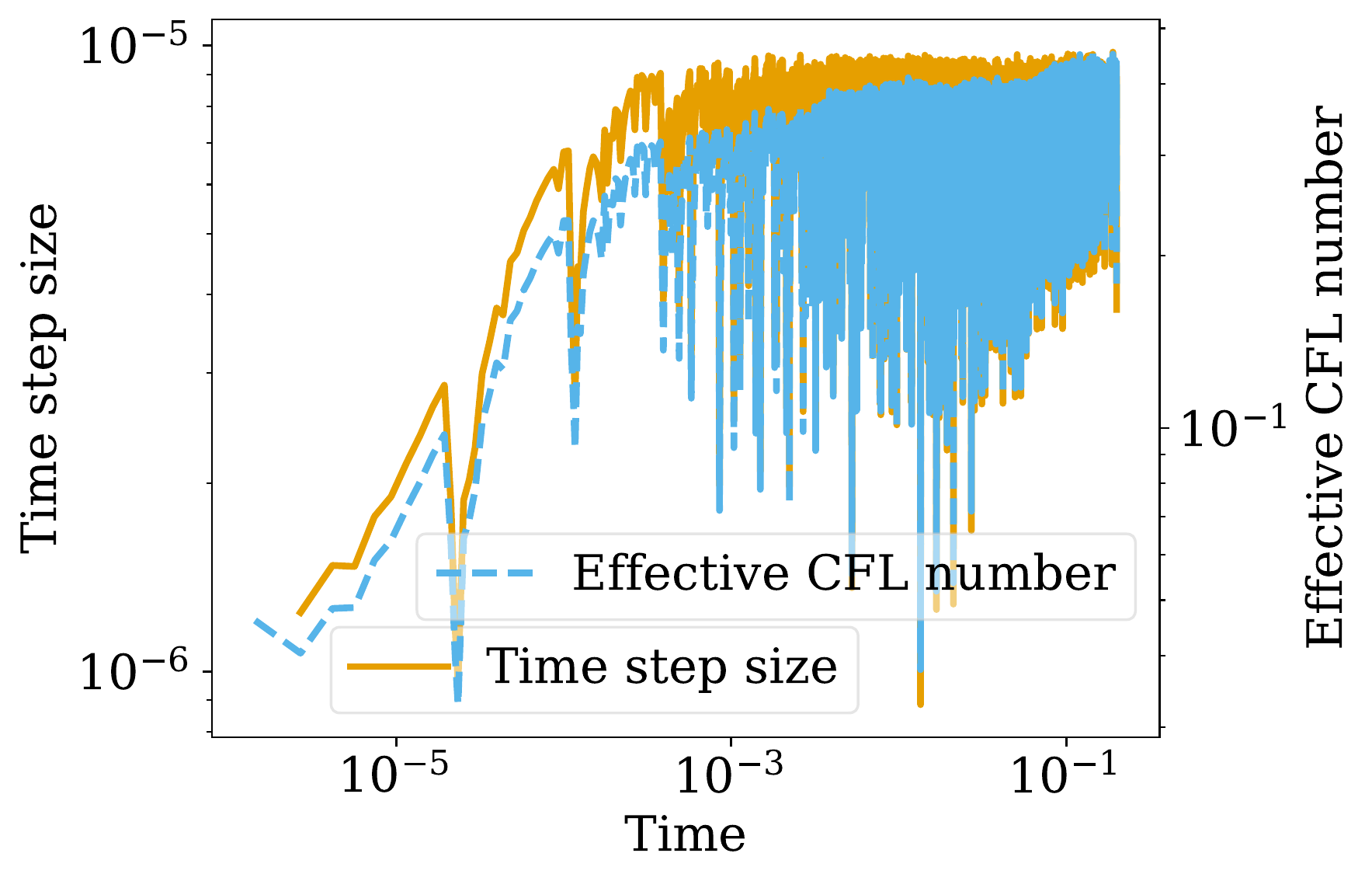}
    \caption{Dependence on time.}
  \end{subfigure}%
  \caption{Time step sizes and effective CFL numbers of \ssp43 applied to
           entropy-dissipative shock capturing semidiscretizations with
           positivity-preserving limiters and adaptive mesh refinement for the
           double Mach reflection problem.}
  \label{fig:double_mach_reflection}
\end{figure}

Figure~\ref{fig:double_mach_reflection} shows the evolution of the time step
sizes and the effective CFL number. There is an initial transient period of \ca
100 time steps where both the time step size and the effective CFL number are
much smaller than in the remaining simulation. This small initial step size is
selected by the automatic detection algorithm and the error-based step size
control, no manual intervention is necessary. The same simulation setup run
with CFL-based step size control crashes in the first few time steps, even
with a small CFL number of $\cfl = 0.01$. The initial transient period can be
captured with even smaller CFL numbers, but running the simulation afterwards
in a reasonable amount of time required adapting the CFL factor. In contrast,
error-based step size control can be used with the same parameters throughout
the simulation.

\subsection{Astrophysical Mach 2000 jet}
\label{sec:astro_jet}

Next, we consider a simulation of an astrophysical Mach 2000 jet using the ideal
compressible Euler equations with ratio of specific heat $\gamma = 5 / 3$
based on the setup of \cite{liu2022essentially}. Specifically, we use
the initial condition
\begin{equation}
  \rho^0 = 0.5, \quad v_1^0 = 0, \quad v_2^0 = 0, \quad p^0 = 0.4172
\end{equation}
in the domain $[-0.5, 0.5]$ equipped with periodic boundary conditions in $y$-direction
and Dirichlet boundary conditions using the initial data in the $x$-direction
except for $t > 0$, $x = -0.5$, and $y \in [-0.05, 0.05]$, where we use
boundary data (denoted by a superscript $b$)
\begin{equation}\label{eq:mach2000_bcs}
  \rho^b = 5, \quad v_1^b = 800, \quad v_2^b = 0, \quad p^b = 0.4172.
\end{equation}
We use an initially uniform mesh with $64$ elements per coordinate direction and
polynomials of degree $\polydeg = 3$. We use the shock capturing approach of
\cite{hennemann2021provably} with the entropy-conservative numerical flux of
\cite{ranocha2020entropy,ranocha2018thesis,ranocha2021preventing} in the volume
terms and a local Lax-Friedrichs flux at interfaces and finite volume subcells.
After every time step, we apply adaptive mesh refinement. Moreover, we
apply the positivity-preserving limiter of \cite{zhang2011maximum} for density
and pressure after every Runge-Kutta stage. We integrate the semidiscrete system
in time using \ssp43 for $t \in [0, 10^{-3}]$. The complete setup can be found
in the reproducibility repository \cite{ranocha2022errorRepro}.

\begin{figure}[!htb]
\centering
  \begin{subfigure}{0.49\linewidth}
    \includegraphics[width=\textwidth]{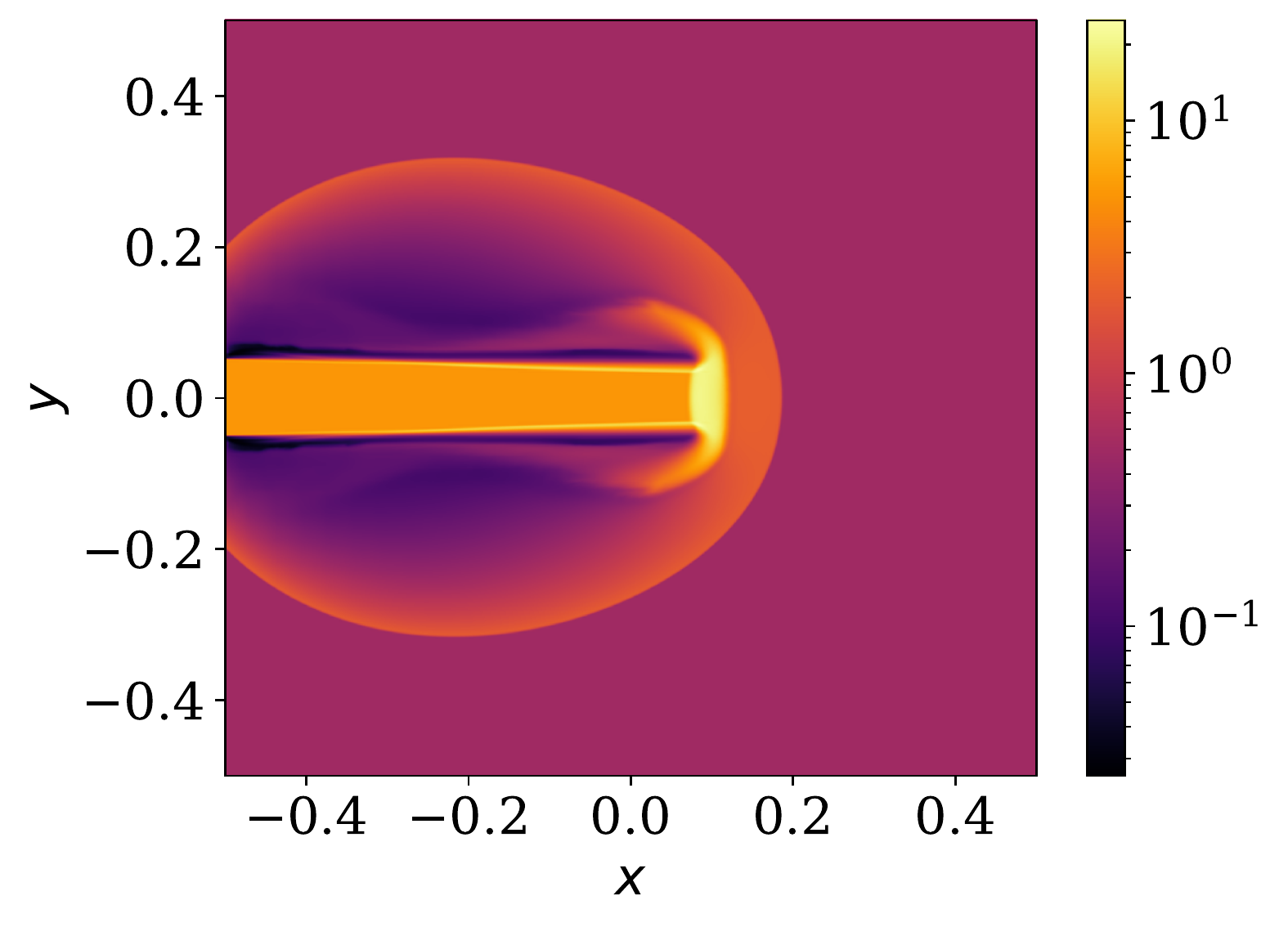}
    \caption{Density.}
  \end{subfigure}%
  \hspace*{\fill}
  \begin{subfigure}{0.49\linewidth}
    \includegraphics[width=\textwidth]{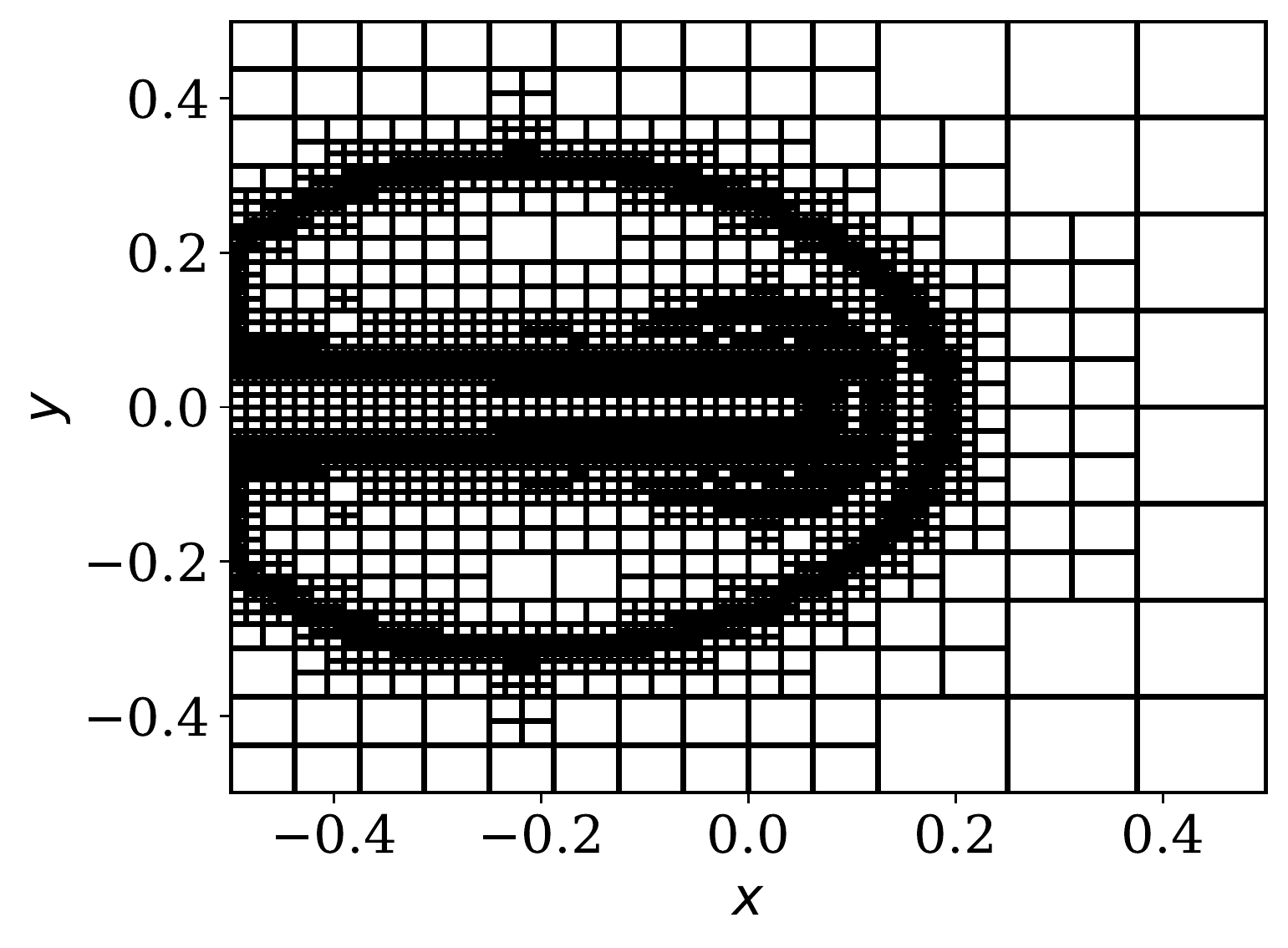}
    \caption{Mesh.}
  \end{subfigure}%
  \caption{Density and mesh at the final time of the simulation of an
           astrophysical Mach 2000 jet using \ssp43 applied to
           entropy-dissipative shock capturing semidiscretizations with
           positivity-preserving limiters and adaptive mesh refinement.}
  \label{fig:astro_jet_snapshot}
\end{figure}

Figure~\ref{fig:astro_jet_snapshot} shows a snapshot of the numerical solution
and the adapted grid at the final time $t = 10^{-3}$. While the initial condition
is given by uniform free stream values, the boundary condition \eqref{eq:mach2000_bcs}
changes the
numerical approximation rapidly in time. Consequently, the mesh is refined over
time and adapts to the solution structures.

\begin{figure}[!htb]
\centering
  \begin{subfigure}{0.49\linewidth}
    \includegraphics[width=\textwidth]{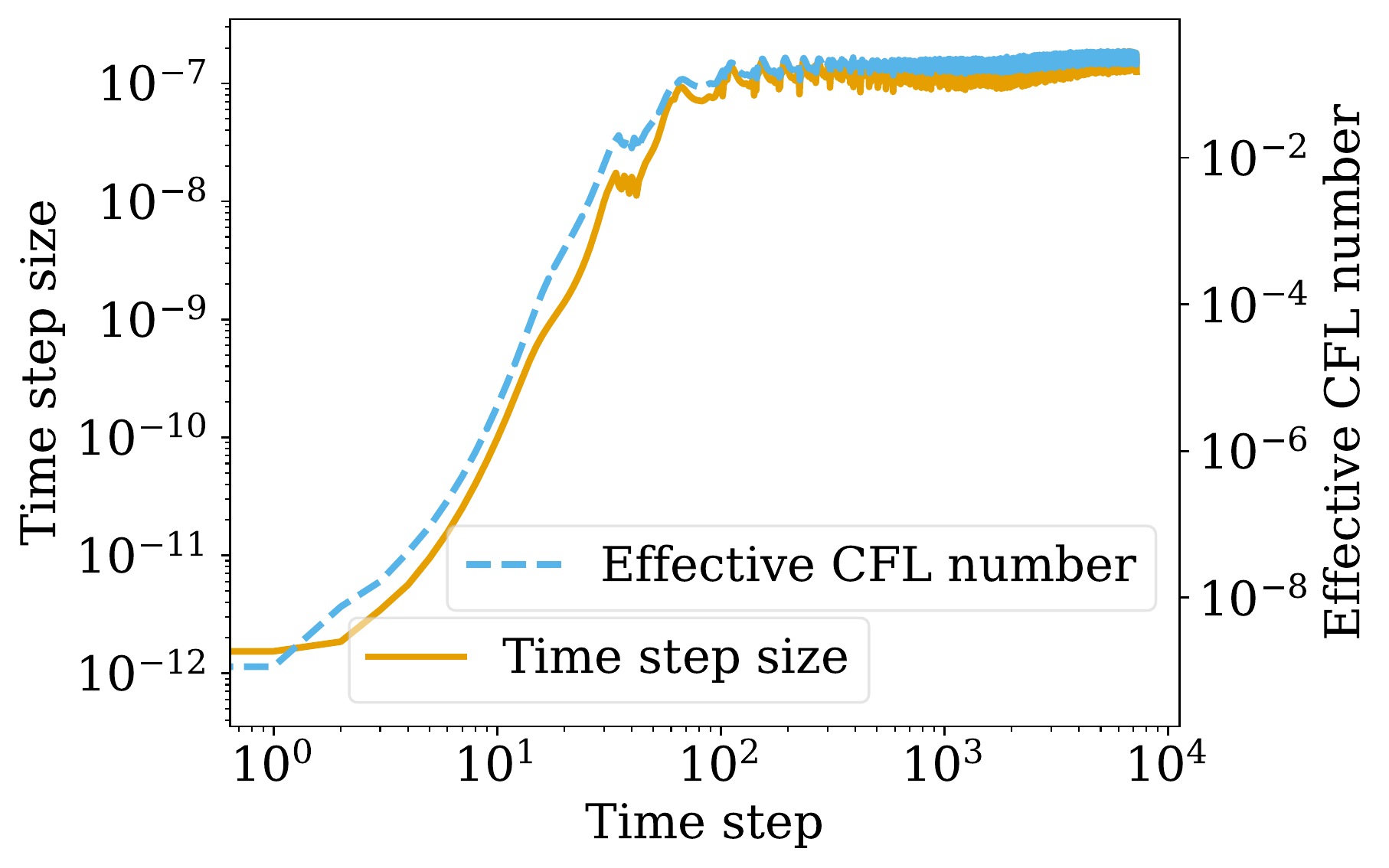}
    \caption{Dependence on time step number.}
  \end{subfigure}%
  \hspace*{\fill}
  \begin{subfigure}{0.49\linewidth}
    \includegraphics[width=\textwidth]{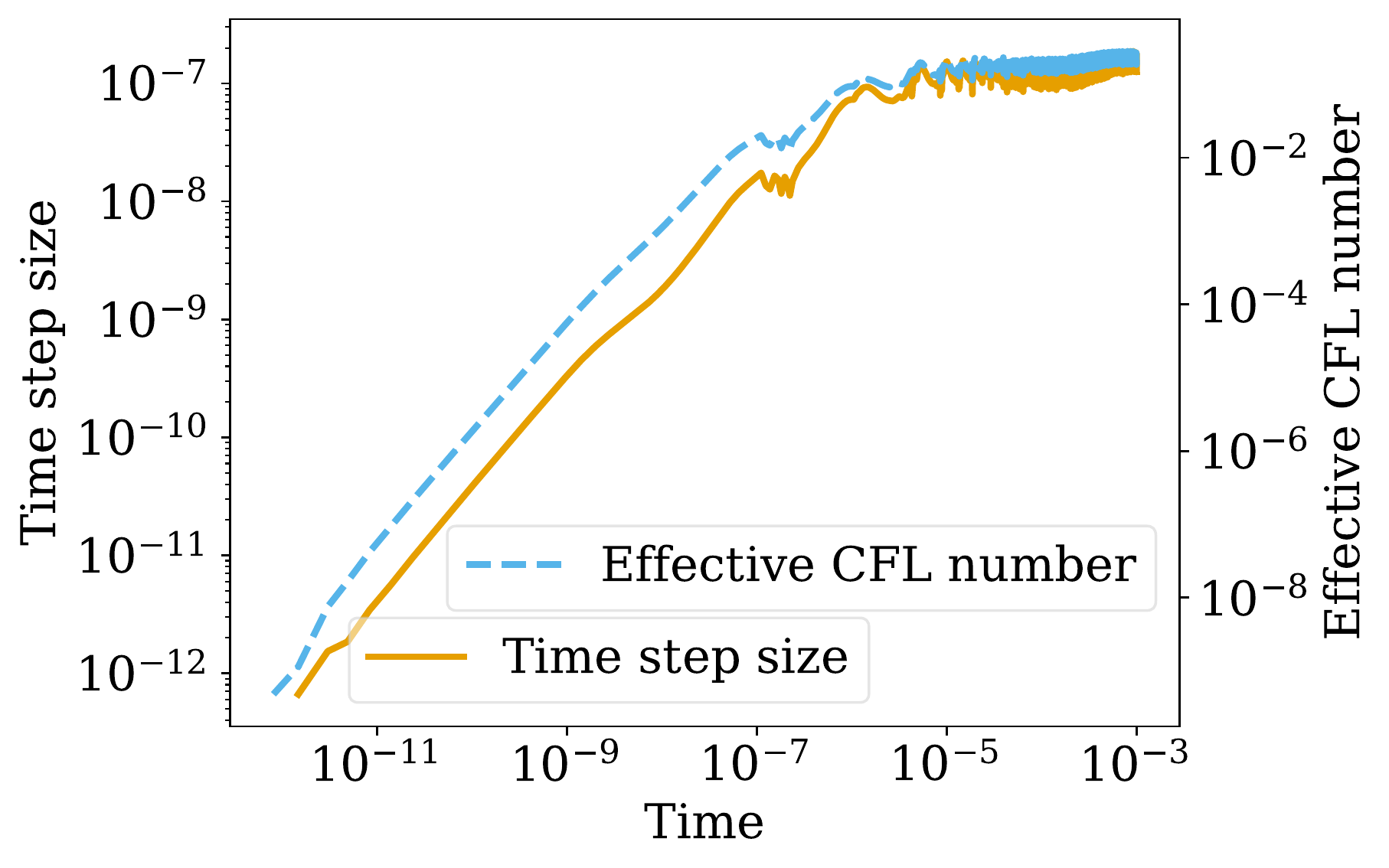}
    \caption{Dependence on time.}
  \end{subfigure}%
  \caption{Time step sizes and effective CFL numbers of \ssp43 applied to
           entropy-dissipative shock capturing semidiscretizations with
           positivity-preserving limiters and adaptive mesh refinement for the
           astrophysical Mach 2000 jet.}
  \label{fig:astro_jet}
\end{figure}

Figure~\ref{fig:astro_jet} shows the evolution of the time step for this
demanding simulation. The initial time step must be very small --- of the order
$\dt \approx 10^{-12}$ --- to allow the simulation to run. For this setup, we
set $\dt = 10^{-12}$ initially; other choices are also possible and are adapted
to similar values automatically by the controller, but the fully automatic
initial guess is not sufficient in this case.

The initial transient period lasts for \ca 100 time steps where the error-based
step size controller increases the time step size gradually. Afterwards, the time
step size plateaus. Handling such an initial transient period is difficult for
standard CFL-based step size controllers. Initially, a CFL factor of
$\nu \approx 10^{-9}$ is required; $\cfl = 10^{-8}$ results in a blow-up of the
simulation. One could implement a more complicated CFL-based controller adapting
the CFL factor in the initial transient period or restart the simulation after
\ca 100 time steps with a larger CFL factor. However, no such additional techniques
are required for error-based step size control, making it easier to use in this
case.

\subsection{Delta wing}
\label{sec:delta_wing}

\begin{figure}[!ht]
\centering
\includegraphics[width=0.90\columnwidth]{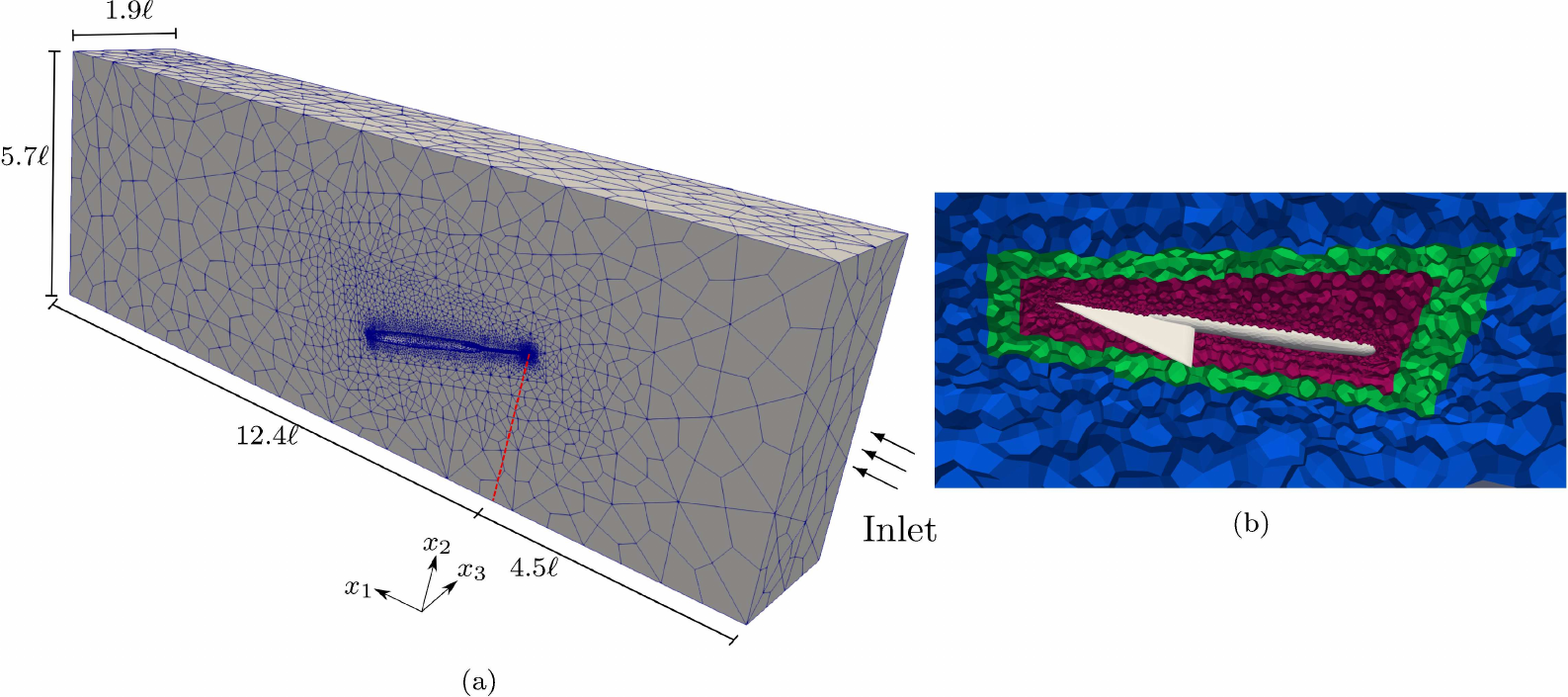}
\caption{Solution polynomial degree distribution, $p$, computational domain and boundary mesh elements
for the \ang{65} swept delta wing test case \cite{parsani2021ssdc}.}
\label{fig:delta_deg_distrib}
\end{figure}

The next test case setup approximates the solution of the compressible Navier-Stokes equations.
We consider an industrially relevant simulation of a \ang{65} swept delta wing at
a Mach number $\mathrm{Ma} = 0.9$, a Reynolds number $\mathrm{Re}=10^{6}$
(based on the mean aerodynamic chord), and setting the angle of attack to $\mathrm{AoA}= \ang{13}$.
The geometry was proposed by Hummel and Redeker \cite{hummel2003vortex}
for the Second International Vortex Flow Experiment (VFE-2) and manufactured based on a NASA wing geometry that served as
reference configuration \cite{chu1996delta65}.
We consider the flow past the leading edge configuration with a medium radius leading edge $r_{LE}/\bar{c}=0.0015$,
where $\bar{c}=\SI{0.653}{m}$. The delta wing has a mean aerodynamic chord of $\ell=\SI{0.667}{m}$, a root chord length of
$c_r = 1.47 \ell$, and a wing span of $b=1.37 \ell$. The sting present in the
wind tunnel testing is kept downstream as part of the geometry up to the
position $x_1/c_r=1.758$, where the $x_1$ Cartesian coordinate points in the
streamwise direction as shown in the left part of Figure~\ref{fig:delta_deg_distrib}.

For this test, we use the $hp$-adaptive entropy stable solver SSDC \cite{parsani2021ssdc}, which is able to deliver
numerical results in good agreement with experimental data.
Entropy stable adiabatic no-slip wall boundary conditions \cite{dalcin2019conservative}
are applied to the wing and
sting surfaces, whereas freestream far-field boundary conditions are applied
at the inlet and outlet planes.
Due to the symmetry of the problem in the span-wise direction, half span of the flow is modeled
through a symmetry boundary condition.
On the rest of the boundary planes, entropy stable inviscid wall boundary conditions
\cite{parsani2015entropy} are prescribed.
The grid consists of $\approx 9.2\times 10^4$ hexahedral elements.
It is subdivided into three blocks, as shown in the right part of Figure~\ref{fig:delta_deg_distrib}, with each block
corresponding to a different solution polynomial degree $p$.
In particular, we use $p=2$ in the far-field region, $p=5$ in the region surrounding the delta
wing and its support, and $p=3$ elsewhere.
Given the degree of the solution and the number of elements in each block, the total number of
degrees of freedom (DOFs) is
$\approx 1.435 \times 10^{7}$.

\begin{figure}[!ht]
\centering
\includegraphics[width=0.7\columnwidth]{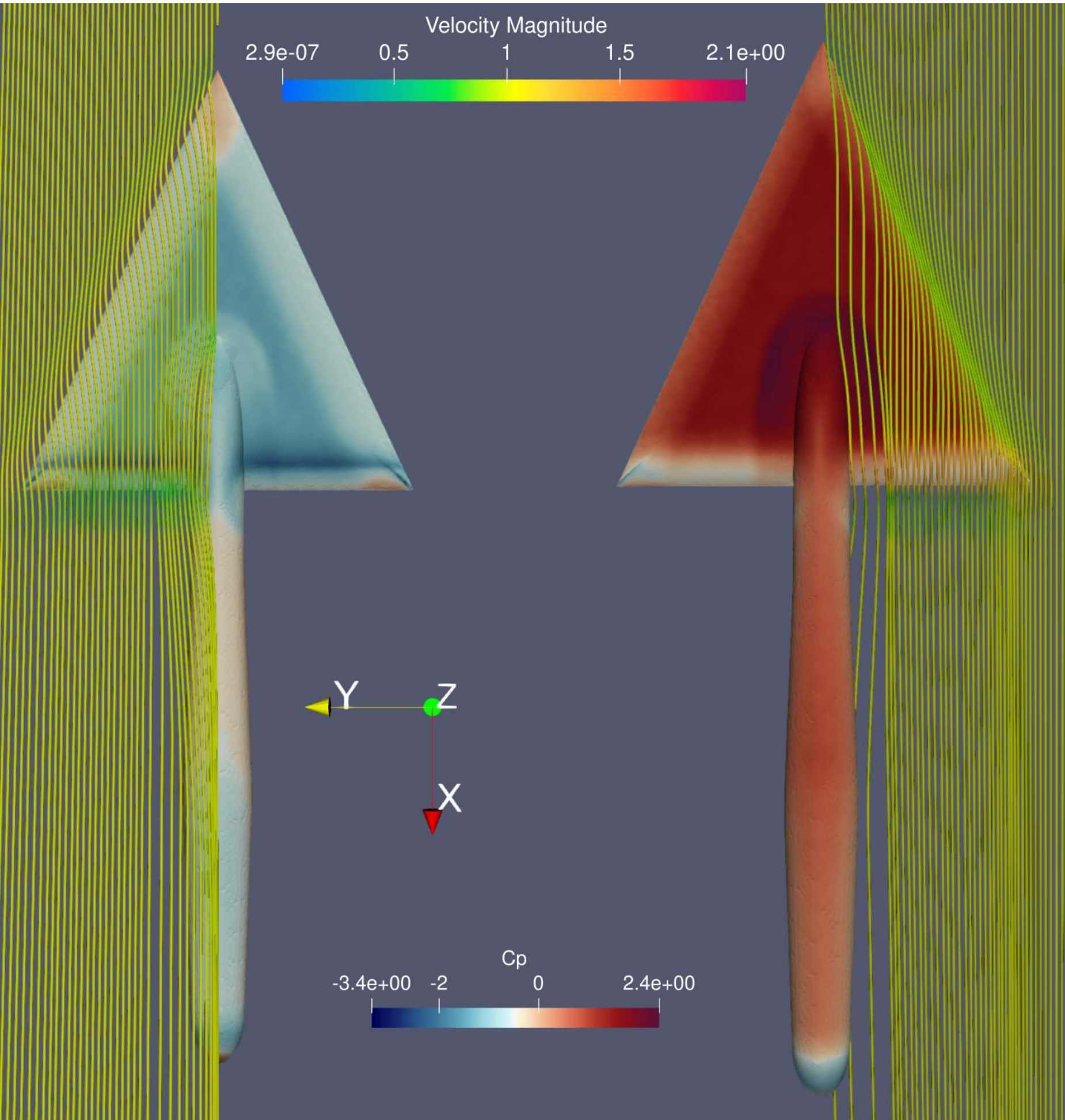}
\caption{Instantaneous pressure coefficient distribution and streamlines colored
by the velocity magnitude at $t = 0.05$
for the \ang{65} swept delta wing test case \cite{parsani2021ssdc} (left: top view, right: bottom
view).}
\label{fig:delta_cp}
\end{figure}

The flow around a delta wing is peculiar.
When the angle of attack exceeds \ang{7}, typically, flow separation occurs at the
leading edge. The roll up of the leading-edge vortices induces low pressure
on the upper surface of the wing and enhances the lift. Figure
\ref{fig:delta_cp} shows the instantaneous pressure coefficient distribution and streamlines colored
by the velocity magnitude at $t = 0.05$.

\begin{figure}[!htb]
\centering
  \begin{subfigure}{0.49\linewidth}
    \includegraphics[width=\textwidth]{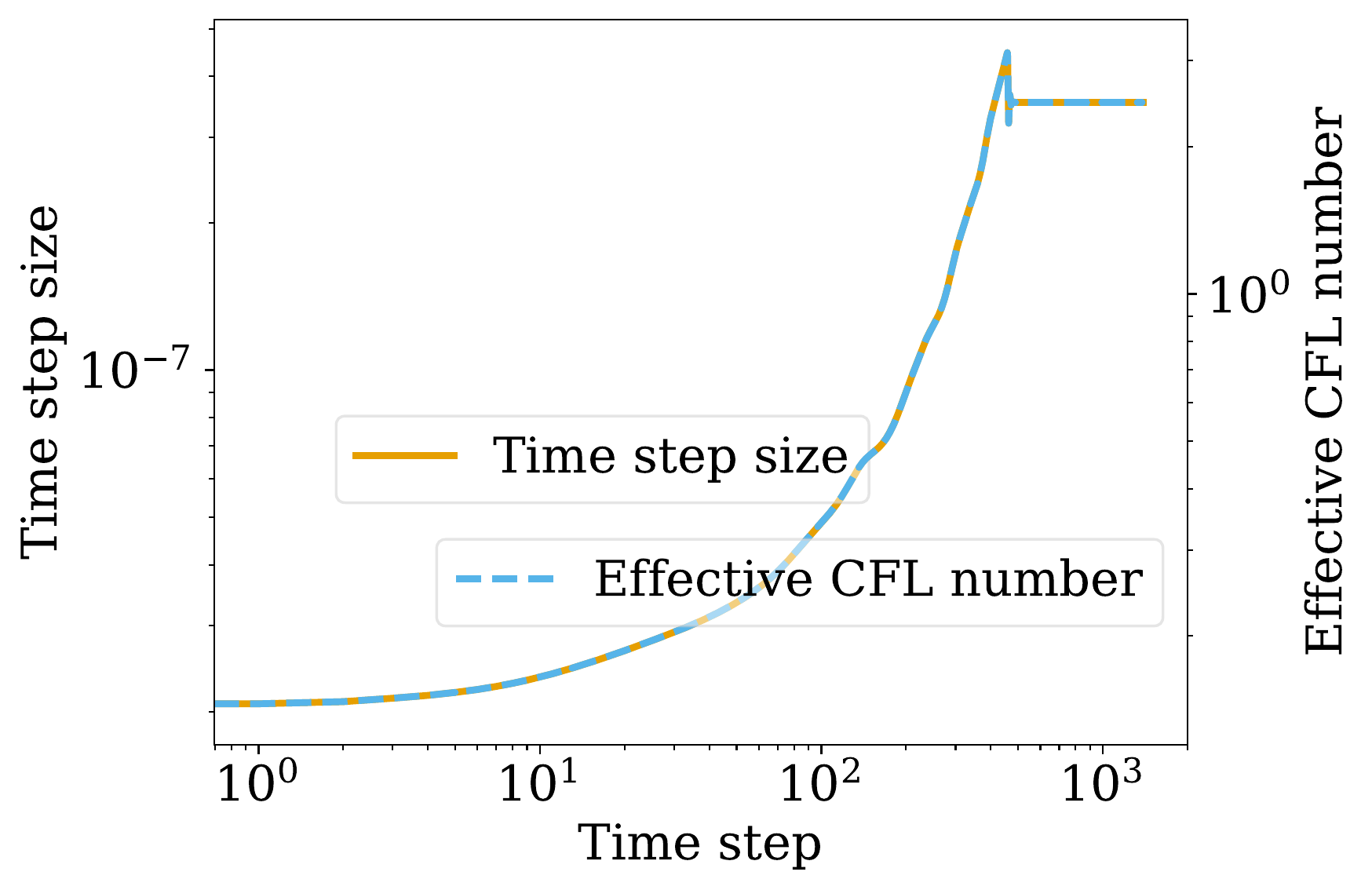}
    \caption{\RK[BS]{3}{3}[][FSAL] (time step number).}
  \end{subfigure}%
  \hspace*{\fill}
  \begin{subfigure}{0.49\linewidth}
    \includegraphics[width=\textwidth]{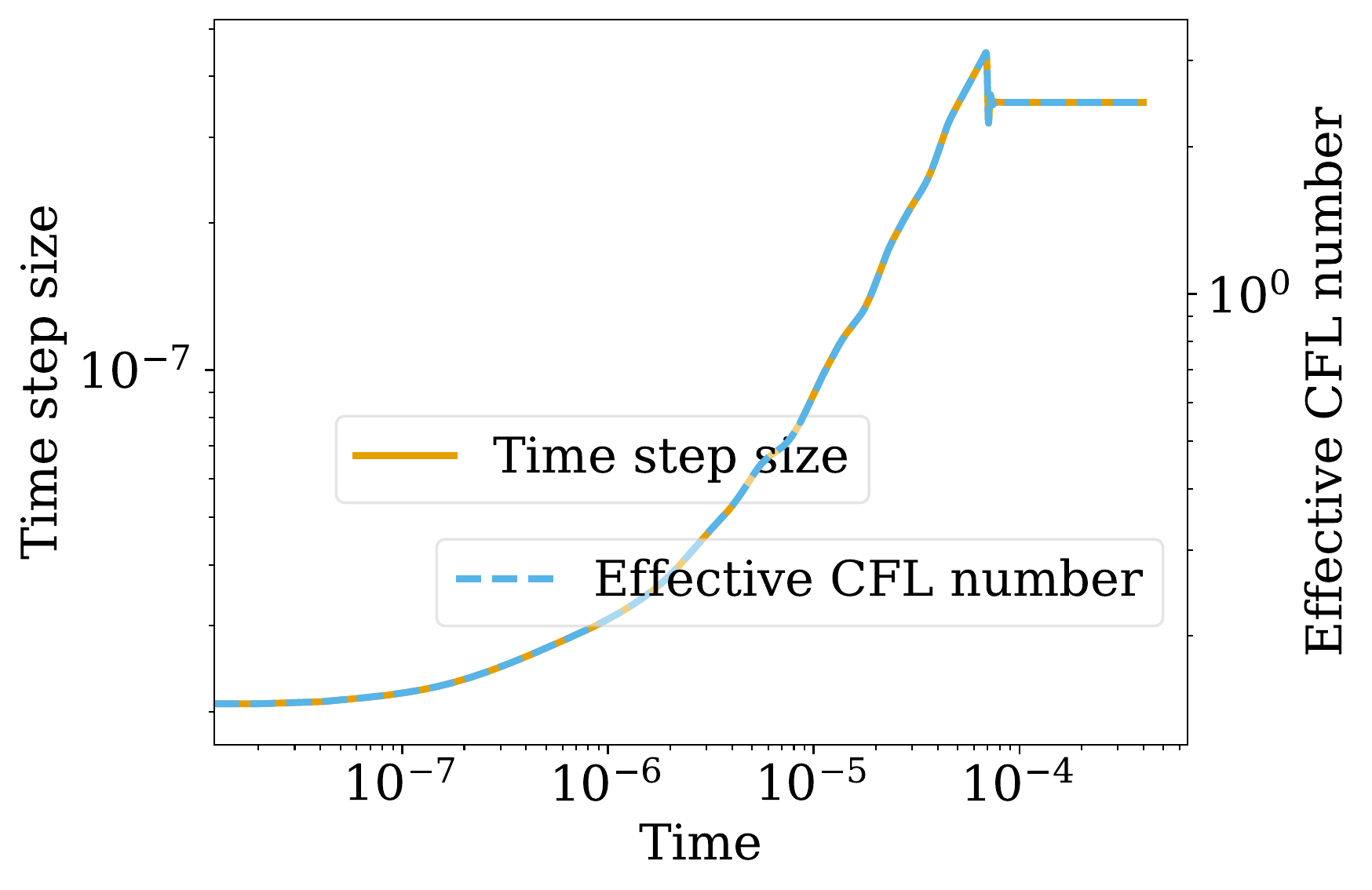}
    \caption{\RK[BS]{3}{3}[][FSAL] (physical time).}
  \end{subfigure}%
  \\
  \begin{subfigure}{0.49\linewidth}
    \includegraphics[width=\textwidth]{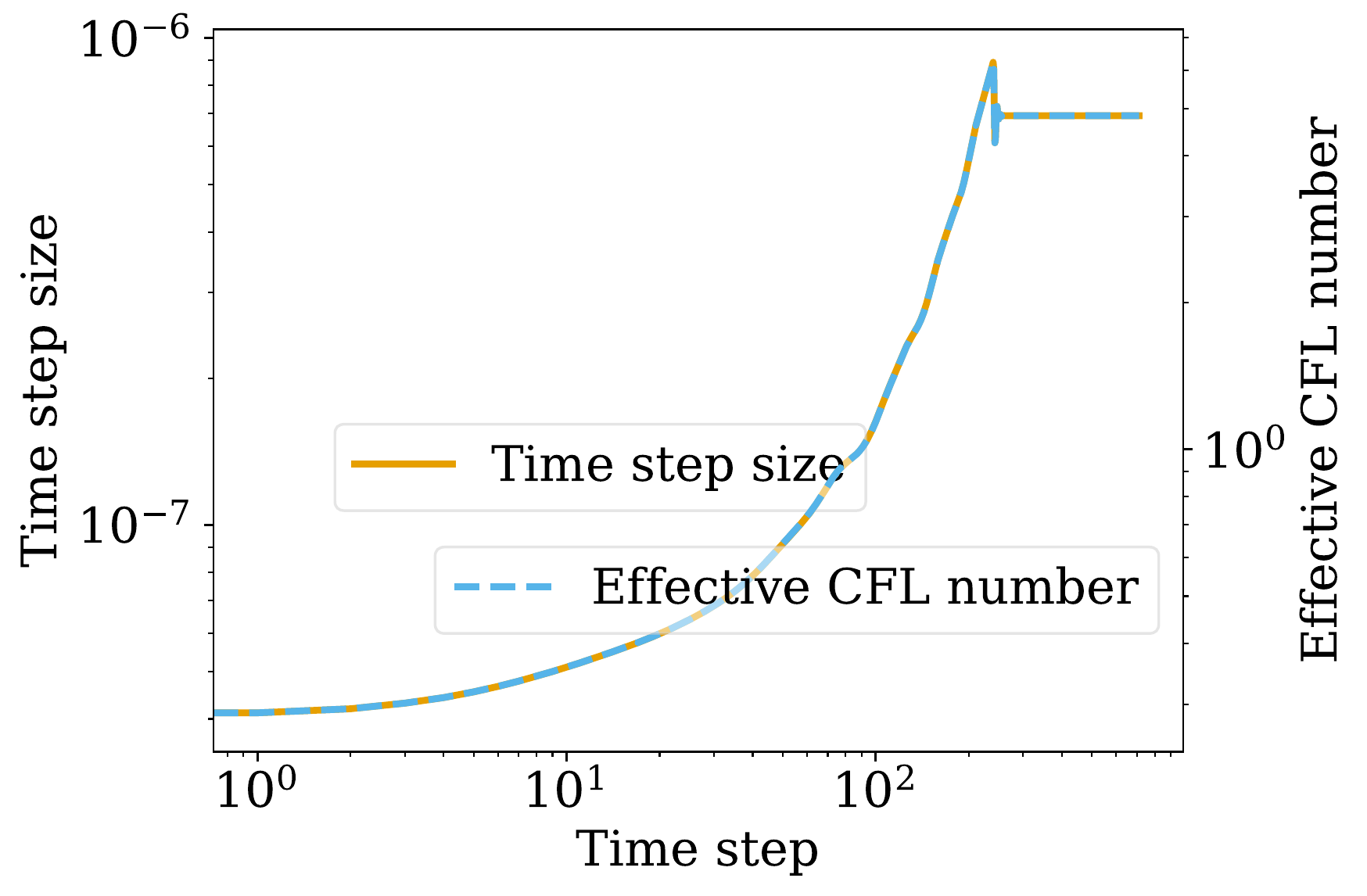}
    \caption{\RK[RDPK]{3}{5}[\ESstarp][FSAL] (time step number).}
  \end{subfigure}%
  \hspace*{\fill}
  \begin{subfigure}{0.49\linewidth}
    \includegraphics[width=\textwidth]{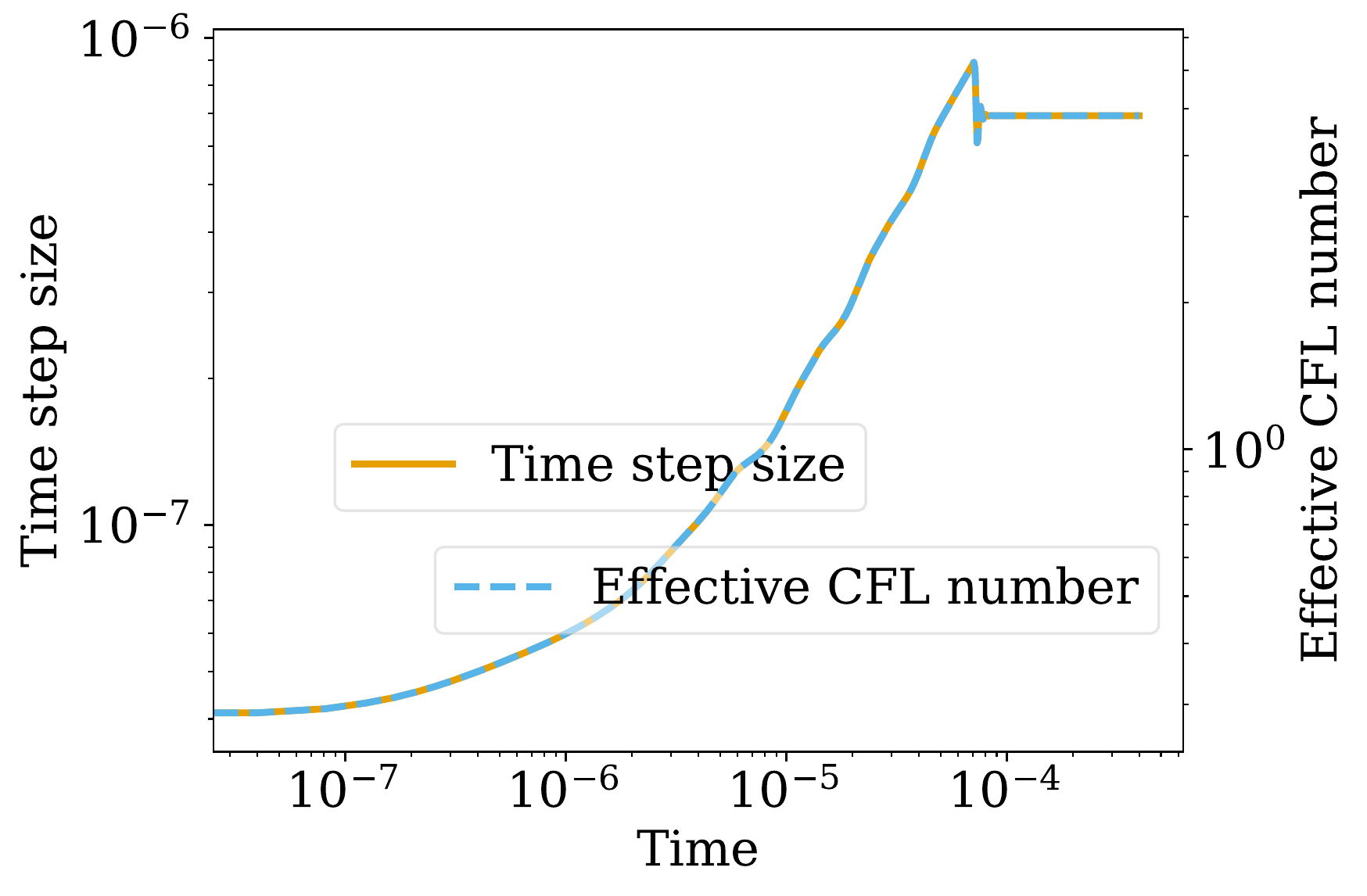}
    \caption{\RK[RDPK]{3}{5}[\ESstarp][FSAL] (physical time).}
  \end{subfigure}%
  \\
  \begin{subfigure}{0.49\linewidth}
    \includegraphics[width=\textwidth]{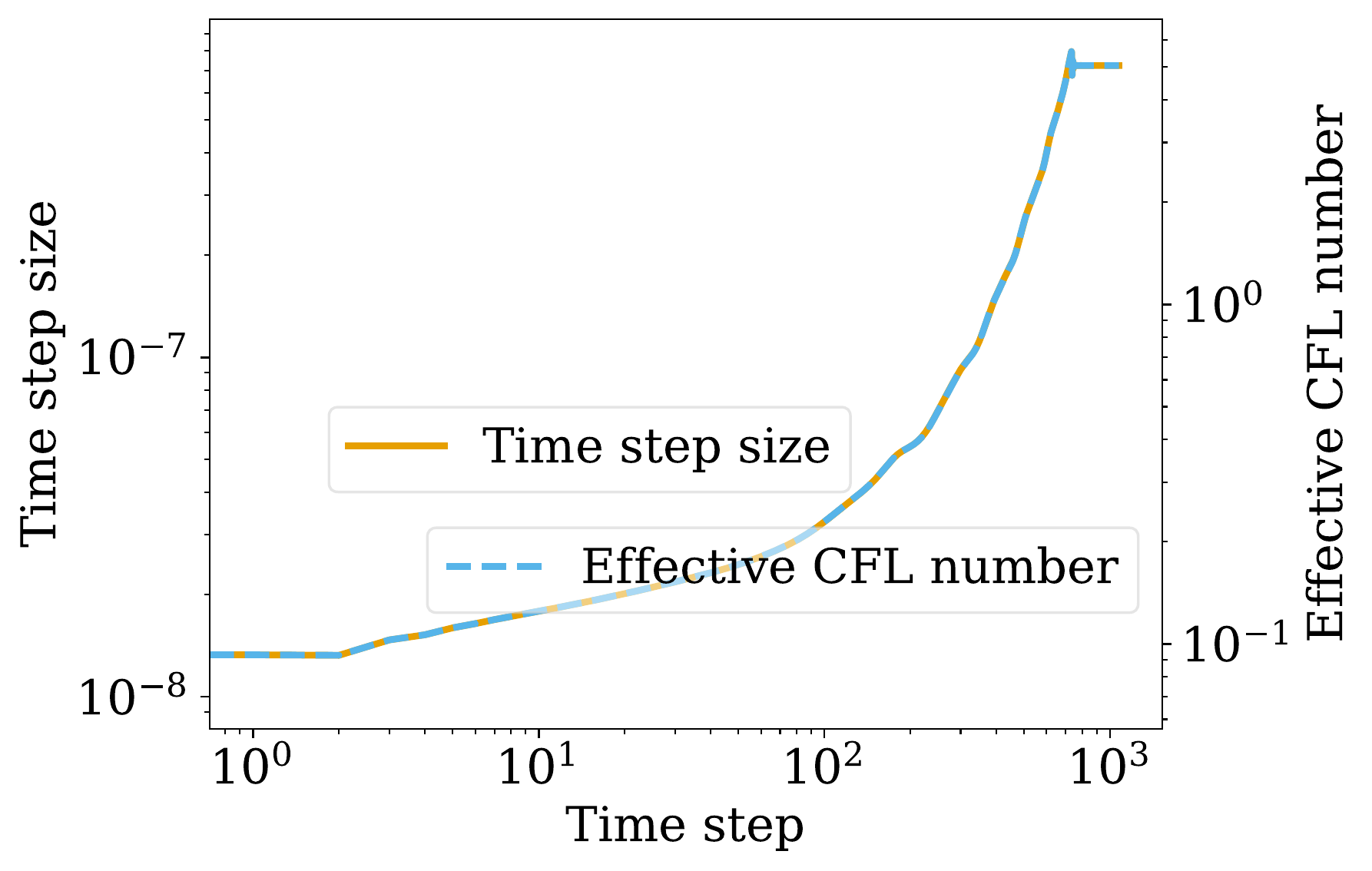}
    \caption{\ssp43 (time step number).}
  \end{subfigure}%
  \hspace*{\fill}
  \begin{subfigure}{0.49\linewidth}
    \includegraphics[width=\textwidth]{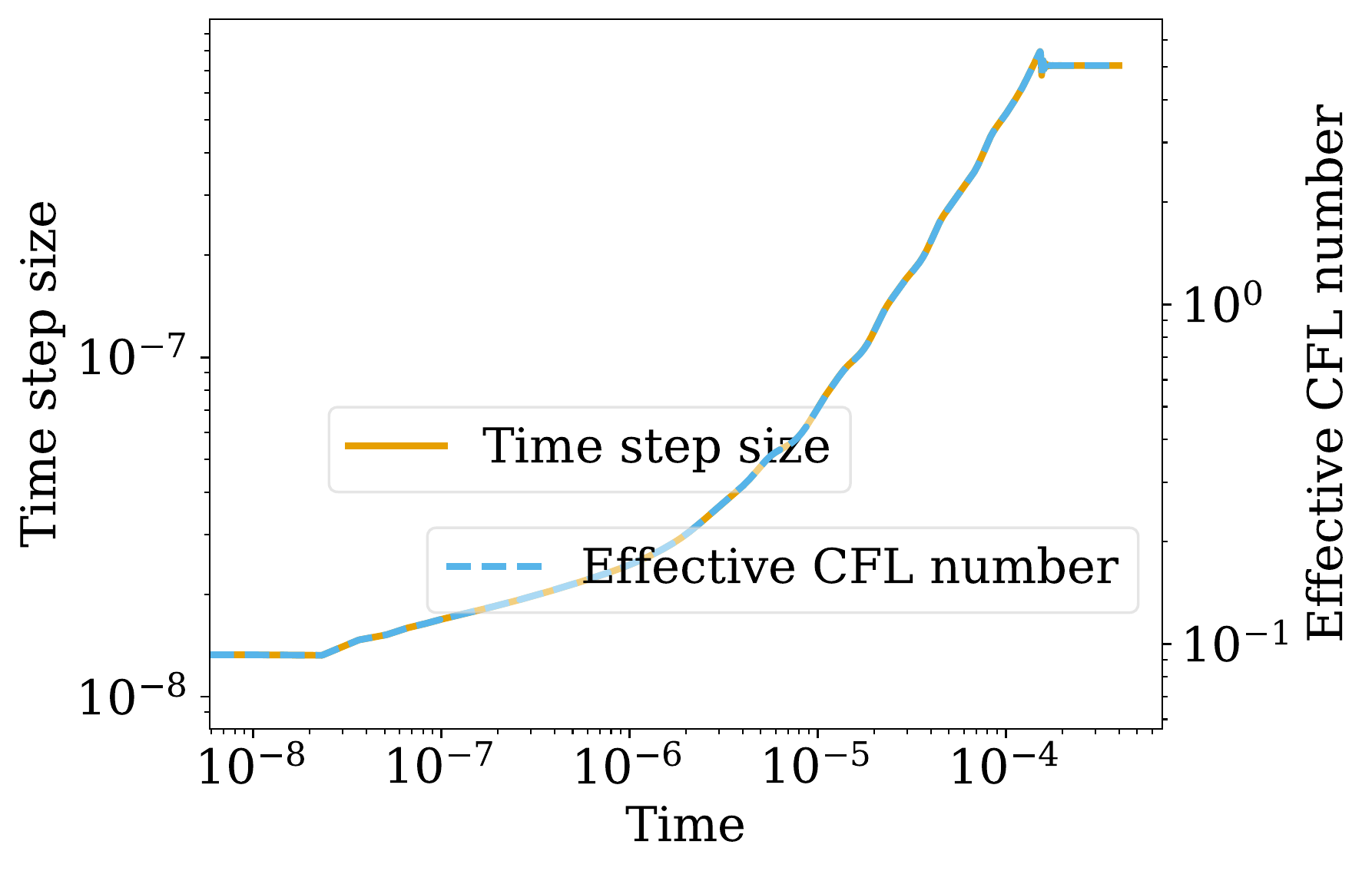}
    \caption{\ssp43 (physical time).}
  \end{subfigure}%
  \caption{Effective CFL numbers and time step sizes depending on the
           time step index (left column) and physical time (right column)
           of some representative RK methods applied to entropy-dissipative
           semidiscretizations of the compressible Navier-Stokes equations for the delta wing
           setup.}
  \label{fig:delta_wing}
\end{figure}

Results of simulations using some representative RK methods are shown in
Figure~\ref{fig:delta_wing}. For this cold start of the simulation with free
stream values, all methods begin with a conservative estimate of the time step
size. The error-based step size controllers increase the time step size in the
first few time steps and reach an asymptotically constant step size after
a few hundred time steps.

The effective CFL number shown in Figure~\ref{fig:delta_wing} is basically
proportional to the time step size itself. Running SSDC with CFL-based step
size control requires a CFL factor that is (at least initially)
\begin{itemize}
  \item \SI{13}{\percent} smaller for \RK[BS]{3}{3}[][FSAL]
  \item \SI{21}{\percent} smaller for \RK[RDPK]{3}{5}[\ESstarp][FSAL]
  \item \SI{32}{\percent} smaller for \ssp43
\end{itemize}
than the asymptotic CFL factor to avoid a blow-up of the simulation.
The required number of right-hand side evaluations and rejected steps are
listed in Table~\ref{tab:delta_wing}.

\begin{table}[!htb]
\sisetup{output-exponent-marker=}
\sisetup{scientific-notation=fixed, fixed-exponent=0}
\centering
\caption{Performance of representative RK methods with error-based step size
         controllers:
         Number of function evaluations (\#FE), accepted steps (\#A),
         and rejected steps (\#R) for the delta wing setup up to time
         $t = 0.015$.}
\label{tab:delta_wing}
\setlength{\tabcolsep}{0.75ex}
\begin{tabular*}{\linewidth}{@{\extracolsep{\fill}}c *2c r@{\hskip 0.5ex}rr@{\hskip 0.5ex}r@{\hskip 1ex}cr@{\hskip 0.5ex}r@{\hskip 0.5ex}r} 
  \toprule
  Scheme & $\beta$ & $\tol$ & \multicolumn{1}{c}{\#FE} & \multicolumn{1}{c}{\#A} & \multicolumn{1}{c}{\#R} \\ 
  \midrule

  \RK[BS]{3}{3}[][FSAL]            & $(0.60, -0.20, 0.00)$
     &  $\tol = 10^{-8}$ & $\num{17641}$ & $5875$ & $4$ \\
  \RK[RDPK]{3}{5}[\ESstarp][FSAL]  & $(0.70, -0.23, 0.00)$
     &  $\tol = 10^{-8}$ & $\num{14963}$ & $2989$ & $3$ \\
  \ssp43                           & $(0.55, -0.27, 0.05)$
     &  $\tol = 10^{-8}$ & $\num{13392}$ & $3344$ & $4$ \\
  \bottomrule
\end{tabular*}
\end{table}

\subsection{NASA common research model}
\label{sec:nasa_crm}

The NASA common research model (CRM) was conceived in 2007. Its aerodynamic
design was completed in 2008, responding to needs broadly expressed within the
US and international aeronautics communities for modern/industry-relevant
geometries coupled with advanced experimental data for applied
computational fluid dynamic validation studies \cite{nasa-crm-history-2019}.
Here, we consider the compressible Navier-Stokes equations for a flow over the
NASA CRM at an angle of attack of \SI{10}{\degree},
a Mach number $\mathrm{Ma} = 0.85$, and a Reynolds number of
$\mathrm{Re} = 5 \times 10^{6}$.

\begin{figure}[!ht]
\centering
    \includegraphics[width=0.8\columnwidth]{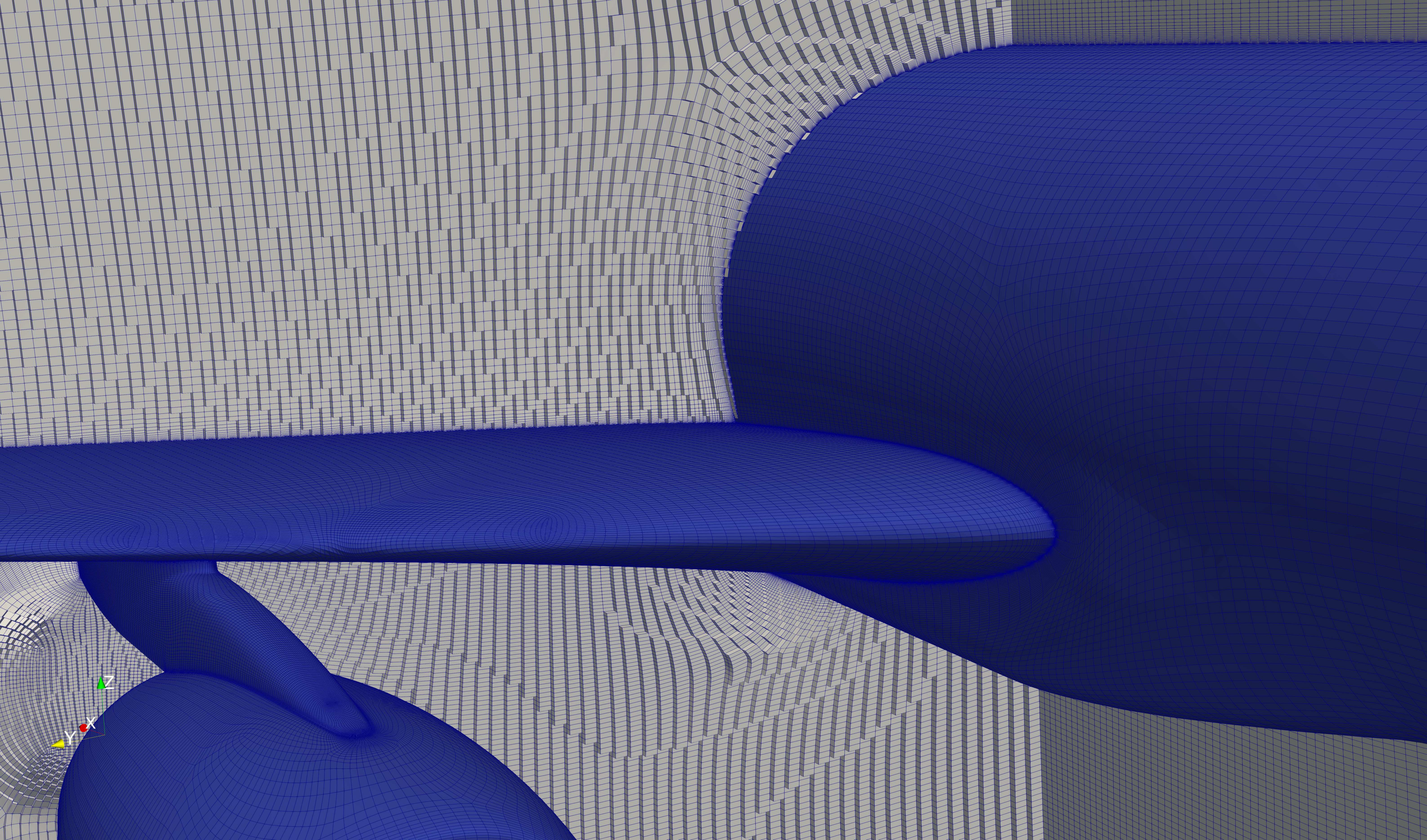}
\caption{Zoom of the mesh for the NASA CRM generated at KAUST.}
\label{fig:zoomCRMgp}
\end{figure}

The computational
domain contains $12.8 \times 10^{6}$ hexahedral elements with a maximum aspect
ratio of approximately $105$.
Entropy stable adiabatic no-slip wall boundary conditions \cite{dalcin2019conservative}
are applied to the aircraft, whereas freestream far-field boundary conditions are applied
at the inlet and outlet planes. We set a solution polynomial degree $p = 3$ in
the whole domain.
Given the degree of the solution and the number of cells, the number of DOFs is
$\approx 7.68 \times 10^{8}$. A zoom of the mesh near the surface of the right wing
and the nacelle is shown in Figure~\ref{fig:zoomCRMgp}.

\begin{figure}[!htb]
\centering
  \begin{subfigure}{0.49\linewidth}
    \includegraphics[width=\textwidth]{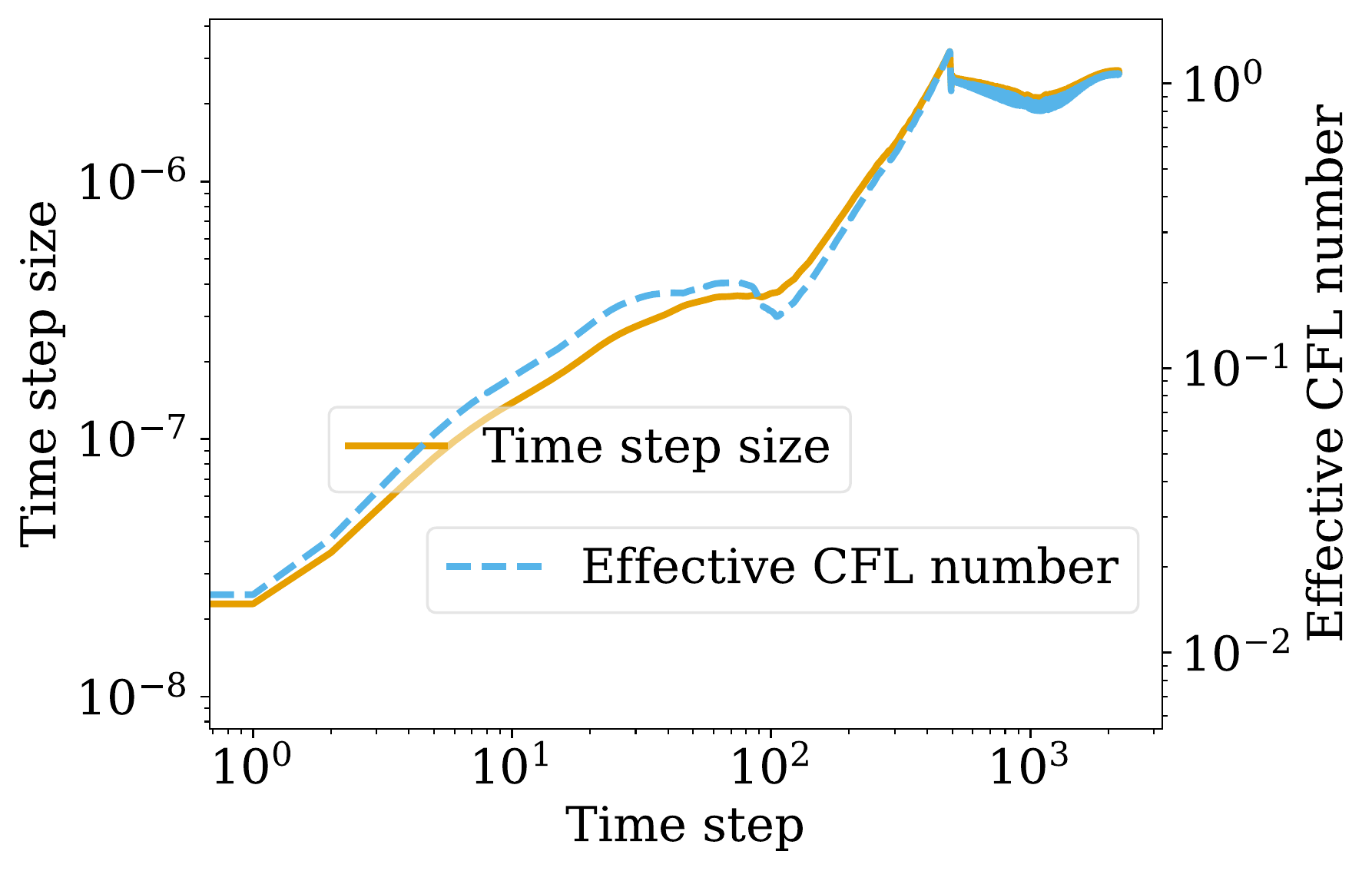}
    \caption{\RK[BS]{3}{3}[][FSAL] (time step number).}
  \end{subfigure}%
  \hspace*{\fill}
  \begin{subfigure}{0.49\linewidth}
    \includegraphics[width=\textwidth]{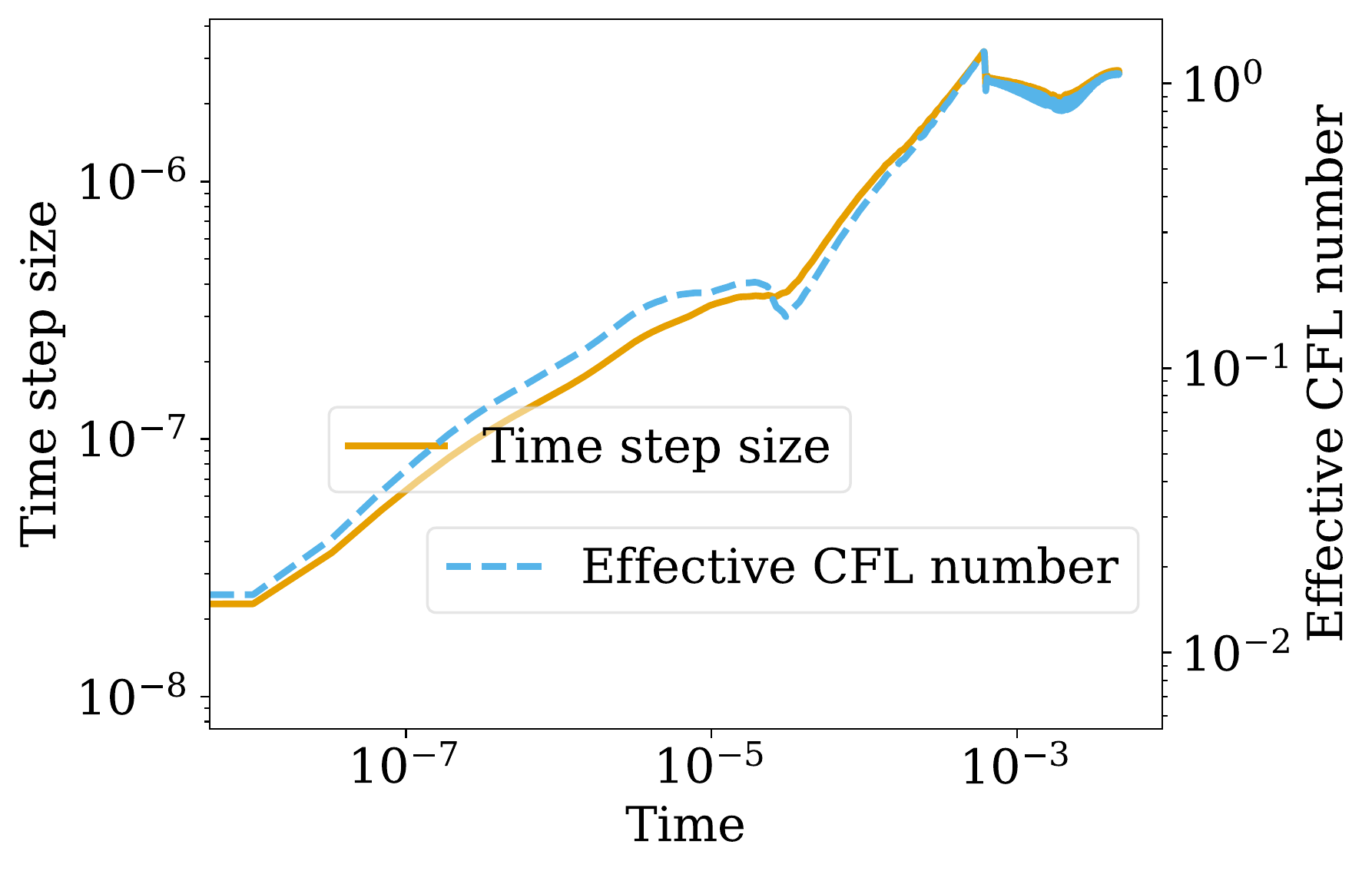}
    \caption{\RK[BS]{3}{3}[][FSAL] (physical time).}
  \end{subfigure}%
  \\
  \begin{subfigure}{0.49\linewidth}
    \includegraphics[width=\textwidth]{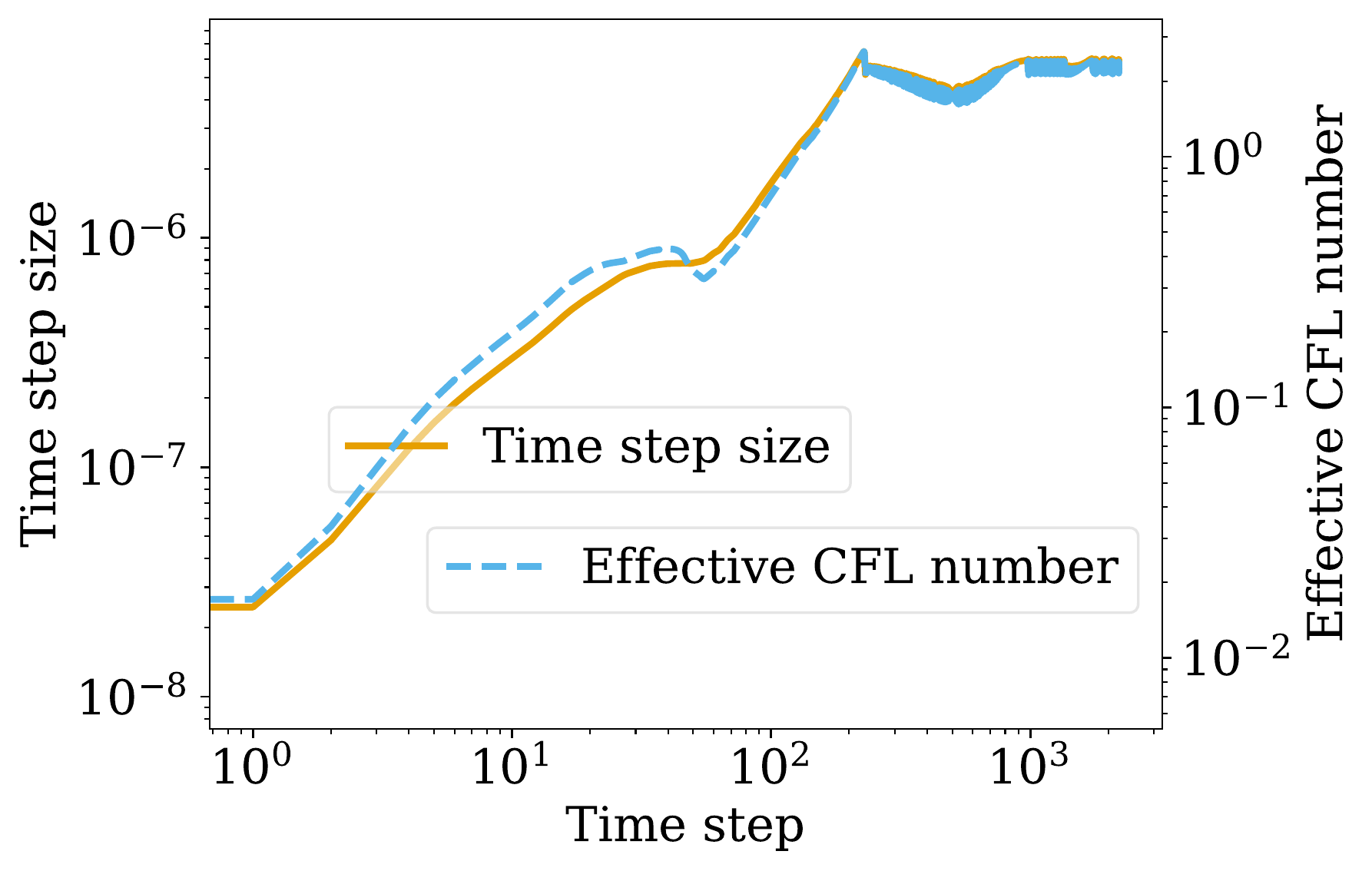}
    \caption{\RK[RDPK]{3}{5}[\ESstarp][FSAL] (time step number).}
  \end{subfigure}%
  \hspace*{\fill}
  \begin{subfigure}{0.49\linewidth}
    \includegraphics[width=\textwidth]{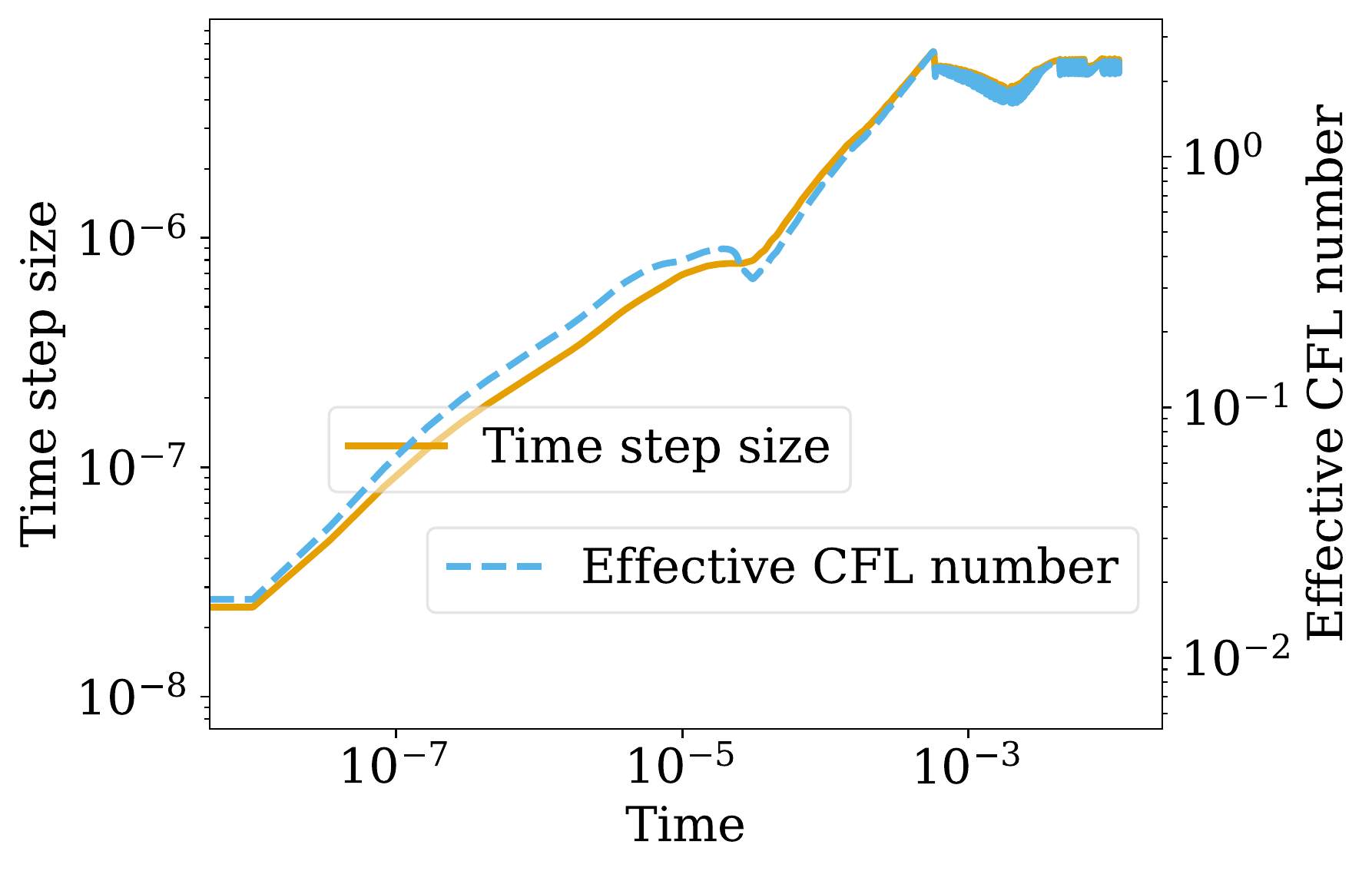}
    \caption{\RK[RDPK]{3}{5}[\ESstarp][FSAL] (physical time).}
  \end{subfigure}%
  \\
  \begin{subfigure}{0.49\linewidth}
    \includegraphics[width=\textwidth]{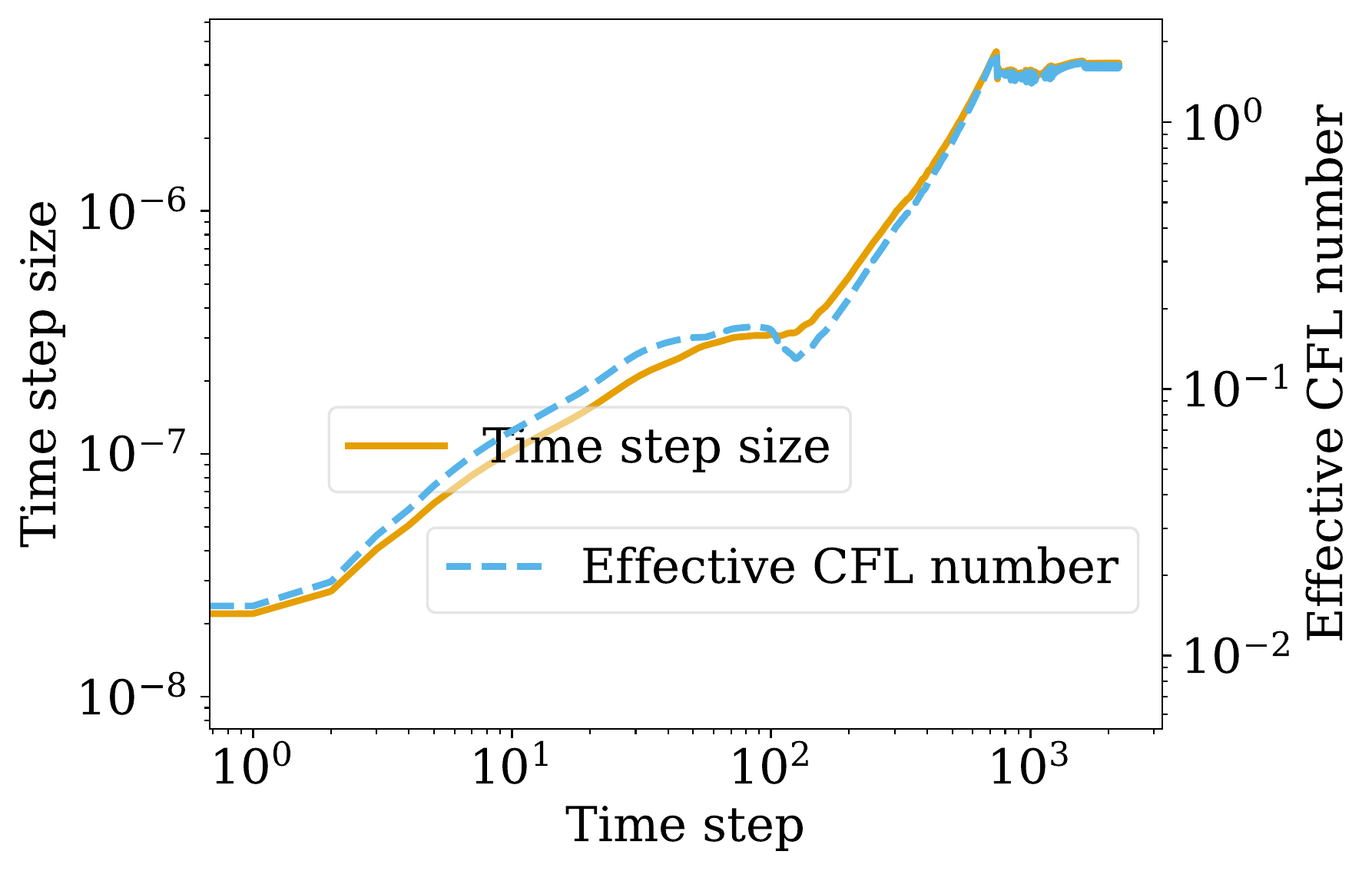}
    \caption{\ssp43 (time step number).}
  \end{subfigure}%
  \hspace*{\fill}
  \begin{subfigure}{0.49\linewidth}
    \includegraphics[width=\textwidth]{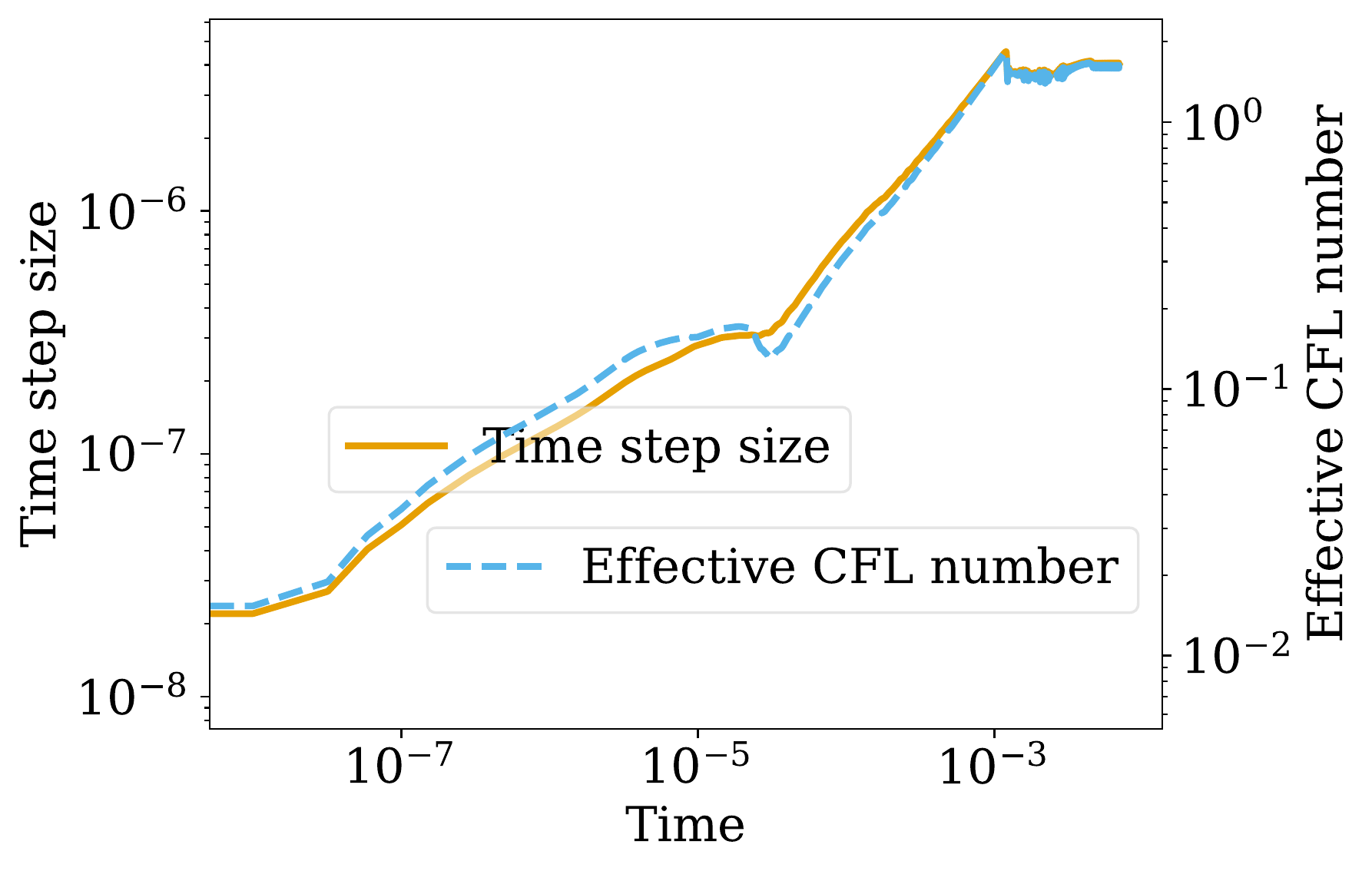}
    \caption{\ssp43 (physical time).}
  \end{subfigure}%
  \caption{Effective CFL numbers and time step sizes depending on the
           time step index (left column) and physical time (right column)
           of some representative RK methods applied to entropy-dissipative
           semidiscretizations of the compressible Navier-Stokes equations for the
           NASA common research model.}
  \label{fig:nasa_crm}
\end{figure}

Results of simulations using some representative RK methods are shown in
Figure~\ref{fig:nasa_crm}. Again, all methods use a conservative estimate of
the initial time step size $\dt$, which then increases $\dt$ monotonically in the first
few hundred time steps and reaches an asymptotically approximately constant step
size afterwards.
The effective CFL number shown in Figure~\ref{fig:nasa_crm} is, again, basically
proportional to the time step size itself. Running SSDC with CFL-based step
size control requires a CFL factor that is (at least initially)
\begin{itemize}
  \item \SI{16}{\percent} smaller for \RK[RDPK]{3}{5}[\ESstarp][FSAL]
  \item \SI{9}{\percent} smaller for \ssp43
\end{itemize}
than the asymptotic CFL factor to avoid a blow-up of the simulation.
For \RK[BS]{3}{3}[][FSAL], the asymptotic CFL factor from error-based step size
control can be chosen for the complete simulation.
The required number of right-hand side evaluations and rejected steps are
listed in Table~\ref{tab:nasa_crm}.

\begin{table}[!htb]
\sisetup{output-exponent-marker=}
\sisetup{scientific-notation=fixed, fixed-exponent=0}
\centering
\caption{Performance of representative RK methods with error-based
         step size controllers:
         Number of function evaluations (\#FE), accepted steps (\#A),
         and rejected steps (\#R) for the NASA common research model up to time
         $t = 0.004$.}
\label{tab:nasa_crm}
\setlength{\tabcolsep}{0.75ex}
\begin{tabular*}{\linewidth}{@{\extracolsep{\fill}}c *2c r@{\hskip 0.5ex}rr@{\hskip 0.5ex}r@{\hskip 1ex}cr@{\hskip 0.5ex}r@{\hskip 0.5ex}r} 
  \toprule
  Scheme & $\beta$ & $\tol$ & \multicolumn{1}{c}{\#FE} & \multicolumn{1}{c}{\#A} & \multicolumn{1}{c}{\#R} \\ 
  \midrule

  \RK[BS]{3}{3}[][FSAL]            & $(0.60, -0.20, 0.00)$
     &  $\tol = 10^{-8}$ & $\num{5874}$ & $1958$ & $0$ \\
  \RK[RDPK]{3}{5}[\ESstarp][FSAL]  & $(0.70, -0.23, 0.00)$
     &  $\tol = 10^{-8}$ & $\num{4610}$ & $922$ & $0$ \\
  \ssp43                           & $(0.55, -0.27, 0.05)$
     &  $\tol = 10^{-8}$ & $\num{5920}$ & $1480$ & $0$ \\
  \bottomrule
\end{tabular*}
\end{table}

\begin{figure}[!ht]
\centering
    \includegraphics[width=0.8\columnwidth]{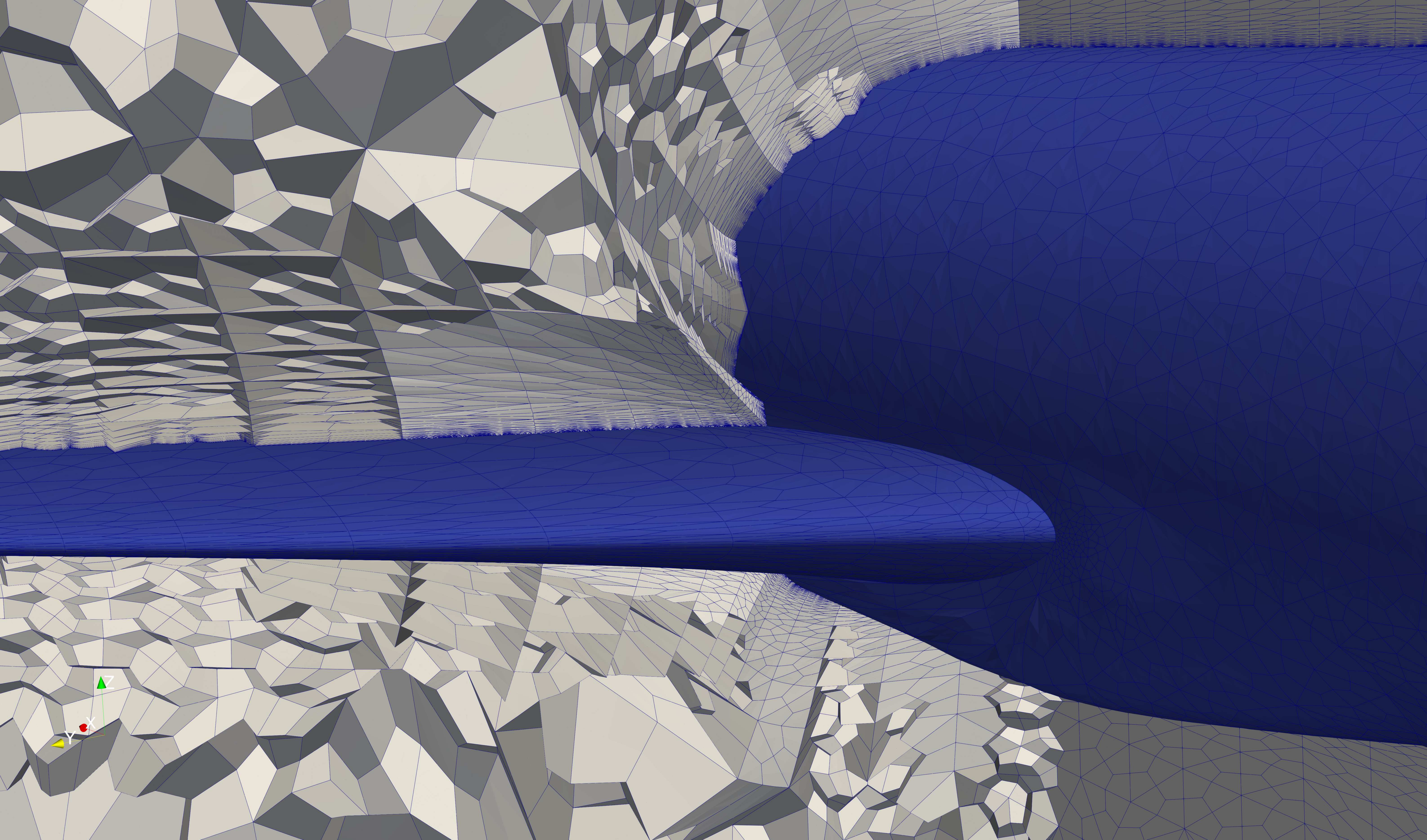}
\caption{Zoom of the CRM-WB-a2.75-Coarse-P2A mesh for the NASA CRM
    provided by the fifth high-order
    workshop \cite{high_order_workshop_5}.}
\label{fig:zoomItsyBitsy}
\end{figure}

To conclude this test case, we highlight that
the error-based controller yields a substantially higher level of
robustness for this industry-relevant test case under a change of mesh structure. To verify that, we simulate the
same flow problem using the mesh illustrated in Figure~\ref{fig:zoomItsyBitsy}. This new mesh is one of the grids provided for the fifth high-order workshop held at Ceneaero, Belgium, and it corresponds to the mesh tagged as ``CRM-WB-a2.75-Coarse-P2A''.
When the CFL-based controller is used, the target CFL number is set to the effective
asymptotic value obtained from the previous
simulations with the mesh shown in Figure~\ref{fig:zoomCRMgp}, \ie approximately
$3.6$ for all
three RK methods (see Figure~\ref{fig:nasa_crm}). The simulations run successfully to the final time for
all the three time integration methods
tested, \ie, \RK[BS]{3}{3}[][FSAL], \RK[RDPK]{3}{5}[\ESstarp][FSAL], and \ssp43,
when the error-based time step
controller is used. On the contrary, all simulations fail when we choose the CFL-based time step
controller. The reason is simple, the effective asymptotic value of the
CFL number for this new mesh is approximately a factor of three smaller than the
value set, \ie, $\approx 1.2$. This indicates that the CFL-based
controller is not sufficiently robust under a change of mesh. In many cases, we
could mitigate the instability by choosing a ``very small value'' for the target
CFL number, but this would
yield a very inefficient and computationally expensive time advancement of the simulation.

\section{Change of variables: dissipative methods using an entropy projection}
\label{sec:GaussSBP}

Whenever a semidiscretization is not strictly a collocation method, there are
some transformations between the basic solution variables evolved in time and
nodal values required to compute fluxes, \eg,
methods using modal coefficients \cite{meister2012comparison,chan2019skew},
entropy variables as primary unknowns \cite{hughes1986new, hiltebrand2014entropy},
staggered grid methods \cite{kopriva1996conservative,parsani2016entropy},
and discontinuous Galerkin difference methods \cite{yan2021entropy}.
Here, we apply entropy-stable Gauss methods \cite{chan2019efficient}
available in Trixi.jl \cite{ranocha2022adaptive,schlottkelakemper2021purely}
with efficient implementations \cite{ranocha2021efficient}. While these
schemes are more complex than collocation methods due to the entropy projection,
they also possess favorable robustness properties for multidimensional simulations of
under-resolved compressible flows with variable density and small-scale features
\cite{chan2022entropy}.

We consider a Kelvin-Helmholtz instability of the 2D compressible Euler equations
of an ideal gas with ratio of specific heats $\gamma = 1.4$ as in \cite{chan2022entropy}.
The domain $[-1, 1]^2$ is equipped with periodic boundary conditions and the
simulation is initialized with
\begin{equation}
\begin{gathered}
  \rho_1 = 1, \quad
  \rho_2 = \rho_1 \frac{1+A}{1-A}, \quad
  \rho^0 = \rho_1 + B(x,y) (\rho_2 - \rho_1), \\
  p^0 = 1, \quad
  v_1^0 = B(x,y) - \frac{1}{2}, \quad
  v_2^0 = \frac{1}{10} \sin(2\pi x),
\end{gathered}
\end{equation}
where $A$ is the Atwood number parameterizing the density contrast and
\begin{equation}
  B(x,y) = \tanh(15 y + 7.5) - \tanh(15 y - 7.5)
\end{equation}
is a smoothed step function. We apply Gauss methods with 32 elements per coordinate
direction and polynomials of degree $\polydeg = 3$. We choose the entropy conservative
flux of \cite{ranocha2020entropy,ranocha2018thesis,ranocha2021preventing} for
the volume terms and a local Lax-Friedrichs flux at surfaces. We apply several
representative RK methods to evolve the resulting entropy-dissipative
semidiscretizations in time $t \in [0, 5]$.

\subsection{Mild initial conditions with low Atwood number}

\begin{table}[!htb]
\sisetup{output-exponent-marker=}
\sisetup{scientific-notation=fixed, fixed-exponent=0}
\centering
\caption{Performance of representative RK methods with default error-based and
         manually tuned CFL-based step size controllers:
         Number of function evaluations (\#FE), accepted steps (\#A),
         and rejected steps (\#R) for the Kelvin-Helmholtz instability with
         low Atwood number $A = 3/7$ and entropy-dissipative Gauss collocation
         semidiscretizations using entropy projections.}
\label{tab:kelvin_helmholtz_low_atwood}
\setlength{\tabcolsep}{0.75ex}
\begin{tabular*}{\linewidth}{@{\extracolsep{\fill}}c *2c r@{\hskip 0.5ex}rr@{\hskip 0.5ex}r@{\hskip 1ex}cr@{\hskip 0.5ex}r@{\hskip 0.5ex}r} 
  \toprule
  Scheme & $\beta$ & $\tol$/$\cfl$ & \multicolumn{1}{c}{\#FE} & \multicolumn{1}{c}{\#A} & \multicolumn{1}{c}{\#R} \\ 
  \midrule

  \RK[BS]{3}{3}[][FSAL]            & $(0.60, -0.20, 0.00)$
     &  $\tol = 10^{-4}$ & $ 7542$ & $  2511$ & $    2$ \\ 
   & &  $\cfl =  0.57  $ & $ 7801$ & $  2600$ &         \smallskip\\ 
  \RK[RDPK]{3}{5}[\ESstarp][FSAL]  & $(0.70, -0.23, 0.00)$
     &  $\tol = 10^{-4}$ & $ 6403$ & $  1278$ & $    2$ \\ 
   & &  $\cfl =  1.10  $ & $ 6856$ & $  1371$ &         \smallskip\\ 
  \RK[RDPK]{4}{9}[\ESstarp][FSAL]  & $(0.38, -0.18, 0.01)$
     &  $\tol = 10^{-5}$ & $ 6537$ & $  722$ & $    4$ \\ 
   & &  $\cfl =  1.85  $ & $ 6922$ & $  769$ &         \smallskip\\ 
  \ssp43                           & $(0.55, -0.27, 0.05)$
     &  $\tol = 10^{-4}$ & $ 5478$ & $  1367$ & $    2$ \\ 
   & &  $\cfl =  1.03  $ & $ 5672$ & $  1418$ &         \\ 
  \bottomrule
\end{tabular*}
\end{table}

Table~\ref{tab:kelvin_helmholtz_low_atwood} shows a summary of the performance statistics
of some representative RK methods for a low Atwood number $A = 3/7$. The CFL
number $\cfl$ were maximized manually up to three significant digits (so that
the simulations did not crash). We used time integration methods with controllers
recommended in \cite{ranocha2021optimized} and a default tolerance $\tol = 10^{-4}$;
we set $\tol = 10^{-5}$ for \RK[RDPK]{4}{9}[\ESstarp][FSAL] since the simulation
crashed for the looser default tolerance.

The total number of function evaluations and time steps is always smaller with
error-based step size control in this example, usually around \SI{5}{\percent}.
In addition, the manual optimization of the CFL factor $\cfl$ required several
simulation runs. Note that the optimal CFL factors are different from
the ones used for the Orzsag-Tang vortex while the same tolerances were used.
For this setup, \ssp43 is the most efficient method, followed by
\RK[RDPK]{3}{5}[\ESstarp][FSAL].

\subsection{Demanding initial conditions with high Atwood number}

Increasing the density contrast makes this under-resolved test case more demanding
\cite{chan2022entropy}. Violations of positivity are more likely to occur and
the sensitivity of the entropy variables for near-vacuum states can cause issues
\cite{chan2018discretely,yan2021entropy}. In particular, the basic variables evolved
in time for the Gauss methods violate positivity properties of the density/pressure
around $t \approx 3$ even for strict CFL factors such as $\cfl = 0.05$ for \ssp43
so that a standard CFL-based step size control is not possible. Similar problems
occur for other RK methods.

\begin{table}[!htb]
\sisetup{output-exponent-marker=}
\sisetup{scientific-notation=fixed, fixed-exponent=0}
\centering
\caption{Performance of representative RK methods with default error-based and
         manually tuned CFL-based step size controllers:
         Number of function evaluations (\#FE), accepted steps (\#A),
         and rejected steps (\#R) for the Kelvin-Helmholtz instability with
         high Atwood number $A = 0.7$ and entropy-dissipative Gauss collocation
         semidiscretizations using entropy projections.}
\label{tab:kelvin_helmholtz_high_atwood}
\setlength{\tabcolsep}{0.75ex}
\begin{tabular*}{\linewidth}{@{\extracolsep{\fill}}c *2c r@{\hskip 0.5ex}rr@{\hskip 0.5ex}r@{\hskip 1ex}cr@{\hskip 0.5ex}r@{\hskip 0.5ex}r} 
  \toprule
  Scheme & $\beta$ & $\tol$/$\cfl$ & \multicolumn{1}{c}{\#FE} & \multicolumn{1}{c}{\#A} & \multicolumn{1}{c}{\#R} \\ 
  \midrule

  \RK[BS]{3}{3}[][FSAL]            & $(0.60, -0.20, 0.00)$
     &  $\tol = 10^{-4}$ & $ 9567$ & $  3185$ & $    3$ \\ 
  \RK[RDPK]{3}{5}[\ESstarp][FSAL]  & $(0.70, -0.23, 0.00)$
     &  $\tol = 10^{-4}$ & $ 8173$ & $  1621$ & $   13$ \\ 
  \RK[RDPK]{4}{9}[\ESstarp][FSAL]  & $(0.38, -0.18, 0.01)$
     &  $\tol = 10^{-6}$ & $10587$ & $  1175$ & $    1$ \\ 
  \ssp43                           & $(0.55, -0.27, 0.05)$
     &  $\tol = 10^{-4}$ & $ 6618$ & $  1640$ & $   14$ \\ 
  \bottomrule
\end{tabular*}
\end{table}

In contrast, error-based step size control works directly in most
cases, as shown in Table~\ref{tab:kelvin_helmholtz_high_atwood}. We only needed
to use a stricter tolerance $\tol = 10^{-6}$ for \RK[RDPK]{4}{9}[\ESstarp][FSAL].
This RK method has more stages than the other methods, making it more difficult
to adapt quickly to changing conditions. Again, \ssp43 is the most efficient
method, followed by \RK[RDPK]{3}{5}[\ESstarp][FSAL].

\section{New systems: convenient exploratory research}
\label{sec:new_systems}

A key feature of error-based time stepping methods is that they require no knowledge of the maximum eigenvalues to select a stable explicit time step.
This feature is particularly convenient when one wants to apply existing numerical methods to new equation systems.
To illustrate this, we consider the shallow water equations augmented with an Exner equation to account for
interactions between the bottom topography and the fluid flow.
The coupled system of equations in two spatial dimensions is
\begin{equation}
\label{eq:SWExner_equation}
\frac{\partial}{\partial t}
\begin{pmatrix}
h\\[0.1cm]
hv_1\\[0.1cm]
hv_2\\[0.1cm]
b
\end{pmatrix}
+
\frac{\partial}{\partial x}
\begin{pmatrix}
hv_1\\[0.1cm]
hv_1^2 + \frac{g}{2}h^2\\[0.1cm]
hv_1v_2\\[0.1cm]
\xi q_1(\vec{v})
\end{pmatrix}
+\frac{\partial}{\partial y}
\begin{pmatrix}
hv_2\\[0.1cm]
hv_1v_2\\[0.1cm]
hv_2^2 + \frac{g}{2}h^2\\[0.1cm]
\xi q_2(\vec{v})
\end{pmatrix}
=
\begin{pmatrix}
0\\[0.1cm]
-gh\frac{\partial b}{\partial x}\\[0.1cm]
-gh\frac{\partial b}{\partial y}\\[0.1cm]
0
\end{pmatrix},
\end{equation}
where $h(x,y,t)$ is the water height, $\vec{v}(x,y,t)=(v_1(x,y,t),v_2(x,y,t))^T$ is the velocity field,
$b(x,y,t)$ is the evolving bottom topography, and $g$ is the gravitational constant.

The last line in \eqref{eq:SWExner_equation} is the Exner equation for the time evolution of the bottom topography $b$.
It contains the constant $\xi = 1/(1-\sigma)$ where $\sigma\in(0,1)$ is the porosity of the bottom material, and
the solid transport discharge $\vec{q}=(q_1,q_2)^T$ as a function of the flow velocity $\vec{v}$.
A closure model is required to describe how the bottom topography couples to the flow and defines the form of the $\vec{q}$ terms.
This closure depends on the characteristics of the bottom material and flow, e.g., the grain size or Froude number.
There exist many closure models for the solid transport discharge proposed and studied
(analytically as well as numerically) in the literature, see e.g.,
\cite{briganti2012efficient,garres2020shallow,grass1981,meyer1948formulas}.
For the present study, we consider the simple and frequently used model due to Grass \cite{grass1981} that models
the instantaneous sediment transport as a power law on the velocity field:
\begin{equation}
\label{eq:Grass_model}
\vec{q}(\vec{v}) = A_g \vec{v} \left(v_1^2+v_2^2\right)^{\frac{m-1}{2}}.
\end{equation}
In \eqref{eq:Grass_model}, the constant $A_g\in[0,1]$ accounts for the porosity of the bottom sediment
layer as well as the effects of the grain size and is usually determined from experimental data.
When $A_g$ is zero there is no sediment transport and \eqref{eq:SWExner_equation} reduces to the standard two-dimensional shallow water equations.
The interaction between the bottom and the water flow is weak when $A_g$ is small and strong as $A_g$ approaches one.
The factor $m\in[1,4]$ in the Grass model may also be determined from data;
however, if one considers odd integer values of $m$, then \eqref{eq:Grass_model} can be
differentiated and the model remains valid for all values of the velocity field $\vec{v}$.
For the remainder of this discussion we will take $m=3$ and the solid transport discharge will take the form
\begin{equation}
\label{eq:Grass_m_3}
\vec{q}(\vec{v}) = A_g \vec{v} \left(v_1^2+v_2^2\right).
\end{equation}

For the particular closure model \eqref{eq:Grass_m_3}, we investigate the flux Jacobian of the
governing system of equations \eqref{eq:SWExner_equation}.
This is because the choice of a stable explicit time step under a CFL condition requires knowledge
of the fastest wave speed of the underlying system.
For this we compute the flux Jacobian in the $x$-direction, which incorporates non-conservative terms
from the right hand side of \eqref{eq:SWExner_equation}, to be
\begin{equation}
\label{eq:fluxJac}
\frac{\partial\vec{f}}{\partial\vec{u}}
=
 \begin{pmatrix}
 0 & 1 & 0 & 0\\[0.15cm]
 gh - v_1^2 & 2v_1 & 0 & gh\\[0.15cm]
 -v_1v_2 & v_2 & v_1 & 0\\[0.15cm]
 -\frac{3\xi A_g v_1}{h}\left(v_1^2+v_2^2\right) & \frac{\xi A_g}{h}\left(3v_1^2+v_2\right) & \frac{2\xi A_g}{h}v_1v_2 & 0
 \end{pmatrix}.
\end{equation}
Note that the analysis and results in the $y$-direction are analogous, so we will only consider \eqref{eq:fluxJac}
for this eigenvalue analysis for time step purposes.
With the chosen Grass closure model it has been shown that all the eigenvalues of \eqref{eq:fluxJac}
are real and the governing system \eqref{eq:SWExner_equation} is hyperbolic \cite{chertock2020operator,diaz2009two}.
The eigenvalues of the Jacobian in the $x$-direction are $v_1$ and the three roots of the characteristic polynomial
\begin{equation}
\label{eq:charPoly}
\lambda^3 - 2v_1\lambda^2 + (v_1^2 - gh - g\xi A_g (3v_1^2+v_2^2)) \lambda + g\xi A_gv_1(3v_1^2+v_2^2) =0.
\end{equation}
It is possible, although unwieldy, to apply Cardano's formula to determine the real, distinct roots of \eqref{eq:charPoly}.
However, an accurate estimate to the
maximum wave speed from \eqref{eq:charPoly} is necessary for traditional CFL-based time stepping methods.

For many applications of interest the characteristic flow speed of the water is much faster than
the speed of the bottom topography \cite{chertock2020operator,hudson2001numerical,uh2021unstructured}.
In particular, for subcritical flow the value of $v_1$ is relatively small, so it is possible to use a
loose bound of the flow speeds with $v_1 \pm \sqrt{gh}$ where $ \sqrt{gh}$
is the surface wave speed. However, this loose bound applied to a CFL time step condition
can result  in excessive numerical diffusion that acts on
the computed (slowly moving) bottom topography.

To approximate the solution of the shallow water equations with an Exner model \eqref{eq:SWExner_equation}
we use a discontinuous Galerkin spectral element
method in space to create the semi-discretization and apply explicit time integration, see \cite{wintermeyer2017entropy} for details.
To couple the moving bottom topography
$b(x,y,t)$ we use an unsteady approach where the water flow
and riverbed are calculated simultaneously. That is, the water
flow can be either steady or unsteady and the changes in the bed update are
considered to be significant, i.e., the wave speed of the bed-updating equation is
of a similar magnitude to the wave speeds of the water flow. For this approach,
the system of flow and bottom evolution equations are discretised simultaneously.
This high-order discontinuous Galerkin method applied to the shallow water equations is well-balanced
for static and subcritical moving water flow \cite{wintermeyer2017entropy,winters2015comparison}.

As an exemplary test case, we consider the evolution of a smooth sand hill bottom topography
in a uniform flow governed by the shallow water equations.
This test case is a two dimensional numerical simulation that was first proposed and analyzed by de Vriend~\cite{de19872dh}.
It has become a common test case for benchmarking purposes in the moving bottom topography literature, e.g.
\cite{benkhaldoun2010two,chertock2020operator,diaz2009two,hudson2001numerical,uh2021unstructured}, as it is easy
to setup and demonstrates the ability of a numerical scheme to capture the morphology of a
bottom topography governed with \eqref{eq:SWExner_equation}.

It involves the evolution of a conical sand dune in a channel with a non-erodible bottom on a square domain with dimensions \SI{1000}{m} $\times$ \SI{1000}{m}.
The initial conditions of the flow variables for this test case given in conservative variables are
\begin{equation}\label{eq:swe_ics}
\begin{pmatrix}
h\\[0.1cm]
hv_1\\[0.1cm]
hv_2\\[0.1cm]
\end{pmatrix}
=
\begin{pmatrix}
10 -b(x,y,0)\\[0.1cm]
10\\[0.1cm]
0\\[0.1cm]
\end{pmatrix},
\end{equation}
and the initial form of the sediment layer is given by
\begin{equation}\label{eq:conical_dune}
b(x,y,0) = \begin{cases}
\sin^2\Bigl(\frac{(x-300)\pi}{200}\Bigr)\sin^2\Bigl(\frac{(y-400)\pi}{200}\Bigr) & \textrm{if  } 300\leq x\leq 500,\, 400\leq y\leq 600, \\[0.1cm]
0 & \textrm{otherwise.}
\end{cases}
\end{equation}
We set the gravitational constant to be $g = \SI{9.8}{m/s^2}$. The porosity of the sediment material is $\sigma = 0.4$ and the free parameter
from the Grass model \eqref{eq:Grass_m_3} is $A_g = 0.001$. The parameter $A_g$ corresponds to the coupling strength and speed of
interaction between the sediment and the water flow \cite{benkhaldoun2010two,hudson2001numerical}.
Taking a small value for $A_g$ models weak coupling and a simulation must be
run for a longer time period in order to observe significant variations in the sediment layer. Therefore, we take the final time of this simulation to
be $t = 100$ hours ($360{,}000$ seconds). The discharge value $hv_1= \SI{10}{m^2/s}$ from \eqref{eq:swe_ics} gives the Froude number
\begin{equation}
\textrm{Fr} = \frac{\sqrt{v_1^2 + v_2^2}}{\sqrt{gh}} \approx 0.1.
\end{equation}
Thus, the flow is in the subcritical regime and we must set boundary conditions
accordingly, see e.g. \cite{nordstrom2022linear}.
For this test case, the boundary conditions are taken to be walls on the top and bottom of the domain, the left is a subcritical inflow,
and the right is a subcritical outflow.
The complete setup can be found in the reproducibility
repository \cite{ranocha2022errorRepro}.

For the spatial discretization we divide the domain $\Omega = [0,1000]^2$ into $900$ Cartesian elements.
In each element we use polynomials of degree $p=5$ in each spatial direction for the discontinuous Galerkin approximation.
This results in $32{,}400$ DOFs for each equation in \eqref{eq:SWExner_equation}.
The classic shallow water equations for water height, momentum, and bottom topography terms are discretized with the well-balanced, high-order
discontinuous Galerkin method described by Wintermeyer et al. \cite{wintermeyer2017entropy}. The additional Exner equation for sediment transport is discretized
with a standard discontinuous Galerkin spectral element approximation.

\begin{figure}[!htb]
\centering
  \begin{subfigure}{0.49\linewidth}
    \includegraphics[width=\textwidth]{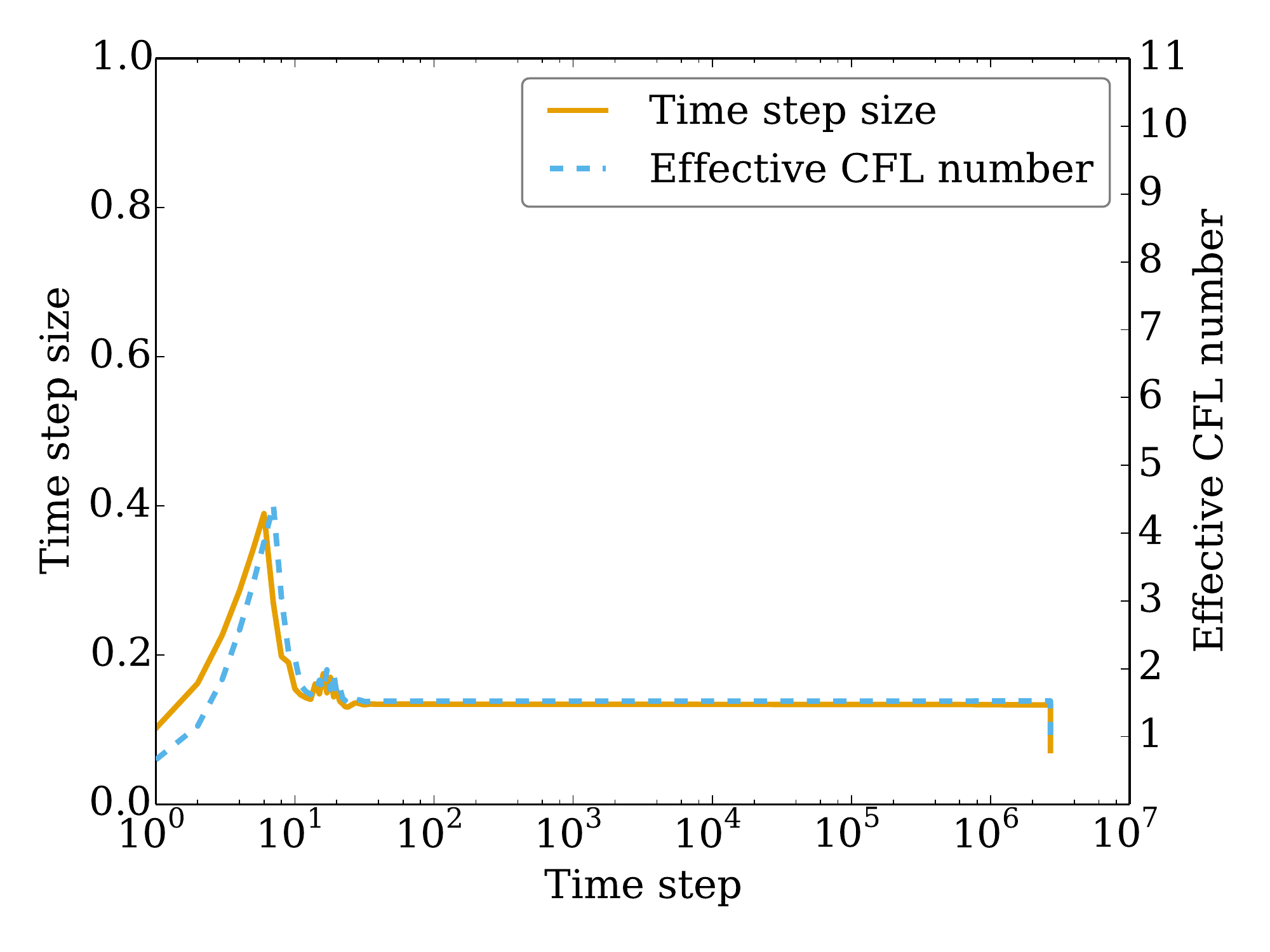}
    \caption{\RK[BS]{3}{3}[][FSAL] (time step number).}
  \end{subfigure}%
  \hspace*{\fill}
  \begin{subfigure}{0.49\linewidth}
    \includegraphics[width=\textwidth]{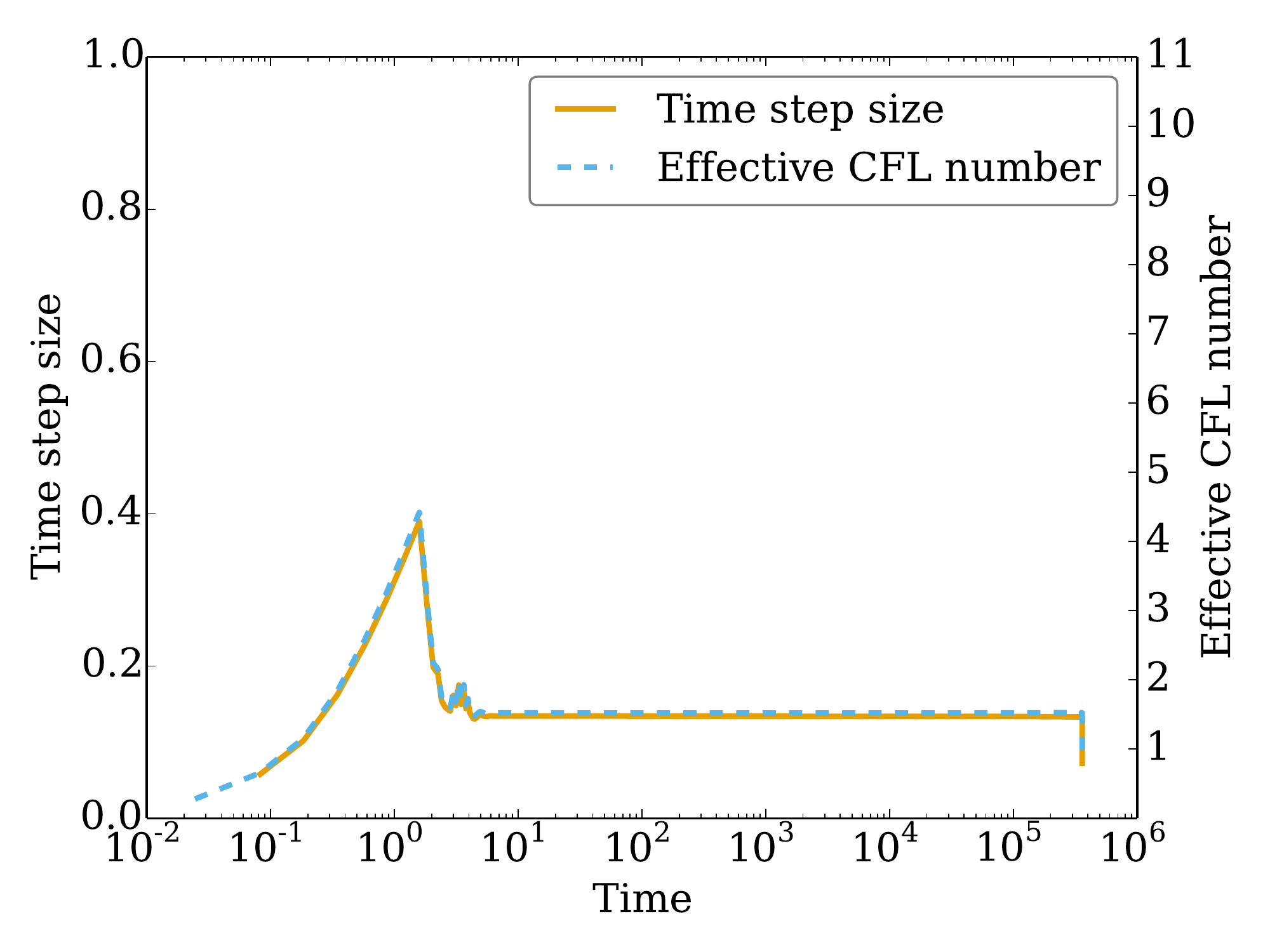}
    \caption{\RK[BS]{3}{3}[][FSAL] (physical time).}
  \end{subfigure}%
  \\
  \begin{subfigure}{0.49\linewidth}
    \includegraphics[width=\textwidth]{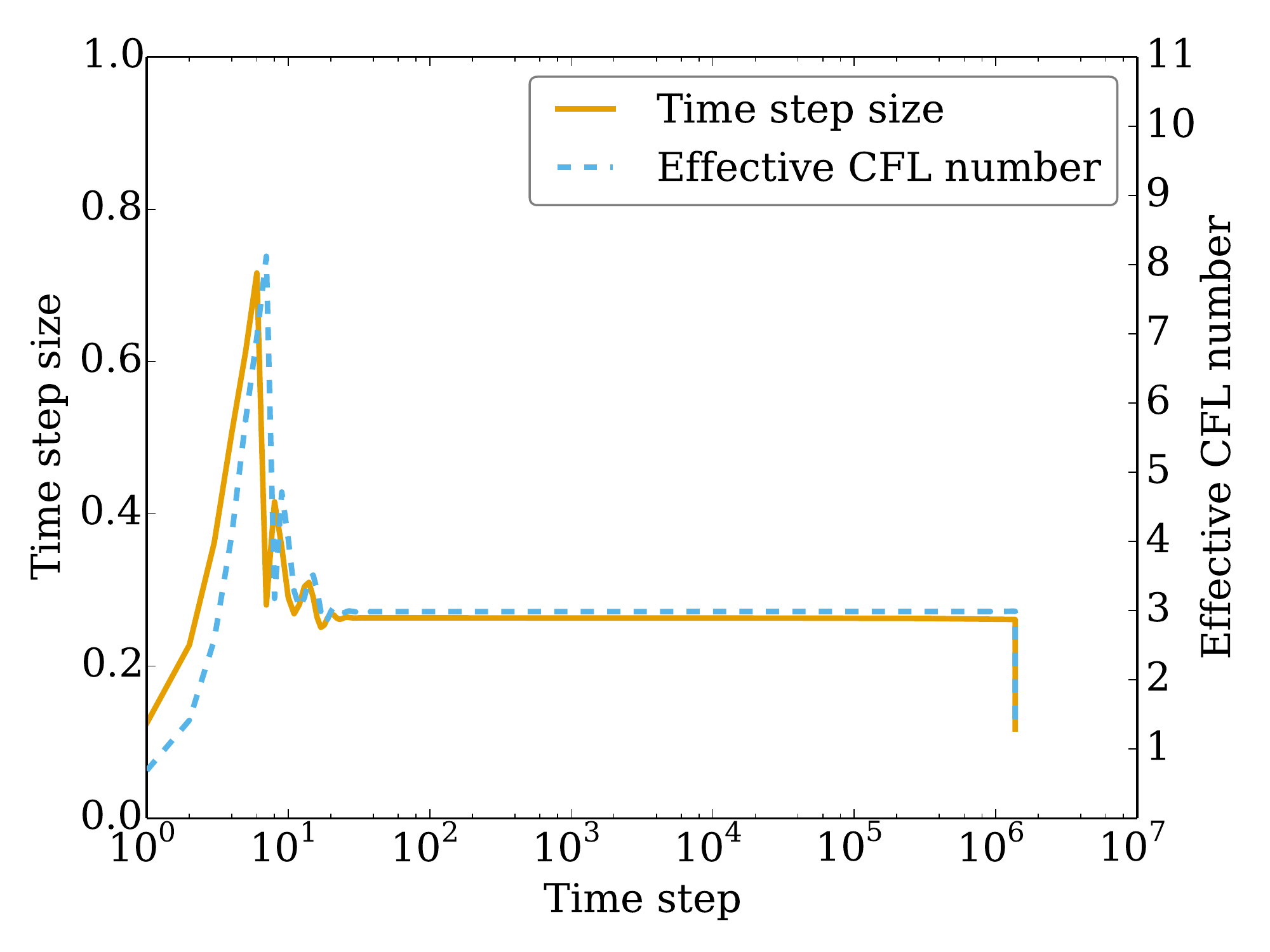}
    \caption{\RK[RDPK]{3}{5}[\ESstarp][FSAL] (time step number).}
  \end{subfigure}%
  \hspace*{\fill}
  \begin{subfigure}{0.49\linewidth}
    \includegraphics[width=\textwidth]{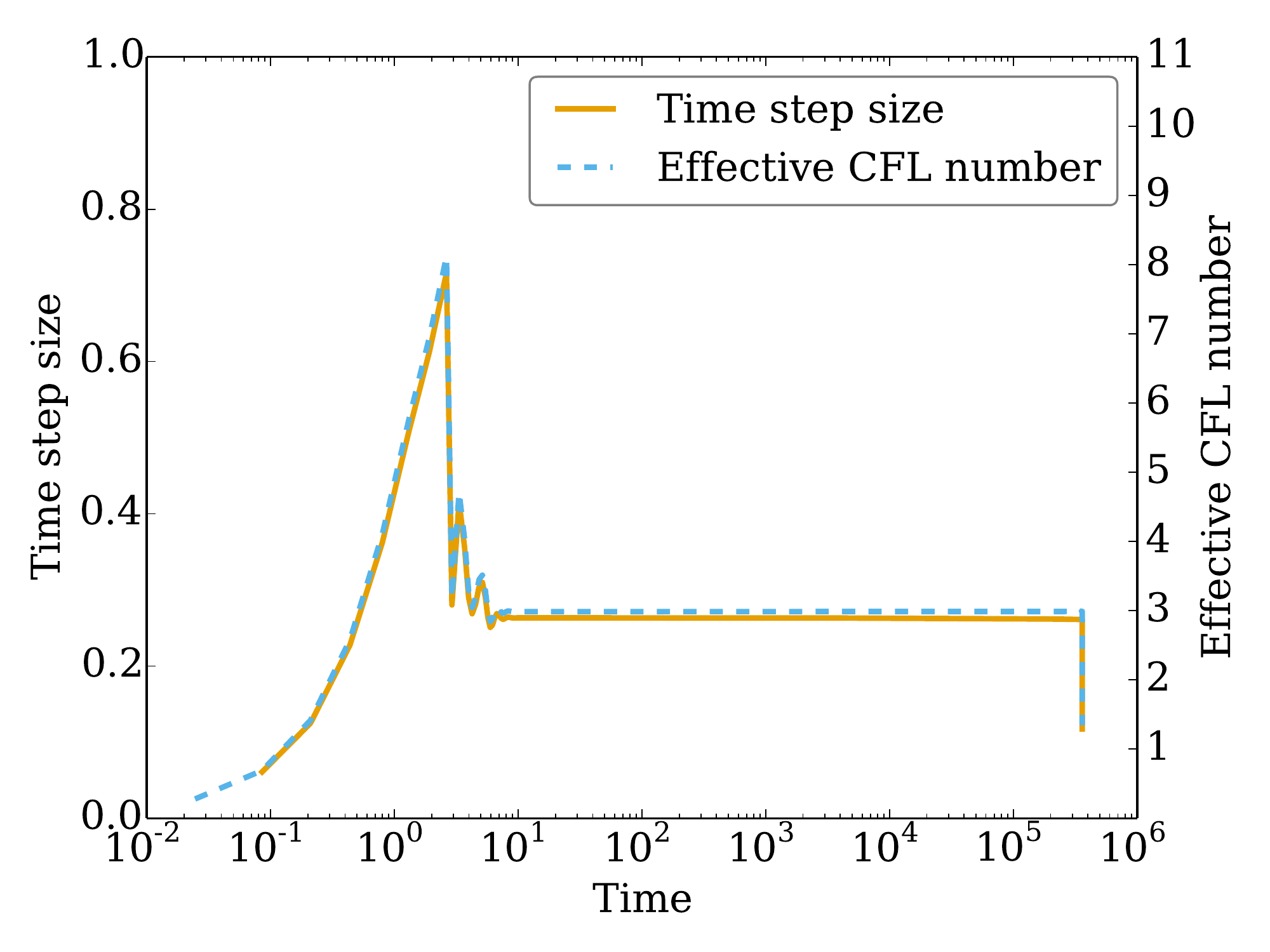}
    \caption{\RK[RDPK]{3}{5}[\ESstarp][FSAL] (physical time).}
  \end{subfigure}%
  \\
  \begin{subfigure}{0.49\linewidth}
    \includegraphics[width=\textwidth]{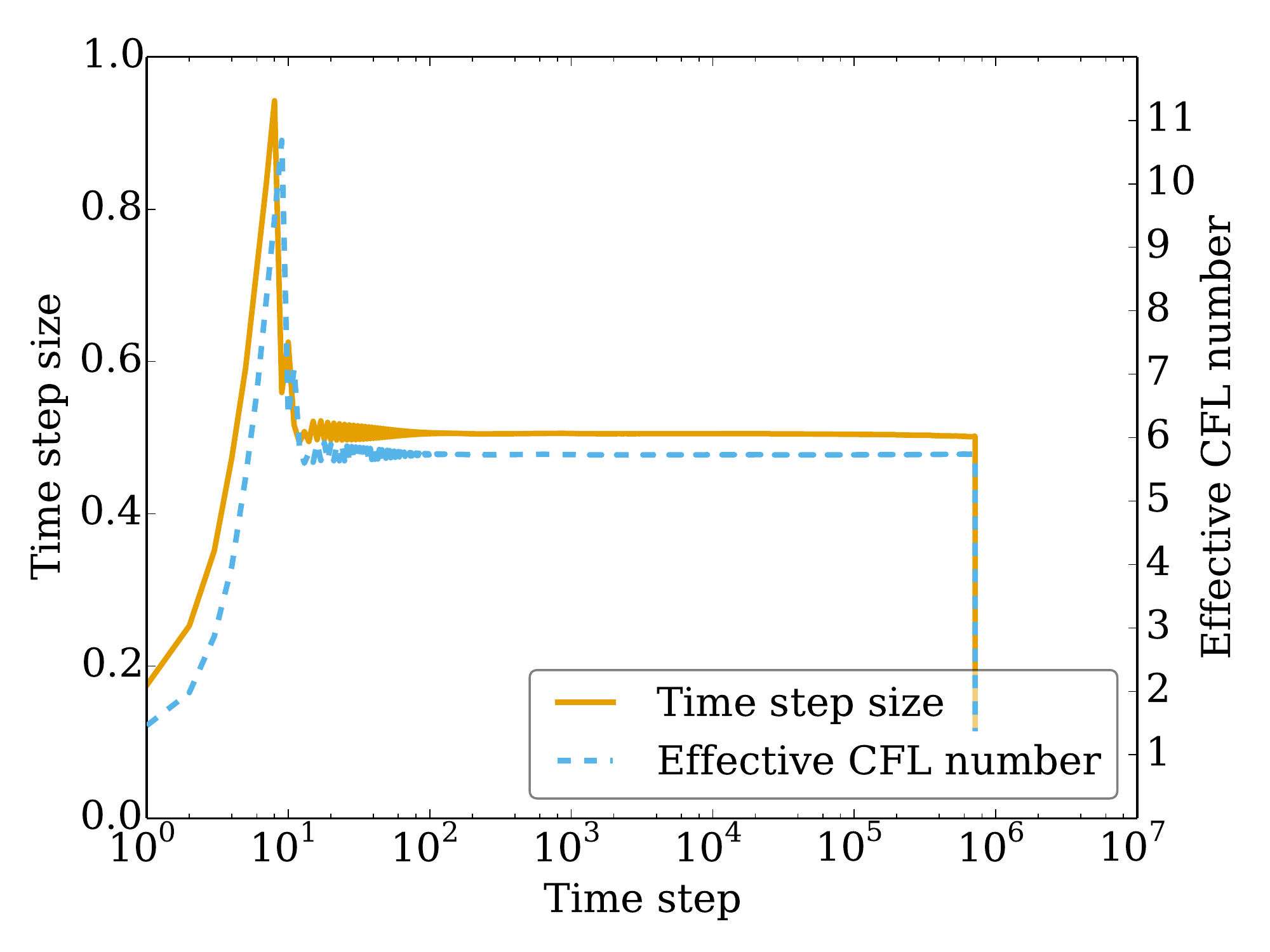}
    \caption{\RK[RDPK]{4}{9}[\ESstarp][FSAL] (time step number).}
  \end{subfigure}%
  \hspace*{\fill}
  \begin{subfigure}{0.49\linewidth}
    \includegraphics[width=\textwidth]{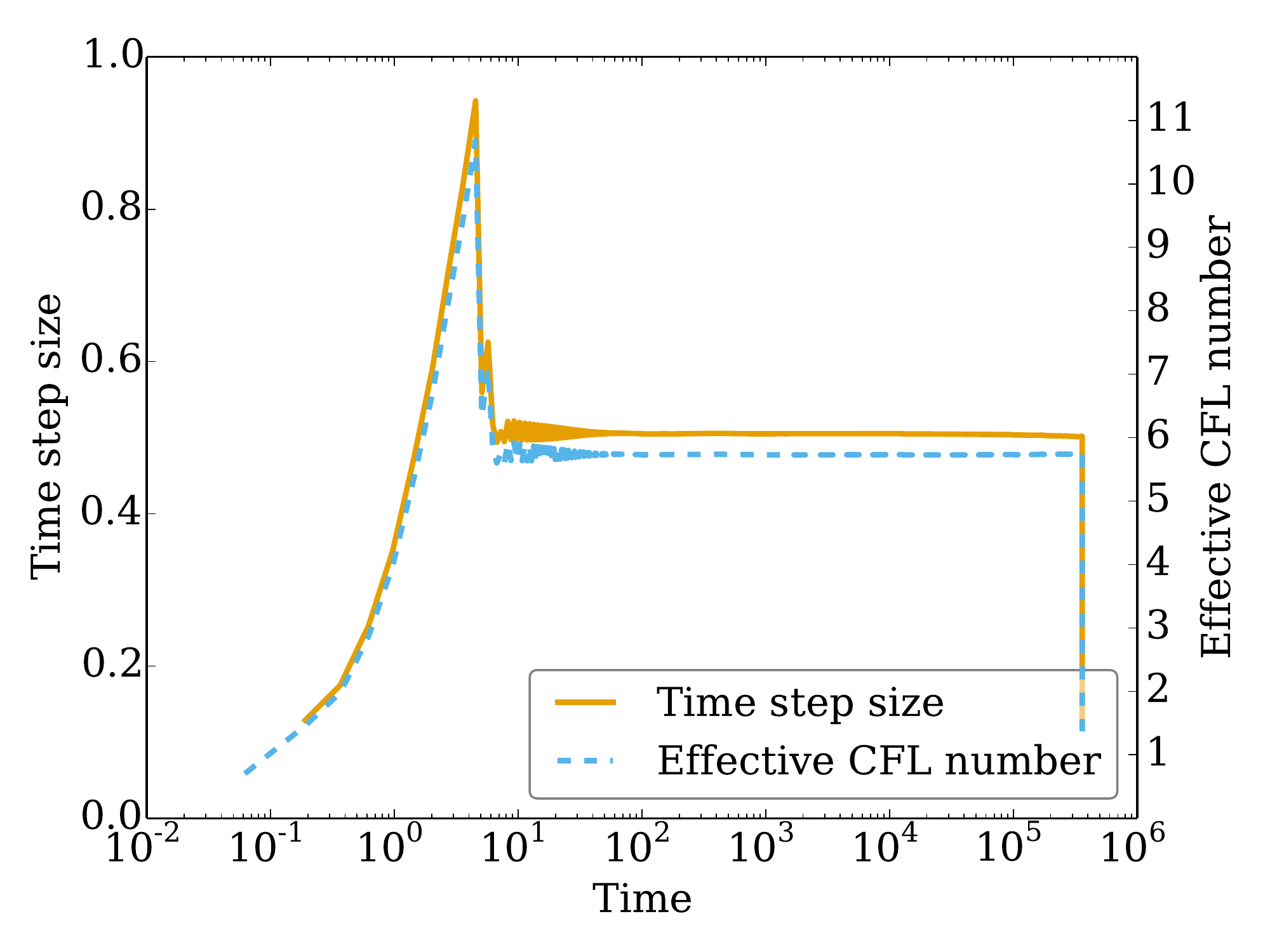}
    \caption{\RK[RDPK]{4}{9}[\ESstarp][FSAL] (physical time).}
  \end{subfigure}%
  \caption{Time step sizes and effective CFL numbers for three error-based RK methods.
           Left column: \# time steps is $x$-axis. Right column: Time is $x$-axis.}
  \label{fig:swe_cfls}
\end{figure}

We first investigate the use of error-based time step control for this shallow water with sediment transport problem setup.
We recorded the time step sizes
and effective CFL numbers after every accepted time step for three representation RK methods.
The tolerance is set to $\tol = 10^{-7}$ for all error-based time
stepping methods considered. All three time integration techniques work directly for this test case.
The behavior of the time step size and
the effective CFL number are presented in Figure~\ref{fig:swe_cfls}.
Initially, the time step size of the three methods increases before the time step settles
to a stable value as the physical problem becomes a quasi-steady flow where the sediment slowly evolves in the presence of the background flow.

We further investigate the effectiveness of the three time integrators with error-based step size control for this test case,
as shown in Table~\ref{tab:swe_exner}. For this quasi-steady state configuration we found that a
stricter tolerance $\tol = 10^{-7}$ for all the integration techniques tested performed the best in the sense
of requiring a lower number of time steps.
Here we find that \RK[RDPK]{4}{9}[\ESstarp][FSAL] is the most efficient
method, followed by \RK[RDPK]{3}{5}[\ESstarp][FSAL].

\begin{table}[!htb]
\sisetup{output-exponent-marker=}
\sisetup{scientific-notation=fixed, fixed-exponent=0}
\centering
\caption{Performance of representative RK methods with default error-based and step size controllers:
         Number of function evaluations (\#FE), accepted steps (\#A),
         and rejected steps (\#R) for a shallow water flow over an evolving sand dune.}
\label{tab:swe_exner}
\setlength{\tabcolsep}{0.75ex}
\begin{tabular*}{\linewidth}{@{\extracolsep{\fill}}c *2c r@{\hskip 0.5ex}rr@{\hskip 0.5ex}r@{\hskip 1ex}cr@{\hskip 0.5ex}r@{\hskip 0.5ex}r} 
  \toprule
  Scheme & $\beta$ & $\tol$/$\cfl$ & \multicolumn{1}{c}{\#FE} & \multicolumn{1}{c}{\#A} & \multicolumn{1}{c}{\#R} \\ 
  \midrule

  \RK[BS]{3}{3}[][FSAL]            & $(0.60, -0.20, 0.00)$
     &  $\tol = 10^{-7}$ & $ 8{,}088{,}687$ & $  2{,}696{,}226$ & $    2$ \\ 
  \RK[RDPK]{3}{5}[\ESstarp][FSAL]  & $(0.70, -0.23, 0.00)$
     &  $\tol = 10^{-7}$ & $ 6{,}873{,}513$ & $  1{,}374{,}701$ & $   1$ \\ 
  \RK[RDPK]{4}{9}[\ESstarp][FSAL]  & $(0.38, -0.18, 0.01)$
     &  $\tol = 10^{-7}$ & $6{,}442{,}851$ & $  715{,}871$ & $    1$ \\ 
  \bottomrule
\end{tabular*}
\end{table}

\begin{figure}[!htb]
  \centering
  \includegraphics[width=0.75\textwidth]{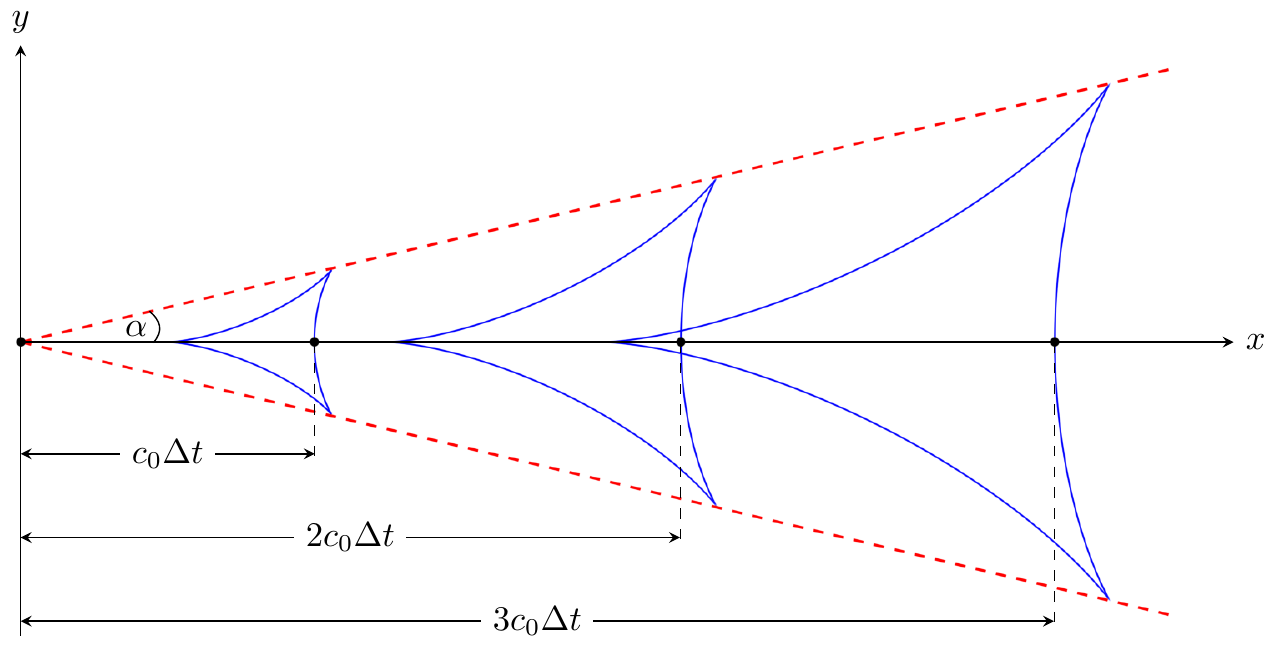}
  \caption{Characteristic evolution of a dune shaped sediment distribution.}
  \label{fig:swe_stars}
\end{figure}

We examine the expected solution behavior as well as simulation results for the sediment transport problem under consideration.
As time evolves, the conical dune \eqref{eq:conical_dune} gradually spreads into a star-shaped pattern.
That is, the sediment will redistribute itself along a triangular shape with a particular spreading angle $\alpha$ and the bulk shape of the dune
will move with a speed $c_0$ as it spreads \cite{de19872dh,hudson2001numerical}. This process is illustrated in Figure~\ref{fig:swe_stars}.
When the Grass model is used and the interaction between the flow and sediment is weak, that is $A_g<0.01$,
de Vriend~\cite{de19872dh} derived an analytical approximation for this spreading angle. In particular, when $m=3$ in the Grass
model \eqref{eq:Grass_m_3}, this spreading angle of a conical dune is
\begin{equation}\label{eq:dune_angle}
\alpha = \tan^{-1}\left(\frac{3\sqrt{3}}{13}\right) \approx \SI{21.787}{\degree}.
\end{equation}

\begin{figure}[!htb]
\centering
  \begin{subfigure}{0.495\linewidth}
    \includegraphics[width=\textwidth]{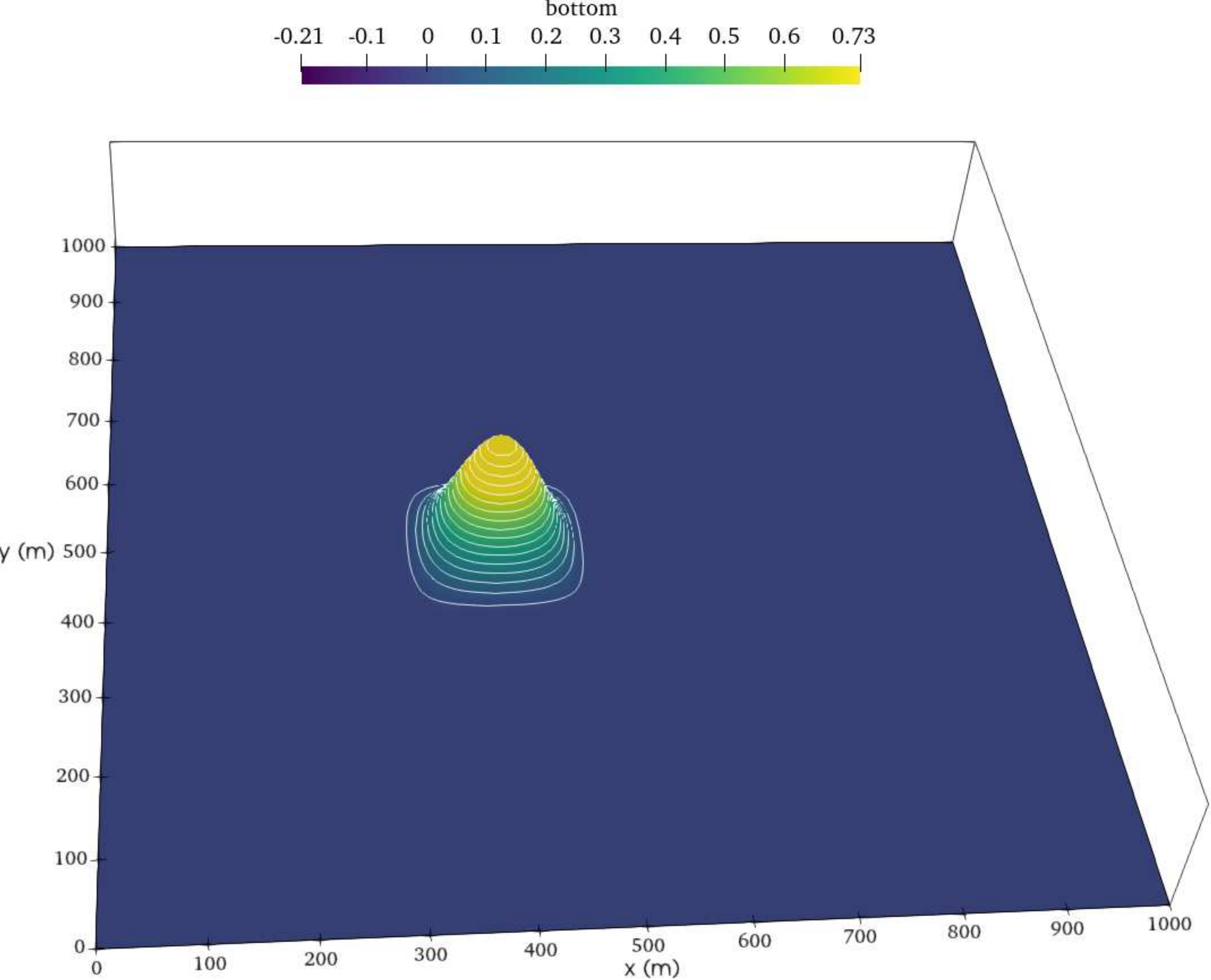}
    \caption{$t=\SI{0}{s}$}
  \end{subfigure}%
  \hspace*{\fill}
  \begin{subfigure}{0.495\linewidth}
    \includegraphics[width=\textwidth]{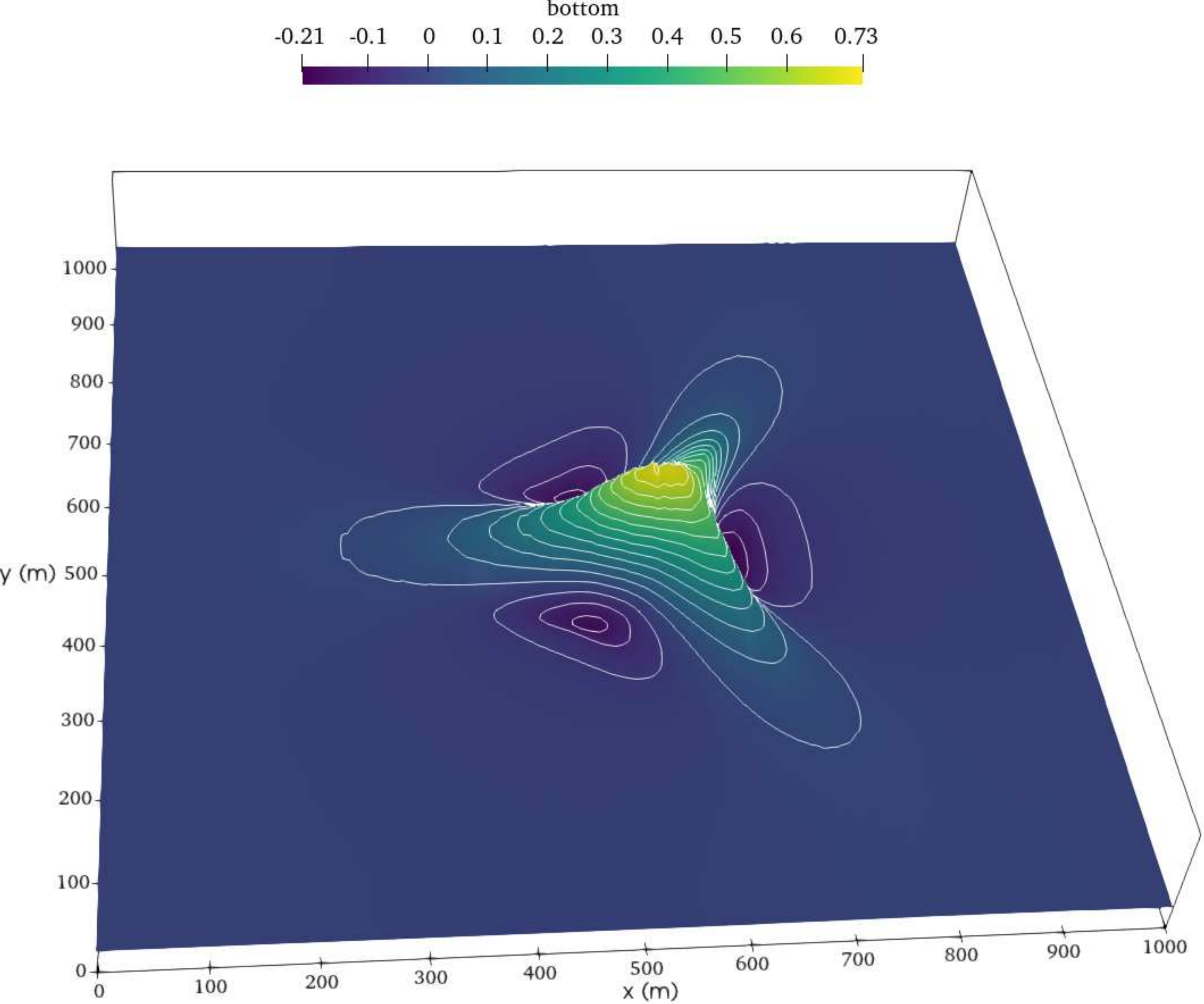}
    \caption{$t=\SI{360000}{s}$}
  \end{subfigure}%
  \caption{Projected lateral view of the approximate solution of the bottom topography $b(x,y,t)$
               that demonstrates expected star-shaped structure of the conical sand dune at the final time.}
  \label{fig:swe_project_sol}
\end{figure}

\begin{figure}[!htb]
\centering
  \begin{subfigure}{0.495\linewidth}
    \includegraphics[width=\textwidth]{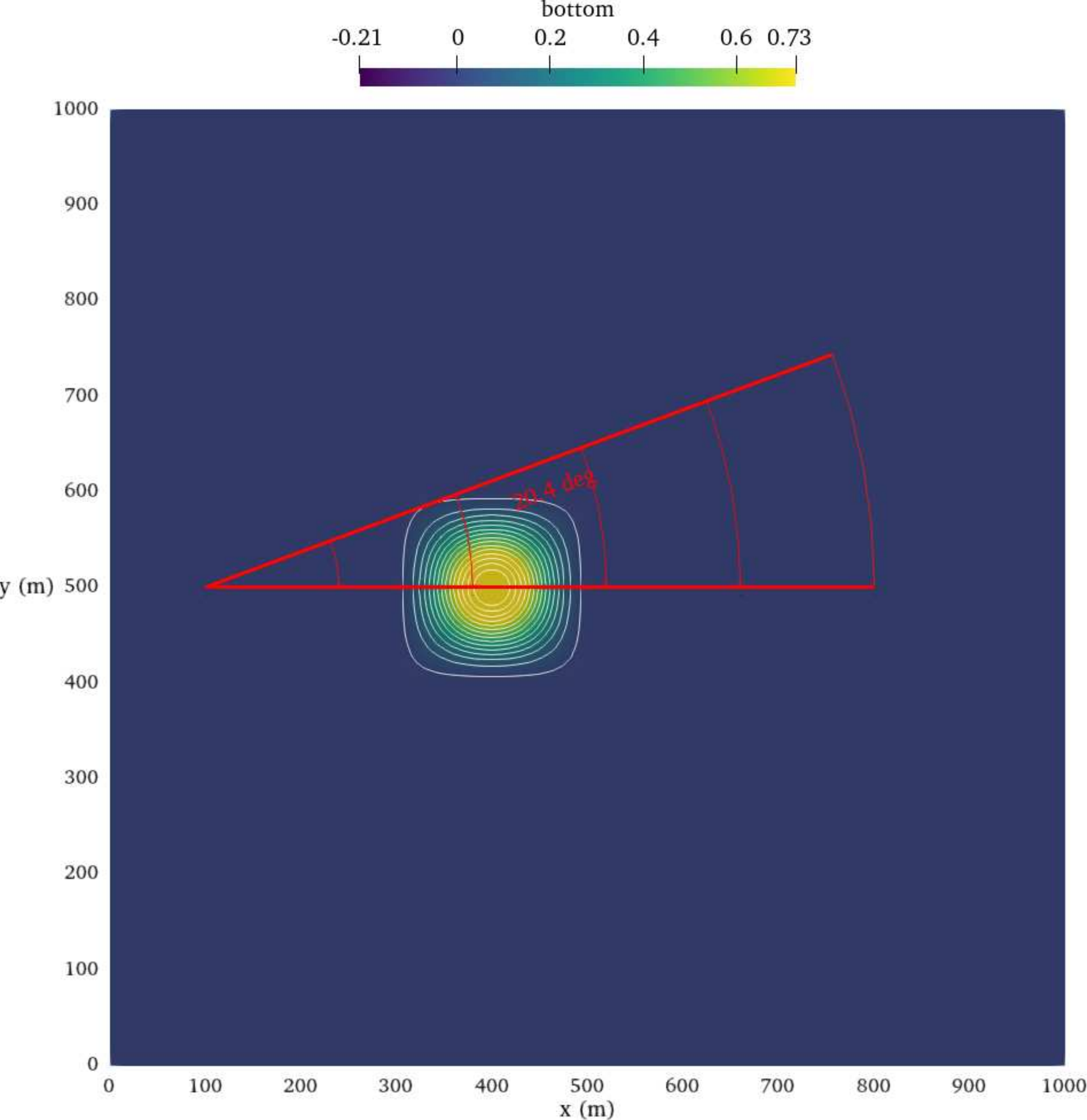}
    \caption{$t=\SI{0}{s}$}
  \end{subfigure}%
  \hspace*{\fill}
  \begin{subfigure}{0.495\linewidth}
    \includegraphics[width=\textwidth]{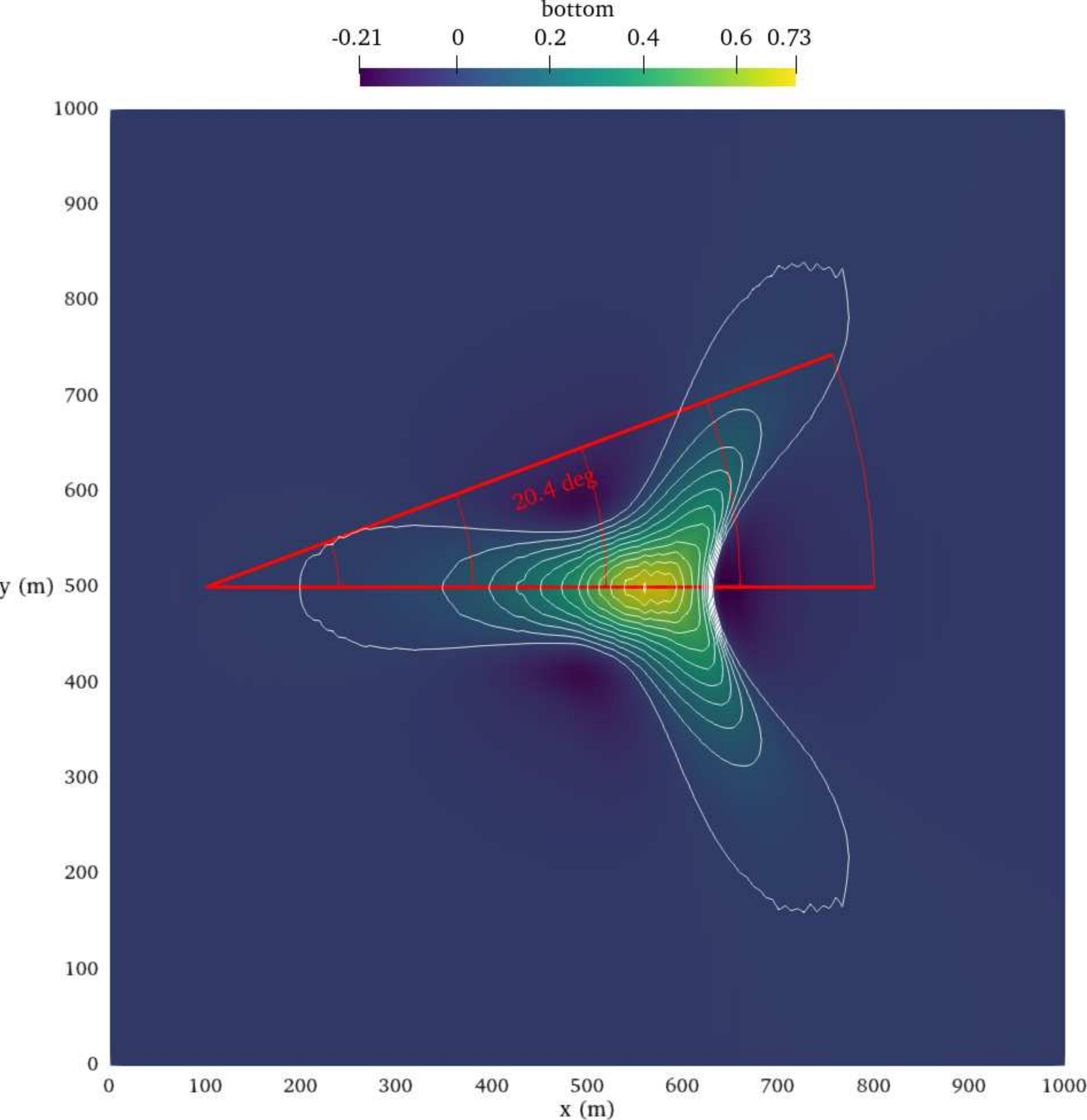}
    \caption{$t=\SI{360000}{s}$}
  \end{subfigure}%
  \caption{Estimation of the spreading angle
           $\alpha_{\textrm{num}} = \SI{20.4}{\degree}$
           for the numerical approximation of a conical sand dune found
           by examining the iso-level $0.0125$ at the initial and final times.
           The theoretical prediction of the spreading angle for this problem
           setup is $\SI{21.787}{\degree}$.}
  \label{fig:swe_angle}
\end{figure}

Next, we present simulation results obtained using \RK[RDPK]{4}{9}[\ESstarp][FSAL], which was found to be the most efficient.
The numerical results from the other two time integration techniques were similar and will not be presented.
Figure~\ref{fig:swe_project_sol} shows the bottom topography at the initial and final times.
From this figure we see that the discretization captures the expected star-shaped pattern of the sediment evolution.
We also estimate the spreading angle of a conical sand dune obtained for the numerical simulation \cite{uh2021unstructured}.
This is done by examining the angle up to the iso-level $0.0125$ over time as shown in Figure~\ref{fig:swe_angle}.
We are able to numerically estimate the spreading angle to be $\alpha_{\textrm{num}} = \SI{20.4}{\degree}$.
The numerical spreading angle of the current method compares well to the theoretical angle \eqref{eq:dune_angle}
as well as to similar angle studies done in the literature
\cite{benkhaldoun2010two,diaz2009two,hudson2001numerical,uh2021unstructured}.

We applied error-based time stepping methods to the shallow water equations augmented
with an Exner equation to model sediment transport.
This served to demonstrate that it is convenient and straightforward to extend existing code for a new system of equations
and obtain a first set of numerical results.
This was particularly useful for the considered shallow water Exner equations because eigenvalue estimates for a general problem
setup are unwieldy to obtain \cite{chertock2020operator,diaz2009two}.
These preliminary results for the evolution and spreading of a conical sand dune compared well to the theoretical solution behavior
as well as results from the literature.
The development of high-order discontinuous Galerkin methods for sediment transport problems is the subject of ongoing research.
Details about the physical problem setup and numerical parameters as well as all code
necessary to reproduce the figures in this section are available in the reproducibility
repository for this article \cite{ranocha2022errorRepro}.

\section{Summary and conclusions}
\label{sec:summary}

We have compared CFL- and error-based step size control of explicit Runge-Kutta
methods applied to discontinuous Galerkin semidiscretizations of compressible
CFD problems. Our numerical experiments demonstrate that error-based step size
control is convenient and robust in a wide range of applications, including
industrially relevant complex geometries with curved meshes in multiple dimensions, shock capturing
schemes, initial transient periods (cold start of a simulation), semidiscretizations
involving a change of variables, and for exploratory research involving new systems
of equations.

There are cures for some problems of CFL-based step size control
discussed in this article. For example, a cold start requiring an initially smaller
CFL factor could run a few hundred time steps with a small CFL factor and use
a larger CFL factor after a restart --- or increase the CFL factor based on a
user-supplied function in the first time steps. A change of linear stability
behavior of shock capturing schemes could be handled by including the shock
capturing indicator in CFL-based step size control. However, this requires sophisticated
additional mechanisms as implemented in some in-house codes.
On the other hand, error-based step size control worked right out of the box
without further modifications for the problems considered here. It does not
require special handling of many edge cases and comes with a reduced sensitivity
with respect to user-supplied parameters.

While we have shown that error-based step size control is convenient and robust
in many applications, including industrially relevant large-scale computations
of compressible CFD problems, it does not come with hard mathematical guarantees.
For example, some invariant domain preserving methods come with some robustness
properties guaranteed under certain CFL constraints; error-based step size control
alone can usually not be proven to satisfy these constraints in all cases.
Nevertheless, based on our experience, it is usually robust for suitably chosen
shock capturing mechanisms.

As a byproduct of this research, we also compared several time integration
methods for hyperbolic conservation laws and compressible flow problems.
For lower Mach number flows or problems with positivity issues and/or strong
shock capturing requirements, \ssp43 is usually the most efficient scheme.
Otherwise, \RK[RDPK]{3}{5}[\ESstarp][FSAL] and \RK[RDPK]{4}{9}[\ESstarp][FSAL]
are usually the most efficient methods, where \RK[RDPK]{3}{5}[\ESstarp][FSAL]
tends to be more robust and \RK[RDPK]{4}{9}[\ESstarp][FSAL] tends to be more
efficient for stricter tolerances. Nevertheless, \RK[BS]{3}{3}[][FSAL] is
still one of the best general-purpose ODE solvers for problems like the ones
considered in this article.

Future research includes the investigation of error-based step size control
with adaptive mesh refinement and stabilization techniques such as
positivity preserving limiters \cite{zhang2011maximum}. These are usually
implemented as callbacks outside of the inner time integration loop and are
thus ``not seen'' by error estimators. Moreover, we would like to optimize
new RK methods including error estimators and step size controllers for
compressible CFD problems in other ranges than \cite{jahdali2021optimized},
e.g., focusing on lower-Mach flows and dissipative effects.

\section*{Acknowledgments}

Andrew Winters was funded through Vetenskapsr{\aa}det, Sweden grant
agreement 2020-03642 VR.
Some computations were enabled by resources provided by the Swedish National
Infrastructure for Computing (SNIC) at Tetralith, partially funded by the Swedish
Research Council through grant agreement no. 2018-05973.
Hugo Guillermo Castro was funded through the award P2021-0004 of
King Abdullah University of Science and Technology.
Some of the simulations were enabled by the Supercomputing Laboratory and the Extreme
Computing Research Center at King Abdullah University of Science and Technology.
Gregor Gassner acknowledges funding through the Klaus-Tschira Stiftung via the project ``HiFiLab''.
Gregor Gassner and Michael Schlottke-Lakemper acknowledge
funding from the Deutsche Forschungsgemeinschaft through the research unit ``SNuBIC'' (DFG-FOR5409).

\printbibliography

\end{document}